\documentclass[12pt,reqno]{amsart}
\usepackage{amsmath,amssymb,amscd,latexsym,amsthm,mathrsfs,comment, amsfonts, euscript}
\usepackage[unicode]{hyperref}
\usepackage{hypbmsec}
\newtheorem{theorem}{Theorem}
\newtheorem{conjecture}{Conjecture}

\theoremstyle{remark}

\newcommand{\F}{\hbox{$\mathbb F$}}
\newcommand{\N}{\hbox{$\mathbb N$}}
\newcommand{\Q}{\hbox{$\mathbb Q$}}
\newcommand{\Z}{\hbox{$\mathbb Z$}}

\renewcommand{\O}{{\EuScript{O}}}
\newcommand{\Nm}{{\mathrm{Nm}}}
\newcommand{\mmod}{{\rm mod~}}
\newcommand{\backddots}{\mathinner{\mkern1mu\raise1pt\vbox
{\kern7pt\hbox{.}}\mkern2mu
\raise4pt\hbox{.}\mkern2mu\raise7pt\hbox{.}\mkern1mu}}

\begin{document}
\begin{center}
{\large\bf ARTIN'S PRIMITIVE ROOT CONJECTURE}\\
{\large\bf - a survey -}
\end{center}\vskip .5cm

\begin{center} PIETER MOREE\\ (with contributions by A.C. Cojocaru, W. Gajda and H. Graves) \end{center}
\begin{center}
To the memory of John L. Selfridge (1927-2010)
\end{center}
\vskip .5cm
\begin{minipage}[c]{12cm}
{\small {\sc Abstract}. One of the first concepts one meets in elementary number
theory is that of the multiplicative order. We give a survey of the literature
on this topic emphasizing the Artin primitive root conjecture (1927). The first part
of the survey is intended for a rather general audience and
rather colloquial, whereas the second part is intended for number theorists and ends
with several open problems. The contributions in the survey on `elliptic Artin' are due to 
Alina Cojocaru.
Wojciec Gajda wrote a section on `Artin for K-theory of number fields', and Hester Graves (together with me) on 
`Artin's conjecture and Euclidean domains'.}
\end{minipage}\vskip .5cm
\tableofcontents
\section{Introduction}
Gauss \cite{Gauss} considered
in articles 315-317 of his Disquisitiones Arithmeticae (1801)
the decimal expansion of numbers of the form $1\over p$ with $p$
prime.  For example, ${1\over 7}=0.{\underline{142857}}142857\ldots$, ${1\over
11}=0.{\underline{09}}0909\ldots$.  

These decimal expansions for $p\ne 2$ and $p\ne 5$ are purely periodic.  It
is easy to see that the period is equal to the smallest positive integer $k$ such
that $10^k\equiv 1({\rm mod~}p)$.  This integer $k$ is the {\it multiplicative
order} of 10 modulo $p$ and is denoted by ${\rm ord}_p(10)$.   The integer $k$
equals the order of the subgroup generated by $10$ in $(\Z/p\Z)^*$,
the multiplicative subgroup of residue classes modulo $p$.  By Lagrange's
theorem the order of a subgroup of a finite group is a divisor of the number
of elements in the group.  Since $(\Z/p\Z)^*$ has
$p-1$ elements, we conclude that ${\rm ord}_p(10)|p-1$.  If 
${\rm ord}_p(10)=p-1$,
we say that $10$ is a {\it primitive root} mod $p$.  The decimal expansion we have
given for $1\over 7$ shows that $10$ is a primitive root mod $7$.  Gauss'
table gives several other such examples and he must have asked himself
whether there are infinitely many such primes, that is, primes $p$ for which the
decimal period is $p-1$.  

Some light on this question can be shed on using the simple observation 
(due to Chebyshev) that if $q=4p+1$ is a prime with $p$ a prime
satisfying $p\equiv 2({\rm mod~}5)$, then $10$ is a primitive root modulo $q$. 
Note that, a priori, ${\rm ord}_q(10)\in \{1,2,4,p,2p,4p\}$.  Now,
using the law of quadratic reciprocity, see, e.g., Lemmermeyer \cite{Lemmermeyer}, one deduces that
$$10^{2p}\equiv\left(10\over q\right)=\left(2\over q\right)
\left(5\over q\right)=(-1)^{q^2-1\over 8}\left(q\over 5\right)
=-\left(4\over 5\right) \equiv -1({\rm mod~}q),$$ 
where we also used that $q\equiv 4({\rm mod~}5)$ and $q\equiv 5({\rm mod~}8)$.
Since no prime divisor of $10^4-1$
satisfies the requirements (and
hence ord$_q(10)\nmid 4$), one sees that ${\rm ord}_q(10)=4p=q-1$.  

The {\it logarithmic integral} ${\rm Li}(x)$ is defined
as $\int_2^x{dt/\log t}$. By partial integration one has that 
${\rm Li}(x)\sim x/\log x$ (i.e. the quotient of the r.h.s. and l.h.s. tends 
to 1 as $x$ tends to infinity). 
The prime number theorem in the form 
$$\pi(x):=\# \{p\le x\}\sim {\rm Li}(x),~x\rightarrow \infty,$$
suggests that the `probability that
a number $n$ is prime' is $1/\log n$.  
(Then we expect $\sum_{2\le n\le x}1/\log n$ primes $\le x$ and asymptotically this is
equal to Li$(x)$.)
Thus for both $n$ and $4n+1$
to be prime and $n\equiv 2({\rm mod~}5)$ we expect a probability of 
$1/({5\log^2n})$, assuming independence of these three events.  Since
they are not, we have to correct by some positive constant $c$ and hence expect
that up to $x$ there are at least $cx/{\log^2x}$ 
(note that $\sum_{2\le n\le x}\log^{-2}n\sim x\log^{-2}x$) primes $p$ such that
$10$ is a primitive root mod $p$. Hence we expect that there are
infinitely many primes $p$ having $10$ as a primitive root mod $p$. This 
conjecture is commonly attributed to Gauss, however, to the author's knowledge there
is no written evidence for it.

Emil Artin in 1927, led by a partial heuristic argument (sketched in \S 4), made the following
conjecture, where for a given $g\in \Q^*$ we define
$${\mathcal{P}}(g)=\{p:{\rm ord}_p(g)=p-1\}{\rm ~and~}{\mathcal{P}}(g)(x)=\#\{p\in{\mathcal{P}}(g):p\le x\}.$$  (In general, if $S$
is a set, we define $S(x)=\#\{s\in S:s\le x\}$.)
\begin{conjecture} {\tt Artin's primitive root conjecture (1927)}.\\
Let $g\in \Q\backslash\{-1,0,1\}$.\\
{\it (Qualitative form)  The set ${\mathcal{P}}(g)$ is infinite 
if $g$ is not a square of a rational number.\\
(Quantitative form) Let $h$ be the largest integer such that
$g=g_0^h$ with $g_0\in\Q$. We have, as $x$ tends to infinity,   
\begin{equation}\label{P(g)(x)}
{\mathcal{P}}(g)(x)=\displaystyle\prod_{q\nmid
h}\left(1-{1\over{q(q-1)}}\right)
\displaystyle\prod_{q|h}\left(1-{1\over{q-1}}\right){x\over{\log x}}
+o\left(x\over{\log x}\right).
\end{equation}}
\end{conjecture}
Here an {\it Euler product} $\prod_q f(q)$ appears, where
$f(q)=1+O(q^{-2})$ (to ensure convergence) and $q$ runs over the primes. We will come across many more Euler products
in this paper. Usually they have an interpretation as a density.\\
\indent Write the main term in (\ref{P(g)(x)}) as $A(h)x/\log x$. In case $h$ is even, 
$A(h)=0$ and, furthermore, clearly ${\mathcal{P}}(g)$ is
finite and hence assertion (\ref{P(g)(x)})  is trivial. Thus the
condition that $g$ is not a square is necessary (and according to
Artin sufficient) for ${\mathcal{P}}(g)$ to be infinite.
For $h$ odd we have that $A(h)$ equals a
positive rational multiple of
$$A(1)=A=\displaystyle\prod_q\left(1-{1\over{q(q-1)}}\right)=0.3739558136192\ldots,$$ the {\it Artin constant}. Thus the quantitative
form implies the qualitative form.\\
\indent The product defining the Artin constant does
not converge quickly.  It turns out that for numerical approximations
it is possible and indeed much more efficient to write
$A=\prod_{k=2}^\infty\zeta(k)^{-e_k}$, with $e_k\in\N$ 
and $\zeta(s)=\sum_{n=1}^{\infty}n^{-s}$ the celebrated 
Riemann zeta function. This then allows one to determine $A$ without too much effort with
a precision of $1000$ decimals, as zeta values can be easily evaluated with high 
numerical precision. This method applies to all 
Euler products that occur in the sequel, cf. \cite{Bachcomplex, Cohen, Finch, Monumi}.\\ 
\indent Usually one
speaks about the Artin primitive root conjecture, rather than Artin's
conjecture since there are various unresolved conjectures due to Artin
(most notably the Artin holomorphy conjecture). Indeed, there are even papers where
both these conjectures make an appearance, e.g. \cite{MP}.\\

The starting point of our analysis of Artin's primitive root
conjecture is the following observation : 
\begin{equation}
\label{uitgangspunt} 
p\in{\mathcal{P}}(g)\iff g^{p-1\over q}\not\equiv 1(\mmod p) 
{\rm ~for~every~prime~}q{\rm~dividing~}p-1.
\end{equation}
\noindent "$\Longrightarrow$" Obvious.\\
"$\Longleftarrow$" Suppose $p\notin{\mathcal{P}}(g)$.  Then $g^{p-1\over
k}\equiv 1(\mmod p)$ for some $k|p-1$, $k>1$.  But this implies that
$g^{p-1\over q_1}\equiv 1(\mmod p)$ for some 
prime divisor $q_1$ of $k$. 
This is a contradiction.\\
\indent  Thus, associated with a prime $p$ we have
conditions for various $q$.  We are going to interchange the role of $p$ and
$q$, that is for a {\it fixed} $q$ we consider the set of primes $p$ such that 
$p\equiv 1(\mmod q)$
and $g^{p-1\over q}\not\equiv 1(\mmod p)$.\\
\indent This set of primes can be considered heuristically (\S 2), but also studied
on invoking (analytic) algebraic number theory and hence we make a brief excursion
to this area (\S 3) before taking up our Artinian considerations in \S 4 again.\\ 

\noindent {\tt Remark 1}. Instead of asking whether infinitely often the period length is
{\it maximal}, one might ask the easier question of whether infinitely often the period length is
{\it even}. Here the answer is yes and much more is known. We come back to this in \S 9.2.\\ 

\noindent {\tt Remark 2}. There is more to decimal expansions
than meets the eye. 
For example, Girstmair \cite{Gi} proved that if a prime $p$ is such that its decimal expansion 
$1/p=0.x_1x_2\ldots x_{p-1}x_1\ldots$ has period $p-1$, then the difference between the sums of its even and odd places, 
namely $(x_2+x_4+\ldots +x_{p-1})-(x_1+x_3+\ldots+x_{p-2})$, is $11h_{\mathbb Q(\sqrt{-p})}$ if 
$p \equiv 3({\rm mod~}4)$ and is zero otherwise, where $h_{\mathbb Q(\sqrt{-p})}$ denotes the so-called class number
of the imaginary quadratic number field $\mathbb Q(\sqrt{-p})$ (number fields are briefly discussed in \S 3.) 
For a variation of this result involving $\sum_{k=1}^{p-1}(-1)^k x_k x_{k+r}$ see Aguirre and Peral \cite{AP}. 
Hirabayashi \cite{hira} generalized Girstmair's formula to all imaginary abelian number fields.
Ram Murty and Thangadurai \cite{MThanga} generalized Girstmair's result to the case where the period is
$(p-1)/r$, where now $r>1$ is also allowed. Here generalized Bernoulli numbers $B_{1,\chi}$ enter the picture.
Other examples are 
provided in \S 9.2 and Khare \cite{Khare}.\\

\noindent {\tt Historical remark}. To the modern number theorist it is a bit strange that Gauss spent
such an effort on the rather recreational topic of decimal periods. However, Bullynck \cite{maarten} points
out that this was a German research topic (1765-1801), and that the young Gauss, being aware of these
developments, placed the whole theory on a firm number-theoretic foundation, thereby solving most of the
problems left by the mathematicians before him. One might have expected him to raise the question of
whether the period is maximal for infinitely many primes. There does not seem to be evidence for this.
As Bullynck (personal communication) points out, this kind of existential question is not typical
for 18th century mathematics.\\ 
\indent In Hasse's mathematical diary (Tagebuch) of the year 1927 
there is (information provided by Roquette) an entry entitled ``Die Dichte der 
Primzahlen $p$, f\"ur die $a$ primitive Wurzel ist (nach 
m\"undlicher Mitteilung von 
Artin 13. IX. 27)" under the date of 27 Sep. 1927. Although the Artin conjecture presumably predates
the 13th of September of 1927, for lack of written evidence of this, the author takes the 13th of September of 1927
to be its birthday.\\
\indent The Artin conjecture has already been formulated before Artin.
It can already be found in a paper of Cunningham \cite{cunningham}, but, in contrast to Artin, Cunningham's
insights into the problem do not seem to go beyond numerical observations. Today Cunningham is best known
for the Cunningham Project which seeks to factor the numbers $b^n\pm 1$ for
$b = 2, 3, 5, 6, 7, 10, 11, 12$ up to high powers $n$, see, e.g., \cite{Wagstaff2}.\\

\noindent {\tt Warning}. This survey aims to be fairly complete, but still it 
reflects the (un)familiarity of its author with various aspects of the topic. He sought
to compensate for this in some cases by inviting a `guest surveyor' to write something on aspects
not so familiar to him, where there is a more extensive literature available.

\section{Naive heuristic approach}

Fix any prime $q$. We try to compute the density of primes $p$ such that 
both $p\equiv 1(\mmod q)$ and
$g^{p-1\over q}\equiv 1(\mmod p)$. 
The {\it prime number theorem for arithmetic progressions} states that, as $x$ tends 
to infinity,
\begin{equation}
\label{pntap}
\displaystyle\pi(x;d,a):=\sum_{p\leqq x\atop p\equiv a(\mmod d)}1
\displaystyle\sim {x\over \varphi(d)\log x}.
\end{equation}
The residue classes $a({\rm mod~}d)$ with $(a,d)=1$ are said to be {\it primitive}.
Note that given an integer $d$, with finitely many exceptions, a prime must be in a primitive residue class modulo $d$.
Since there are $\varphi(d)$ primitive residue classes modulo $d$, (\ref{pntap}) says
that the primes are asymptotically equidistributed over the primitive residue
classes modulo $d$ (in view of the prime number theorem). Dirichlet's theorem (1837) is the weaker statement that any primitive
residue class contains infinitely many primes.\\
\indent By (\ref{pntap}) $p\equiv 1(\mmod q)$ holds true for primes $p$
with frequency $1/\varphi(q)={1/(q-1)}$.  
Recall Fermat's little theorem
which asserts that $g^{p-1}\equiv 1({\rm mod~}p)$ if $p\nmid g$.
Using this we infer, in case $p\nmid g$, that $g^{p-1\over
q}$ is a solution of $x^q\equiv 1(\mmod p)$.  We expect there to be $q$
solutions and we want a solution to be $1$ modulo $q$.  Thus we expect to be
successful with probability $1\over q$, except when $q|h$.  Then
$g^{p-1\over q}=g_0^{h{p-1\over q}}\equiv 1(\mmod p)$, trivially. If
we assume that these events are independent, then the
probability that both events occur is $1\over q(q-1)$ 
if $q\nmid h$ and $1\over q-1$ otherwise.

By (\ref{uitgangspunt}) the above events should not occur for any
$q$ in order to ensure that $p\in{\mathcal{P}}(g)$.
This suggests a natural density of 
$$\displaystyle\prod_{q\nmid h}\left(1-{1\over q(q-1)}\right)
\displaystyle\prod_{q|h}\left(1-{1\over q-1}\right)=A(h)$$ for such primes and
hence we expect (\ref{P(g)(x)}) to hold true.

\section{Algebraic number theory}

In order to understand Artin's approach to his conjecture we need some facts
about algebraic number theory.  We say $\alpha$ is {\it algebraic} over
$\Q$ if it satisfies an equation of the form $f(\alpha)=0$ with
$f(x)\in\Z[x]$.  We say $\alpha$ is an {\it algebraic integer} over
$\Q$ if it satisfies an equation of  the form $f(\alpha)=0$, where
$f(x)\in\Z[x]$ is monic.  A {\it number field} over $\Q$ is a
field obtained by adjoining finitely many algebraic numbers
$\alpha_1,\ldots,\alpha_s$ to $\Q$.  That is, $K=\Q(\alpha_1,\ldots,\alpha_s)$.  By $\Q(\alpha_1,\ldots,\alpha_s)$ we
denote the smallest field containing 
$\alpha_1,\ldots,\alpha_s$ and $\Q$. The
so called {\it theorem of the primitive element} asserts that for
a given number field $K$ there exists an algebraic integer $\alpha$ such
that $K=\Q(\alpha)$.  The element $\alpha$ is said to be a {\it primitive
element}.  For an algebraic integer $\alpha$ let $f_\alpha(x)$ be the
unique (as it turns out) monic polynomial $f_\alpha(x)$ of smallest degree
such that $f_\alpha(\alpha)=0$.
Then $K=\Q(\alpha)$ is isomorphic to $\Q[x]/(f_\alpha(x))$. 
The {\it degree} $[K:\Q]$ of $K$ is defined as deg$f_\alpha(x)$. The
{\it compositum} of two fields $K$ and $L$ is defined as the smallest
field containing both $K$ and $L$.\\

\noindent {\tt Example}. Let  $\alpha=i$.  Then $f_i(x)=x^2+1$. Note that  $\Q[x]$
contains elements of the form $\sum_{j=0}^ma_jx^j$. 
These elements correspond with an element in $\Q[x]/(x^2+1)$
obtained by replacing every $x^2$ by $-1$.  Thus, as a set,  
$\Q[x]/(x^2+1)\cong\{a+bx\mid a,b\in\Q\}$.  Similarly, as a
set, $\Q(i)=\{a+bi\mid a,b\in\Q\}$.  It is not difficult to
see that a field isomorphism between $\Q(i)$ and $\Q[x]/(x^2+1)$ is given by sending $x$ to $i$.\\
\indent For a general number field $K$, we have the following 
picture :
$$\begin{array}{rll}
\Q(\alpha)=K&\supset&{\mathcal{O}}_K\\
\cup&&\cup\\
\Q&\supset&\Z
\end{array}$$ 
Here ${\mathcal{O}}_K$ is the {\it ring of integers} of the field $K$.  It
plays a role analogous to that of $\Z$ in $\Q$ (in our example
 ${\mathcal{O}}_K$ is $\Z[i]=\{a+bi|a,b\in \Z\}$ the {\it ring of
Gaussian integers}). 
The {\it main
theorem of arithmetic} states that, up to the order of factors, factorisation
into primes in $\Z$ is unique.  One might hope for a ring inside $K$
with similar properties.  Let us consider the set of algebraic integers over
$\Q$ that are inside $K$.  This set is actually a ring, the ring of
integers, ${\mathcal{O}}_K$, of $K$.  It usually does not have unique factorization in
elements.  It, however, has unique factorization into terms of {\it prime
ideals} (up to order of factors).  The prime ideal $(p)$ in $\Z$ turns out to factor as
${\mathfrak P}_1^{e_1}\ldots {\mathfrak P}_g^{e_g}$ inside ${\mathcal{O}}_K$.  The equivalent statement
of $p=\# \Z/p\Z$ in $\Z$, reads in  ${\mathcal{O}}_K$:
$p^{f_i}=\# {\mathcal{O}}_K/{\mathfrak P}_i{\mathcal{O}}_K$.  Here $f_i$ is called
the {\it degree} of the prime ideal ${\mathfrak P}_i$.  
More generally, 
$\# {\mathcal{O}}_K/{\mathfrak{a}}{\mathcal{O}}_K$, notation $N\mathfrak{a}$ is 
the {\it norm} of the ideal $\mathfrak{a}$.
We have the relationship
$\sum_{i=1}^ge_if_i=n$, with $n$ the degree of $K$.  
We say that a rational prime $p$ splits
completely in ${\mathcal{O}}_K$ if $e_i=f_i=1$ for $i=1,\ldots,g$.  Note that
in this case $g=n$.  If a prime $p$ splits completely in the ring of
integers ${\mathcal{O}}_K$, with $K\cong\Q[x]/(f(x))$, then over
$\Z/p\Z$, $f(x)$ splits into $n$ distinct linear factors.

\subsection{Analytic algebraic number theory}
Let $\mathfrak{O}$ run through the non-zero integral ideals of ${\mathcal{O}}_K$,
with $K$ a number field.  For $\Re(s)>1$, we define
$\zeta_K(s)=\sum_{\mathfrak{O}}N\mathfrak{O}^{-s}$ and
elsewhere by analytic continuation.  
Note that if $K=\Q$, then $\zeta_K(s)=\zeta(s)$, the usual Riemann zeta-function.
For $\Re(s)>1$ we have an Euler
product, 
$\zeta_K(s)=\prod_{\mathfrak P}(1-
N{\mathfrak P}^{-s})^{-1}$, where
the product runs over all prime ideals $\mathfrak P$ of  
${\mathcal{O}}_K$. 
In 1903, Landau proved the {\it Prime Ideal Theorem} to the effect that
$$\pi_K(x):
=\sum_{N{\mathfrak P}\leqq
x}1\sim {\rm Li}(x),~~ x\to\infty.$$  Assuming the {\it Riemann
Hypothesis} for the field $K$, that is, assuming that all the non-trivial
zeros of $\zeta_K(s)$ are on $\Re(s)={1\over 2}$, one obtains the much
sharper estimate
\begin{equation}\label{pik(x)}
\pi_K(x)={\rm Li}(x)+O_K(\sqrt{x}\log x),
~~x\to\infty,
\end{equation}
where the implied constant depends at most on $K$. (For Landau's big O-notation and related
notation like $\Omega(x)$ we refer to any introductory book on analytic number theory, e.g., 
\cite{Apostol, mova, tenenbaum}.)\\
\indent  In the
analysis of Artin's primitive root conjecture an infinite family of fields
will play a role and we need an error term in (\ref{pik(x)}) that is
explicit in its dependency on
$K$.  Such a result was folklore and finally written down by Lang 
\cite{Lafolk}:
\begin{equation}
\label{langfolklore} 
\pi_K(x)={\rm Li}(x)+O(\sqrt{x}\log(x^{[K:{\mathbb Q}]}|D_K|)),
\end{equation} 
where $D_K$ 
denotes the
discriminant of the field $K$.\\
\indent The function $\pi_K(x)$ is heavily biased
towards prime ideals of degree one:
\begin{eqnarray*}
\pi_K(x)&=&\sum_{N{\mathfrak P}\leqq x\atop N{\mathfrak P}=p}1
+\sum_{N{\mathfrak P}\leqq x\atop N{\mathfrak P}=p^2}1+\ldots\\
&=&\sum_{N{\mathfrak P}\leqq x\atop N{\mathfrak P}=p}1
+O([K:\mathbb Q]\sum_{p^j\leqq x\atop j\geqq 2}1)
=\sum_{N{\mathfrak P}\leqq x\atop N{\mathfrak P}=p}1+O([K:\mathbb Q]\sqrt{x}).
\end{eqnarray*}

In Artin's problem one is specifically interested in so called {\it normal
fields}.  For such a field all prime ideals of degree one, with at most
finitely many exceptions, come from rational primes $p$ that split
completely.  Throwing out these exceptions we have the following picture :
$$\begin{array}{c}
{\mathfrak P}_1{\mathfrak P}_2\dots{\mathfrak P}_n\\
\ddots\vdots\backddots\\
p_1
\end{array}\qquad
\begin{array}{c}
{\mathfrak P}'_1{\mathfrak P}'_2\dots{\mathfrak P}'_n\\
\ddots\vdots\backddots\\
p'_1
\end{array}$$
(A $[K:\Q]$-fold covering of the rational primes.)
Thus, for a normal field we have 
$$\pi_K(x)=[K:\Q]\sum_{p\leqq x\atop
p{\rm ~splits~completely~in~}K}1+O([K:\mathbb Q]\sqrt{x}\log x),$$ 
that is, by the Prime Ideal Theorem, 
\begin{equation}\label{sum{p leqq x}1}
\begin{array}{cll}
\displaystyle\sum_{p\leqq x\atop p{\rm ~splits~completely~in~}K}1
\sim{1\over[K:\Q]}{x\over\log x},~~x\to\infty.
\end{array}
\end{equation}
This is a particular case of a very important result called the 
{\it Chebotarev
density theorem}.  For a nice introductory account see Lenstra and
Stevenhagen \cite{LS}. Another recommendable paper (research level)
on this topic is a paper by Serre \cite{Serre}.\\ 

\noindent {\tt Example}. (Cyclotomic fields). A {\it
cyclotomic field} $K$ is a field of the form $K=\Q(\zeta_n)$ with
$\zeta_n=e^{2\pi i/n}$ and $n$ a natural number.  One can
show that $$f_{\zeta_n}(x)=\displaystyle\prod_{a=1\atop
(a,n)=1}^n(x-\zeta_n^a)=\Phi_n(x)\in\Z[x].$$  The polynomial
$\Phi_n(x)$ is the $n$th {\it cyclotomic polynomial}.  From this we deduce that
$[\Q(\zeta_n):\Q]={\rm deg~}\Phi_n(x)=\varphi(n)$.  The cyclotomic
fields are normal.  It can be shown that 
$p$ splits completely in $K$ iff
$p\equiv 1(\mmod n)$ iff $x^n-1$ factors completely over $\Z/p\Z$.
 
Applying (\ref{sum{p leqq x}1}) we deduce that 
$$\displaystyle\sum_{p\leqq x\atop p\equiv 1(\mmod n)}1
\displaystyle\sim {x\over \varphi(n)\log x},~~~x\to\infty,$$
a particular case of (\ref{pntap}).  The same result can be deduced for every 
primitive residue class modulo $n$ on invoking Chebotarev's density theorem. Thus
the Chebotarev density theorem implies the prime number theorem for arithmetic progressions
(\ref{pntap}).

\section{Artin's heuristic approach} 
We have now assembled the required preliminaries to continue the story on
the progress made on Artin's primitive root conjecture.  Remember that we
are interested in the set of all primes $p$ satisfying 
$$p\equiv 1(\mmod q),~g^{p-1\over q}\equiv 1(\mmod p),$$
where $q$ is any fixed prime.
By what we have learned about cyclotomic fields the primes $p\equiv
1(\mmod q)$ are precisely those for which the equation $x^q\equiv 1(\mmod p)$
has $q$ distinct solutions mod $p$.  Claim: if $g^{p-1\over
q}\equiv 1(\mmod p)$, then $x^q\equiv g(\mmod p)$ has a solution mod $p$.  To
see this, write $g=\gamma^a$ with $\gamma$ a primitive root mod $p$.  Then
$\gamma^{a{p-1\over q}}\equiv 1(\mmod p)$.  A power of a primitive root can
only be congruent to $1(\mmod p)$ if the exponent is a multiple of $p-1$,
hence $q|a$.  This implies $a=bq$.  Then $(\gamma^b)^q\equiv g(\mmod p)$ and
this proves the claim.  Let $\alpha_1,\ldots,\alpha_q$ be
the distinct mod $p$ solutions of $x^q\equiv 1(\mmod p)$.  Then we see that
$\alpha_1\gamma^b,\ldots,\alpha_q\gamma^b$ are all distinct solutions of
$x^q\equiv g(\mmod p)$. We conclude
that a prime which satisfies $p\equiv 1(\mmod q)$
and $g^{p-1\over q}\equiv 1(\mmod p)$ splits completely in the number field
\begin{equation}
\label{fielddefinition}
k_q:=\Q(\zeta_q,g^{1\over q}).
\end{equation}  
(Note that it does not make a difference
which $q$th root of $g$ we take, since we always end up with the same
field.)  We leave it to the reader to show that if $p$ splits completely in
$k_q$, then $p\equiv 1(\mmod q)$ and $g^{p-1\over q}\equiv 1(\mmod p)$.  Thus
we have 
$$p\equiv 1({\rm mod~}q),~g^{p-1\over q}\equiv 1({\rm mod~}p)\iff p {\rm ~splits~completely~in~}k_q.$$
Let us denote the degree of $k_q$ by $n(q)$.  It is not difficult to show
that $k_q$ is normal.  Hence we may apply (\ref{sum{p leqq x}1}) to deduce
that the number of primes $p\le x$ that do not split completely in $k_q$
equals $(1-1/n(q))x/\log x$, asymptotically.
So one expects $\prod_q(1-1/n(q))$ as the
density of primes $p$ for which $g$ is a primitive root mod $p$.\\
\indent  This
heuristic argument was thought to be plausible until about 1960, when the Lehmers 
\cite{Lehmers} made some
numerical calculations that did not seem to always match with Artin's
heuristic.
Then, in 1968, Heilbronn realized that the events `$p$ does not
split completely in $k_q$' are not 
necessarily independent as $p$ and $q$ range
through all primes and
published a corrected quantitative conjecture. Artin, however, made
this correction much earlier, namely in 1958 in a letter to the Lehmers
in a response to a letter of the Lehmers regarding his numerical work. 
Artin
did not publish his corrected conjecture, nor did the Lehmers refer to
Artin in their paper \cite{Lehmers}, although they give the correction
factor. 
As late as 1964 Hasse provided a correction factor in the 1964 edition
of his book \cite{Hasseboek} that is incorrect if $g\equiv 1({\rm mod~}4)$ is
{\it not} a prime.
For some excerpts of the correspondence between Artin and
the Lehmers,  see Stevenhagen \cite{Steven}.\\
\indent Take for example $g=5$, thus $h=1$ (with $h$ defined as in Conjecture 1).  Then
$k_2=\Q(\sqrt{5})\subseteqq k_5=\Q(\zeta_5,5^{1/5})$, i.e. $k_2$
is a subfield of $k_5$ (since
$\sqrt{5}=\zeta_5-\zeta_5^2-\zeta_5^3+\zeta_5^4$).
Now if $K\subseteqq L$ and $p$ splits completely in
$L$, then $p$ must split completely in $K$.  This means that the condition
`$p$ does not split completely in $k_2$' implies the condition `$p$ does
not split completely in $k_5$'.  So, assuming that
there are no further dependencies
between the various conditions on $q$, we expect that
$$
{\mathcal{P}}(5)(x)\sim \prod_{q\ne 5}\left(1-{1\over n(q)}\right){x\over \log x}\sim
\prod_{q\ne 5}\left(1-{1\over{q(q-1)}}\right)
{x\over \log x },~~x\rightarrow\infty.
$$
Note that the Euler product involved equals $20A/19$.
This turns out to be more consistent with the numerical data than
Artin's prediction $A$.
\section{Modified heuristic approach (\`a la Artin)}
Recall that $k_q$ is defined in (\ref{fielddefinition}) for
prime $q$. Let $m$ be squarefree. We define $k_m$ to be the
compositum of the fields $k_q$ with $q|m$ and prime.  It can
be shown that $k_m=\Q(\zeta_m,g^{1/m})$. We are interested in the density of primes $p$ that do not
split completely in any $k_q$ (if
this exists).
We can try to compute it by inclusion-exclusion.  So
we start writing down $1-{1\over n(2)}-{1\over n(3)}\cdots$.  However, we
have counted double here the primes $p$ that split completely in both $k_2$
and $k_3$.  It can be shown that the primes that split completely in both
$k_n$ and $k_m$ are precisely the primes that split completely in
$k_{{\rm lcm}(n,m)}$.  Thus we have to add ${1\over n(6)}$.  Continuing in
this way we arrive at the heuristic
\begin{equation}                                                                  
{\mathcal{P}}(g)(x)\sim\sum_{k=1}^\infty
{\mu(k)\over[\Q(\zeta_k,g^{1/k}):\Q]}{x\over\log x},~~x\rightarrow\infty.
\end{equation}
Here $\mu$ is the {\it M\"obius function}. It is defined by $\mu(1)=1$,
and if $n>1$, then $\mu(n)=(-1)^s$, where $s$ denotes the number of distinct prime factors
of $n$ if $n$ is squarefree, and $\mu(n)=0$ otherwise.
If we are only interested in the primes $p$ that do not split completely in any of a finite
set of number fields, then there would be little
to prove (by invoking Chebotarev's theorem).  It is the 
{\it infinitely} many $q$ that makes the problem hard.

Artin's heuristic amounts to assuming that the fields $k_{q_1}$ and $k_{q_2}$
are {\it linearly disjoint} over $\Q$ (i.e. 
$k_{q_1}\cap k_{q_2}=\mathbb Q$), for distinct primes $q_1$ and $q_2$.  
This has as a consequence that 
$[k_{q_1q_2}:\Q]=[k_{q_1}:\Q][k_{q_2}:\Q]$ and in
general for squarefree $k$ that $n(k)=k\varphi(k)/(k,h)$.  Thus,
according to Artin we should have 
$$\sum_{k=1}^\infty{\mu(k)\over n(k)}
=\sum_{k=1}^\infty{\mu(k)(k,h)\over k\varphi(k)}
=\prod_{q\nmid h}\left(1-{1\over q(q-1)}\right)
\prod_{q|h}\left(1-{1\over q-1}\right)=A(h),$$ 
where in the derivation of the last equality we used the fact 
(see, e.g., \cite[Theorem 1.9]{mova})
that if
$f(k)$ is a multiplicative function, then
$$\sum_{k=1}^\infty|f(k)|<\infty\Rightarrow
\sum_{k=1}^\infty f(k)=\prod_q(1+f(q)+f(q^2)+\cdots).$$
(A function $f$ on the integers is called multiplicative if when
$(n,m)=1$, then $f(nm)=f(n)f(m)$.)
It can be shown that if $g\ne -1$ and $g$ is not a square, then
$$\displaystyle\sum_{k=1}^\infty{\mu(k)\over n(k)}=c_gA,$$ where $c_g>0$ is a
rational number depending on $g$.\\
\indent Let $g\ne -1$ or a square.  Hooley \cite{Hooley1} proved in 1967
that if the Riemann Hypothesis holds for the number fields
$\Q(\zeta_k,g^{1/k})$ with $k$ squarefree, then we have 
\begin{equation}
\label{classicalhooley} 
{\mathcal{P}}(g)(x)={x\over\log x}\sum_{k=1}^\infty
{\mu(k)\over[\Q(\zeta_k,g^{1/k}):\Q]}
+O_g\left({x\log\log x\over \log^2x}\right)
\end{equation} 
and he explicitly evaluated the latter sum, say $\delta(g)$, as 
\begin{equation}
\label{expli}
\delta(g)=\left\{\begin{array}{rl}
A(h) & {\rm ~if~}d\not\equiv 1({\rm mod~}4);\\  
\left(1-\mu(|d|)\prod_{q|d\atop q|h}{1\over q-2}\prod_{q|d\atop p\nmid h}{1\over q^2-q-1}\right)A(h) & {\rm ~otherwise},
\end{array}\right.
\hbox{}
\end{equation}
with $d$ the discriminant of $\mathbb Q(\sqrt{g})$.
For convenience we denote
the assumption on the fields  $\Q(\zeta_k,g^{1/k})$ made by Hooley as
Hooley's Riemann Hypothesis (HRH) and the quantity
${x\log\log x/\log^2x}$ as $R(x)$.

\section{Hooley's work}

For simplicity we assume $g=2$ and sketch Hooley's proof \cite{Hooley1}.
The restriction $g=2$ allows us to focus exclusively on the analytical
aspects, the algebraic aspects being discussed elsewhere in this survey.\\
\indent A large part of
his paper is devoted to proving an estimate for 
$$\hbox{$P(x,l)=\#\{p\leqq x: p$ splits completely in $k_l\}$},$$
for which the implied constant in the error term does not depend on $l$. 
His result on this is a special case of the result 
(\ref{langfolklore}) of Lang \cite{Lafolk} (proven a few
years later):
\begin{equation}
\label{hohoho}
P(x,l)={{\rm Li}(x)\over [k_l:\mathbb Q]}+O(\sqrt{x}\log(lx)),
\end{equation}
where we assume the Riemann Hypothesis for $k_l$.  Another quantity needed
is 
$$\hbox{$N(x,\eta)=\#\{p\leqq x: p$ does not split completely in
$k_q,~ \forall q\leqq\eta\}$}.$$
By inclusion-exclusion we have $N(x,\eta)=\sum_{l'}\mu(l')P(x,l'),$
where $l'$ ranges over all divisors of $\prod_{q\leqq
\eta}q$. Note that, for $\eta$ large enough, 
$$\displaystyle\prod_{q\leqq\eta}q=e^{\sum_{q\leqq\eta}\log q}\leqq
e^{2\eta},$$
where we used that 
$\sum_{q\leqq\eta}\log q\sim\eta$ (which is equivalent with
the prime number theorem).  
Observe that
${\mathcal{P}}(2)(x)=N(x,x-1)$.  The problem in 
estimating $N(x,\eta)$ is that by using inclusion-exclusion and (\ref{hohoho})
we can only estimate $N(x,\eta)$ for rather small $\eta$, whereas
one would like to take $\eta=x-1$ and so we are forced to work with $N(x,\eta)$ for $\eta$
rather smaller than $x$.  We will actually choose $\eta={1\over 6}\log x$ and 
thus $l'<e^{2\eta}=x^{1/3}$ (for $x$ large enough). Let us
introduce a third quantity $M(x,\eta_1,\eta_2)=\#\{p\leqq x:p$ splits
completely in some $k_q,~\eta_1<q<\eta_2\}.$  It is easy to see that
$N(x,\xi_1)-M(x,\xi_1,\xi_2)-M(x,\xi_2,\xi_3)-M(x,\xi_3,x-1)\leqq
{\mathcal{P}}(2)(x)\leqq N(x,\xi_1)$, that is, 
$${\mathcal{P}}(2)(x)=N(x,\xi_1)+O(M(x,\xi_1,\xi_2))
+O(M(x,\xi_2,\xi_3))+O(M(x,\xi_3,x-1)).$$  
We will choose $\xi_1={1\over 6}\log x$,
$\xi_2=\sqrt{x}\log^{-2}{x}$ and $\xi_3={\sqrt{x}\log{x}}$.  Note
the small gap between $\xi_2$ and $\xi_3$.  This is the `hard region'
and unfortunately it seems out of reach to close it, and work say with 
$${\mathcal{P}}(2)(x)=N(x,\xi_1)+O(M(x,\xi_1,\xi_2))
+O(M(x,\xi_2,x-1)).$$  
Now using the estimate for $P(x,l)$ one easily arrives at   
\begin{eqnarray*}
N(x,\xi_1)&=&{\rm Li}(x)\sum_{l'\leqq x^{1\over 3}}{\mu(l')\over l'\varphi(l')}
+O\left({x\over\log^2{x}}\right)\\
&=&A {\rm Li}(x)+O\left({x\over\log^2{x}}\right)
\end{eqnarray*}
Again using the estimate for $P(x,q)$ we arrive at
$$M(x,\xi_1,\xi_2)\leqq\sum_{\xi_1\leqq q\leqq\xi_2}P(x,q)
=O\left({x\over\log^2{x}}\right).$$

In the region $(\xi_2,\xi_3)$ we use the fact that the primes $p$
counted by $P(x,q)$ certainly satisfy $p\equiv 1(\mmod q)$.  Thus
$P(x,q)\leqq\pi(x,q,1)$.  Then using the Brun-Titchmarsh estimate 
$$\pi(x,q,1)<{2x\over \varphi(q)\log(x/q)},~~~1\leqq q<x,$$
we obtain
$$M(x,\xi_2,\xi_3)\leqq\sum_{\xi_2\leqq q\leqq\xi_3}P(x,q)\leqq 
\sum_{\xi_2\leqq q\leqq\xi_3}\pi(x,q,1)
=O({x\over \log x}\sum_{\xi_2\le q\le \xi_3}{1\over q}).$$
Using a result of Mertens to the effect that
$$\sum_{p\le x}{1\over p}=\log \log x + {\rm constant}+O({1\over \log x}),$$
it then follows that $$\sum_{\xi_2\le p\le \xi_3}{1\over p}=O\Big({\log \log x\over \log x}\Big)$$ and hence
$M(x;\xi_2,\xi_3)=O(R(x))$.
Finally, if a prime $p$ is counted by $M(x,\xi_3,x-1)$, then its order is
less than ${x\over \xi_3}={\sqrt{x}\over \log x}$ and thus it divides
$2^m-1$ for some $m<{\sqrt{x}/\log{x}}$.  Now
$$2^{M(x,\xi_3,x-1)}<
\prod_{p~{\rm counted~by~}M(x,\xi_3,x-1)}p
\leqq\displaystyle\prod_{m<{\sqrt{x}\over\log{x}}}(2^m-1).$$
Thus $$M(x,\xi_3,x-1)<\displaystyle\sum_{m<{\sqrt{x}\over\log{x}}}m
=O\left({x\over\log^2{x}}\right).$$  On gathering all terms we obtain
${\mathcal{P}}(2)(x)=A{\rm Li}(x)+O(R(x)).$\\
\indent Repeating this argument for arbitrary $g\ne -1$ or a square,
we arrive at (\ref{classicalhooley}) under HRH, which implies 
the truth of the
qualitative form of Artin's primitive root conjecture and a
modified form of the quantitative version. This completes the
sketch of the proof.\\

\indent Vinogradov \cite{Vi} has shown that unconditionally
\begin{equation}
\label{wijn}
{\mathcal P}(g)(x)\le \delta(g)\pi(x)+O_g\Big({x(\log \log x)^2\over \log^{5/4}x}\Big),
\end{equation}
where $\delta(g)$ is the density of primitive roots as determined by Hooley (given by (\ref{expli})).
Vinogradov's proof contains some small errors that have been corrected by Van der Waall \cite{Robert}.
Wiertelak \cite{wiedewiet} has established a version of this result 
(with a larger error term, however), where he made the $g$ dependence of the error 
explicit.

\subsection{Unconditional results}
If one asks for unconditional results, the state of affairs is quite
appalling.  There is no known number $g$ for which ${\mathcal{P}}(g)$
is known to be infinite!
In 1986, however, Heath-Brown \cite{HB} proved 
(improving on earlier fundamental work by Gupta, Ram Murty and Srinivasan \cite{GM, MSri}) a 
result which implies 
that there are
at most two primes $q_1$ and $q_2$ for which ${\mathcal{P}}(q_1)$ and
${\mathcal{P}}(q_2)$ are finite and at most three squarefree numbers 
$s_1,s_2$ and $s_3$ for which ${\mathcal{P}}(s_1)$, ${\mathcal{P}}(s_2)$
and ${\mathcal{P}}(s_3)$ are finite.\\
\indent Remember the observation that if $q=4p+1$ and $p\equiv 2(\mmod 5)$, then
$q\in {\mathcal{P}}(10)$.  Sieve theory cannot prove presently that there
are infinitely many such primes. 
However, the so called
lower bound sieve method combined with the Chen-Iwaniec switching method gave
rise to $\gg x\log^{-2}{x}$ primes $p\leqq x$ such that either
$p-1=2^{e}q$ for some prime $q$ or $p-1=2^{e}q_1q_2$ with
$p^\alpha<q_1<p^\delta$ for some $\alpha>{1\over 4}$ and $\delta<{1\over
2}$.  This together with a rather elementary argument, is then enough to
establish Heath-Brown's result.\\
\indent Reznikov and Moree \cite{RM} established a variant of Heath-Brown's result
and used it to make considerable progress regarding a conjecture of Lubotzky and
Shalev to the extent that for all $d$ sufficiently large the number of subgroups
of index $d$ in the fundamental group $\pi_1(M)$ of a hyperbolic three manifold $M$ is at least 
$e^{dc(M)}$, for some positive constant $c(M)$ depending at most on $M$.\\
\indent Skolem has raised the question of whether all the integral solutions of
$X_1X_2-X_3X_4=1$ can be obtained from a fixed polynomial solution by letting
the variables run through $\mathbb Z$ and expressed his belief in favour of a negative
answer. Zannier \cite{Zannier} has shown that if one not only considers solutions in integers, but
those in $S$-integers, where $S$ is a finite set of places, then for a suitable
finite $S$ the answer to Skolem's question is positive. His proof makes use of a fundamental
lemma, the proof of which makes use of a variation of Heath-Brown's arguments.

\subsubsection{Variants for quadratic fields} Narkiewicz \cite{nar} proved the unconditional result, valid for a certain
subclass of abelian number fields, that
if $a_1$, $a_2$, $a_3$ are multiplicatively independent integers in ${\O}_K$ 
with their norms satisfying certain conditions, then
at least one of the numbers $a_1$, $a_2$, $a_3$ is a primitive root for 
infinitely many prime ideals of $K$ of the first degree. The proof uses
ideas from the papers by Gupta and Ram Murty \cite{GM} and Heath-Brown  \cite{HB}.\\  
\indent Given a nontrivial unit in a real quadratic
number field, for a rational prime $p$ which is {\it inert} in the field the maximal order of the unit 
modulo $p$ is $p+1$. An extension of Artin's conjecture is that there are infinitely many such inert primes for which 
this order is maximal. This is known at present only under GRH.
Unconditionally, J. Cohen \cite{JC1} showed that for any choice of 7 units in different real quadratic fields 
satisfying a certain simple restriction, at least one  of these units satisfies 
this version of Artin's conjecture.\\
\indent Given an algebraic number from a quadratic field, for a rational prime $p$ which is inert in the field the maximal order of 
  the unit modulo $p$ is $p^2-1$. An extension of Artin's conjecture is that there are infinitely many 
  such inert primes for which this order is maximal. 
J. Cohen \cite{JC2} recently showed that, given
a quadratic field $K$, for any choice of 85 algebraic numbers satisfying a certain simple restriction, at least 
one of them has order modulo $p$ at least $(p^2-1)/24$ for infinitely many
inert primes $p$ in $K$.

\section{Probabilistic model}
Is there a probabilistic model for a prime
$p$ to be such that $g$ is a primitive root mod $p$? That is, does there exist
a function $f_g(p)$ such that
\begin{equation}
\label{model} 
\sum_{p\in S,~p\leqq x}f_g(p)\sim
\sum_{p\in S,~p\leqq x\atop\hbox{
\scriptsize $g$ is a primitive root mod $p$}}1,
\end{equation} 
where $S$ must be a set of primes having a positive natural density.\\
\indent  There are $\varphi(p-1)$
primitive roots mod $p$ and thus
one could try to take $f_g(p)={\varphi(p-1)/(p-1)}$, assuming
that the primitive roots are
randomly distributed over
$1,\ldots,p-1$. 
Results of Elliott \cite{Elliott}, Cobeli and Zaharescu \cite{CZ-1}, Rudnick and Zaharescu \cite{RZ} and
Wang and Bauer \cite{WB2} indeed all suggest that the distribution
of the primitive roots over $1,\ldots,p-1$ is to a large extent governed by 
Poisson processes. Thus we would hope
to find something like 
\begin{equation}\label{sum(phi(p-1)over(p-1))}
\sum_{p\leqq x}{\varphi(p-1)\over p-1}\sim{\mathcal{P}}(g)(x),~~
x\rightarrow\infty.
\end{equation}
However, by Hooley's theorem there are under HRH $g_1$ and $g_2$ such that
${\mathcal{P}}(g_1)(x)\nsim{\mathcal{P}}(g_2)(x)$.  Thus
(\ref{sum(phi(p-1)over(p-1))}) cannot be true for all $g$.  Let us not be deterred by
this and try to evaluate $\sum_{p\leqq x}{\varphi(p-1)\over
p-1}$.  Note that ${\varphi(n)\over n}=\prod_{p|n}(1-{1\over p})=\sum_{d|n}{\mu(d)\over d}.$ Thus 
$$
\sum_{p\leqq x}{\varphi(p-1)\over p-1}
=\sum_{p\leqq x}\sum_{d|p-1}{\mu(d)\over d}
=\sum_{d\leqq x}{\mu(d)\over d}\sum_{p\le x\atop p\equiv 1(d)}1
=\sum_{d\leqq x}{\mu(d)\over d}\pi(x,d,1).$$
Let $c_1$ be a fixed real number. Then the estimate
\begin{equation}
\label{siwa}
\pi(x;d,a)={{\rm Li}(x)\over \varphi(d)}
+O(xe^{-c_2\sqrt{\log x}})
\end{equation}
holds uniformly for all integers $a$ and $d$ such that
$(a,d)=1$ and $1\le d\le \log ^{c_1} x$, with $c_2$ some
positive constant. This is a well-known result due to Siegel
and Walfisz and sharpens (\ref{pntap}). On applying it we obtain
$$
\sum_{p\leqq x}{\varphi(p-1)\over p-1}
={\rm Li}(x)\sum_{d\leqq\log^{c_1} x}{\mu(d)\over d\varphi(d)}+
O\left({x\over \log^{c_1}x}\right)=A{\rm Li}(x)+
O\left({x\over \log^{c_1}x}\right).$$
(This result is Lemma 1 in Stephens \cite{Stephens}, an earlier argument along
the same lines is given in Pillai \cite{pi0} who considered $\sum_{p\le x}\varphi(p-1)$.)\\
\indent  Recall that Artin's heuristic answer was not always correct because it
failed to take into account quadratic interactions.  So let us try to
incorporate this into our model.  In order that $g$ is a primitive root mod
$p$ it is necessary that $(g/p)=-1$ and this is a quadratic
condition.  In particular if $g$ has to be a primitive root it must be a
non-square mod $p$.  There are $(p-1)/2$ non-squares and ${2\varphi(p-1)/(p-1)}$ is their density.  
Thus it makes sense to try to
evaluate the sum in the identity below. Assuming HRH, 
by a computation a little more complicated than for the sum with
the condition $(g/p)=-1$ omitted: 
$$2\sum_{p\leqq x,~({g\over p})=-1}{\varphi(p-1)\over p-1}
={\mathcal{P}}(g)(x)+O(R(x)).$$  

In case $h>1$, with $h$ as in Conjecture 1, equality does not always hold.  The effect we have to take
into account there is that if $(p-1,h)>1$, then
$$g^{p-1\over(p-1,h)}\equiv
\left(g_0^{h\over(p-1,h)}\right)^{p-1}\equiv 1(\mmod p),$$ and
thus $g$ cannot be a primitive root mod $p$.

Here is a proposal for $f_g(p)$:
\begin{equation}
\label{modelexplicit}
f_g(p)=\left\{
\begin{array}{cl}
{2{\varphi(p-1)\over p-1}} 
&\hbox{if $({q\over p})=-1$ and $(p-1,h)=1$;}\\
0 &\hbox{otherwise}
\end{array}\right.
\end{equation}
It can be shown assuming HRH (see \cite{MoreeAP}) that
\begin{equation}
\label{mmodel}
\sum_{p\leqq x}f_g(p)={\mathcal{P}}(g)(x)
+O(R(x)),
\end{equation}
and thus if $S$ is the set of all primes, then (\ref{modelexplicit}) gives a probabilistic model that works.\\  
\indent  On noting that the naive heuristic fails the traditional approach was to show that it holds true on
average:
$${1\over N}\sum_{g\leqq N}{\mathcal{P}}(g)(x)
={1\over N}\sum_{p\leqq x}M_p(N),$$  
where $M_p(N)$ is the number of integers $g\leqq N$ that are primitive
roots modulo $p$.  Clearly $M_g(N)=\varphi(p-1)\{{N\over p}+O(1)\}$.  Thus
$${1\over N}\sum_{g\leqq N}{\mathcal{P}}(g)(x)
=\sum_{p\leqq x}{\varphi(p-1)\over p-1}
+O\left({x^2\over N\log x}\right)
+O\left({x\over \log^D x}\right),$$  
provided that $N>x^2\log^{D-1}{x}$.
Less trivial unconditional results in this direction were obtained by 
Goldfeld \cite{Goldfeld} and
Stephens \cite{Stephens}, see also Li \cite{Li2}. (For an analogue for
number fields see Egami \cite{Egami}.)\\
\indent It is relatively easy to compute the first sum in
(\ref{model}) with $S=\{p:p\equiv a({\rm mod~}f)\}$ 
and one obtains a relatively complicated but
completely explicit Euler product.  The question is whether 
the true behaviour of $${\mathcal{P}}(g,f,a)(x):=\#\{p\leqq
x:p\in {\mathcal{P}}(g),~~p\equiv a(\mmod f)\}$$ works with this. 

Lenstra's work \cite{Lenstra}, which introduced Galois theory
into the subject, implies the following result on ${\mathcal{P}}(g,f,a)(x)$. 

\begin{theorem}
\label{lenstrarekenkundigerij}
Let $1\le a\le f,~(a,f)=1$. Let
$\sigma_a$ be the automorphism of $\mathbb Q(\zeta_f)$ determined by
$\sigma_a(\zeta_f)=\zeta_f^a$. Let $c_a(n)$ be $1$ if the restriction
of $\sigma_a$ to the field $\mathbb Q(\zeta_f)\cap \mathbb Q(\zeta_n,g^{1/n})$
is the identity and $c_a(n)=0$ otherwise. Put
$$\delta(a,f,g)=\sum_{n=1}^{\infty}{\mu(n)c_a(n)\over 
[\mathbb Q(\zeta_f,\zeta_n,g^{1/n}):\mathbb Q]}.$$
Then, assuming RH for all 
number fields $\Q(\zeta_f,\zeta_n,g^{1/n})$ with
$n$ squarefree,
$${\mathcal{P}}(g,f,a)(x)=\delta(a,f,g){x\over \log x}+
O_{g,f}(R(x)).$$
\end{theorem}
 Inspired by the relatively easy Euler product coming from quadratic
heuristics, the author set out to explicitly evaluate $\delta(a,f,g)$
and managed to find an Euler product for
it (see the 2008 paper \cite{oudekoe} which appeared earlier as a MPIM-preprint in 1998):
\begin{theorem} 
\label{main1}
Let $g$ be an integer with $|g|>1$. Let $h$ be the largest integer such that $g$ is an $h$-th power
of an integer. Write $g=g_1g_2^2$ with $g_1,~g_2$ integers and $g_1$ squarefree.
Put
$$A(a,f,h)=\prod_{q|(a-1,f)}(1-{1\over q})
\prod_{q\nmid f\atop p|h}(1-{1\over q-1})\prod_{p\nmid f\atop p\nmid h}
\left(1-{1\over q(q-1)}\right)$$
if $(a-1,f,h)=1$ and $A(a,f,h)=0$ otherwise.
Let 
$$\beta={g_1\over (g_1,f)}{\rm ~and~}\gamma_1=
\begin{cases}
(-1)^{\beta-1\over 2}(f,g_1) & {\rm ~if~}\beta{\rm ~is~odd};\cr
1 & {\rm ~otherwise}.
\end{cases}
$$ We have
$$\delta(a,f,g)={A(a,f,h)\over \varphi(f)}\left(1-
({\gamma_1\over a}){\mu(|\beta|)\over 
\prod_{q|\beta,~q|h}(q-2)\prod_{p|\beta,~p\nmid h}(q^2-q-1)}\right)$$
in case $g_1\equiv 1({\rm mod~}4)$ or $g_1\equiv 2({\rm mod~}4)$ and
$8|f$ or $g_1\equiv 3({\rm mod~}4)$ and $4|f$ and
$$\delta(a,f,g)={A(a,f,h)\over \varphi(f)},$$
otherwise. Here $({\cdot\over \cdot})$ denotes the Kronecker symbol.
\end{theorem}  
This formula is the same as the one found with quadratic heuristics! 
As a consequence we have the following result \cite{MoreeAP}, which
generalizes (\ref{mmodel}):
\begin{theorem}
Assume RH for all number fields $\Q(\zeta_f,\zeta_n,g^{1/n})$ with $n$
squarefree. 
Then
$$\sum_{p\leqq x\atop p\equiv a(\mmod f)}f_g(p)
={\mathcal{P}}(g,f,a)(x)+O_{g,f}(R(x)).$$
\end{theorem}
A rather
more insightful derivation of the Euler product
for $\delta(a,f,g)$ is given in a paper by Lenstra, Moree and Stevenhagen
\cite{LMS}, as a consequence of a more explicit approach 
involving quadratic characters in the line
of Lenstra's earlier work. This approach is the most powerful thus far
in finding explicit Euler products for primitive root densities. For a preview of this work
see Stevenhagen \cite{Steven}. The upshot is that the density for a whole
class of Artin type problems can be always written as
$$\delta=(1+\prod_q E_q)\prod_q A_q,$$
where the $E_q$ are (real)-character averages and hence satisfy $-1\le E_p\le 1$. Only finitely
many of them are distinct from 1. The infinite product $\prod_q A_q$ can be regarded as the
`canonical' density of the problem at hand. In this formulation the situations where $\delta=0$
can be relatively easily computed (from the expression of the density
as an infinite sum this can be quite hard, cf. the formula for
$\delta(a,f,g)$ in Theorem \ref{lenstrarekenkundigerij}). For the original Artin problem a very easy computation
yields $\prod_q A_q=A(h)$ and $1+\prod_q E_q=E$ and one finds (\ref{expli}) again, see \cite{Steven}. In
the same paper the method is applied to deal with near-primitive roots (see \S 9.7.3), improving on 
earlier treatments of this problem. In particular, the first satisfactory
treatment of the vanishing in the near-primitive root problem is given. In a sequel to this paper, authored
by Moree and Stevenhagen \cite{Most}, the method is used to deal with some higher-rank variations of Artin's primitive root
conjecture. Again this leads to the quickest known proofs of Euler product forms of the relevant densities.\\
\indent  Rodier \cite{Rodier}, in connection with
a coding theoretical result involving Dickson polynomials, was interested in the set of primes $p$
such that $p\equiv a({\rm mod~}28)$, with $a\in \{-1,3,19\}$. Assuming equidistribution over
the congruence classes modulo 28 one would expect (as
did Rodier) that the density of this set of primes
is $3A/\varphi(28)=A/4$. Moree \cite{morod} computed the density of primes $p$ (on GRH) such that
$2$ is a primitive root modulo $p$, $({p\over l_j})=\epsilon_j$, 
for $1\le j\le s$ and $p\equiv \epsilon_0({\rm mod~}4)$,
with $\epsilon_0,\ldots,\epsilon_s\in \{\pm 1\}$ and the $l_j$ distinct odd primes. Applying
this result with $\epsilon_0=-1$, $s=1$, $\epsilon_1=-1$ and $l_1=7$, it then follows that
Rodier's set has density $21A/82$ on GRH, contradicting his conjecture. Of course one can also
use Theorem \ref{main1} to arrive at this density.\\
\indent One can ask for which $f$ and $g$ we have equidistribution over the congruence
classes modulo $f$ of the primes $p$ for which $g$ is a primitive root modulo $p$, that
is for which $f$ and $g$ does there exist $\delta$ such that $\delta(a,f,g)=\delta$ for
all $a$ with $(a,f)=1$? The explicit formula for $\delta(a,f,g)$ allows one to determine the $f$
such that the primes in ${\mathcal{P}}(g)$ are equidistributed in the
various primitive residue classes mod $f$. This was done earlier by Moree
\cite{Moree3} using an Euler product for $\delta(1,f,g)$ (which
is easier to obtain) and some Galois theoretic arguments. Using the Euler product
for $\delta(a,f,g)$ given in \cite{oudekoe} a shorter proof can be given (in the same 
article).
Basically $f$ must be a power of two in case $h=1$.  Also $2^\alpha3^\beta$
can occur, the smallest positive integer $g$ for which this happens is $g=21^7$.

\section{The indicator function}
\subsection{The indicator function and probabilistic models}
Let $G$ be a finite cyclic group of order $n$ and put $f(g)=1$ if
$G$ is generated by $g$ and $f(g)=0$ otherwise. It was already
observed a long time ago (see, e.g., \cite[Satz 496]{Landauvor} and
\cite{hinz-1} for a version over number fields) that
\begin{equation}
\label{indicator}
f(g)={\varphi(n)\over n}\sum_{d|n}{\mu(d)\over \varphi(d)}
\sum_{{\rm ord}\chi=d}\chi(g),
\end{equation}
where the sum is over the characters $\chi$ of $G$ having order $d$.
Using this, one deduces that
$$
{\mathcal{P}}(g)(x)=\sum_{p\le x}{\varphi(p-1)\over p-1}\sum_{d|p-1}
{\mu(d)\over \varphi(d)}\sum_{{\rm ord}\chi=d}\chi(g),$$
where the inner sum is over the characters of $(\mathbb Z/p\mathbb Z)^*$ of
order $d$.
Now write ${\mathcal{P}}(g)(x)=\sum_{d\le x}S_d(x)$, with
$$S_d(x)={\mu(d)\over \varphi(d)}\sum_{p\le x\atop p\equiv 1({\rm mod~}d)}
{\varphi(p-1)\over p-1}\sum_{{\rm ord}\chi=d}\chi(g).$$
Note that
$$S_1(x)=\sum_{p\le x}{\varphi(p-1)\over p-1}$$
and that
$$S_2(x)=-({g\over p})\sum_{p\le x}{\varphi(p-1)\over p-1},$$
with $({g\over p})$ the Legendre symbol.
The term
$S_d(x)$ can be written as a linear combination of terms
$${1\over \pi(x)}
\sum_{p\le x,~p\equiv 1({\rm mod~}d)\atop p\equiv 1({\rm mod~}\delta)}
\left({g\over p}\right)_d,$$
where $({g\over p})_d$ denotes the $d$-th power residue symbol, cf.
Lemmermeyer \cite[p. 111]{Lemmermeyer}.
In case $g$ is squarefree these terms
are $o(S_1(x))$ for every $d>2$ and $\delta$.
Hence we expect that
$${\mathcal{P}}(g)(x)\sim (S_1(x)+S_2(x))\sim 2\sum_{p\le x\atop 
({g\over p})=-1}
{\varphi(p-1)\over p-1},~~x\rightarrow\infty.$$
A slightly more complicated argument leads to the conjecture
$${\mathcal{P}}(g)(x)\sim 2\sum_{p\le x,~(p-1,h)=1\atop ({g\over p})=-1}
{\varphi(p-1)\over p-1},~~x\rightarrow\infty,$$
where $g$ is not required to be squarefree.
The problem in proving this (and hence a quantitative Artin's primitive
root conjecture) is that we have `too many error terms' and can only
achieve this if it can be established that there is sufficient cancellation,
which is presently out of reach.\\
\indent The analog of the latter conjecture for the general index $t$ case
(cf. \S 9.7.3) is far from easy to infer 
by ad hoc methods (as we
managed to do for $t=1$ in the previous section), but can be
merely computed using the approach sketched here. Furthermore, under the
usual RH proviso this conjecture can then be shown to be true
\cite[Theorem 1]{MoreeNP}. A rather easier proof of this result is given
in \cite{Mojapan}.\\
\indent The unconditional results on Artin's conjecture on average
in \cite{Goldfeld, Stephens}, were obtained on using the indicator
function and estimates for character sums.
\subsection{The indicator function in the function field setting}
\indent We now give another application of the identity (\ref{indicator}).
Let $\F_q$ be the finite field with $q=p^n$ elements. 
Recall that its multiplicative group is cyclic.
Let $a(x)$ be 
a polynomial in the polynomial ring $\F_p[x]$. We are interested in the
number of irreducible polynomials $p(x)\in \F_p[x]$ such that
$a(x)$ generates $(\F_p[x]/(p(x)))^*$. Recall that if $p(x)$ is of degree
$n$, then
$\F_p[x]/(p(x))$ is isomorphic with
$\F_{p^n}$, where the isomorphism is given
explicitly as follows: for $h(x)\in \F_p[x]$, we write, using the
Euclidean algorithm,
$$h(x)=p(x)q(x)+r(x),~{\rm with~either~}r(x)=0{\rm ~or~}0\le {\rm deg} \; r
<{\rm deg}\;  p=n.$$
Let $\rho\in \F_{p^n}$ be a root of $p(x)$. Then
$$h(\rho)=r(\rho)=a_0+a_1\rho+\ldots+a_{n-1}\rho^{n-1},~a_i\in \F_p$$
describes all the elements of $\F_{p^n}$. Thus, $a(x)$ generating
$(\F_p[x]/(p(x)))^*$ is equivalent to $a(\rho)$ generating
$\F^*_{p^n}$.\\
\indent Hence, to count the number of irreducible $p(x)$ of degree $n$ 
for which $a(x)$ is a generator is tantamount to counting the number
of $\rho$ of degree $n$ for which $a(\rho)$ generates $\F^*_{p^n}$.
Indeed, since each $p(x)$ has $n$ roots, we find that the number of
irreducible polynomials of degree $n$ in $\F_p[x]$ such that $a(x)$
generates $(\F_p[x]/(p(x)))^*$  is equal to the number of elements
$\rho\in \F_{p^n}$ of degree $n$ such that $a(\rho)$ generates
$\F^*_{p^n}$ divided by $n$. The latter number is by (\ref{indicator}) equal to
$${1\over n}\sum_{\rho\in {\mathbb F}_{p^n}\atop {\rm deg}\rho=n}
{\varphi(p^n-1)\over p^n-1}\sum_{d|p^n-1}{\mu(d)\over \varphi(d)}
\sum_{{\rm ord}\chi=d}\chi(a(\rho)).$$
By inclusion-exclusion we find that the number of irreducible polynomials
of degree $n$ over $\F_q[x]$ equals
$${1\over n}\sum_{\rho\in {\mathbb F}_{p^n}\atop {\rm deg}\rho=n}1
={1\over n}\sum_{d|n}\mu(d)p^{n\over d}={p^n\over n}+
O\left({p^{n\over 2}\over n}\right).$$
(For $t\ge 2$ we have 
$\sum_{d|n}\mu(d)t^{n/d}>t^n-\sum_{j=0}^{[n/2]}t^j>0$ and hence that at least one irreducible polynomial of degree
$n$ over $\F_q[x]$ exists.)
Thus the contribution from the main term (the term with $d=1$) is
$${p^n+O(p^{n\over 2})\over n}{\varphi(p^n-1)\over p^n-1},$$
and the error term is
$${\varphi(p^n-1)\over n(p^n-1)}
\sum_{d|p^n-1\atop d>1}{\mu(d)\over \varphi(d)}
\sum_{{\rm ord}\chi=d}\sum_{\rho\in {\mathbb F}_{p^n}\atop {\rm deg}\rho=n}
\chi(a(\rho)).$$
Applying Weil's estimate for character sums and imposing
the appropriate conditions on $a(x)$, the truth of Artin's conjecture in this
situation is then established, cf. Pappalardi and Shparlinski \cite{PS}.\\
\indent A little knowledge of the history of Weil's estimate makes clear
that we do not always need to invoke this powerful and deep result. Weil,
in a period when he felt depressed, turned to the works of Gauss and 
started reading haphazardly. He soon found himself awestruck by
Gauss' approach in counting solutions mod $p$ of homogenous equations
of the form $x_0^m+\ldots+x_s^m$. This approach involves Gauss sums (and
Jacobi sums). Weil then started musing about generalisations of this
beautiful work and the rest is mathematical history.... 
Thus we would expect that in case 
$a(x)=x^m+c$ the required estimate for the error term can be obtained using only Gauss sums over finite fields.
This is indeed so, see Jensen and Ram Murty \cite{Jensen}.\\
\indent The situation described here is a particular case of the
Artin primitive root conjecture for function fields. Hasse, who
was at the cradle of Artin's conjecture (Artin first formulated it in 
a discussion with Hasse), remained very interested in it throughout
his life and tried to stimulate people to work on it, cf.
\cite{Hasse}. One of
them was his PhD student Bilharz \cite{Bilharz},
 who in 1937 proved
Artin's conjecture for function fields assuming the RH for the zeta function
corresponding to the function field. Indeed, much
later this was proved by Weil and
hence Artin's primitive root conjecture is true for function fields!
The result of Bilharz has been subsequently generalized by 
Hsu \cite{Hsu} to Carlitz
modules, a particular case of 
rank one Drinfeld modules, and in a subsequent paper to rank one Drinfeld modules \cite{Hsuyu}. 
(In the latter case the authors work with an $A$-field $K$ having an injective homomorphism $\gamma:A\rightarrow K$,
in a more recent paper \cite{YY} the case where $\gamma$ is not-injective is addressed.)
Chen et 
al.\,\cite{CKY2} generalized Bilharz's theorem to all one-dimensional tori over global function fields having
a finite constant field.
Rosen  in his well-written book \cite{Rosen} devotes
a whole chapter to the Artin  primitive root conjecture in the function field setting.
\section{Some variations of Artin's problem}
\subsection{Elliptic Artin {\it (by A.C. Cojocaru)}}
Artin's problem is pertinent to many generalizations and variations. 
A natural generalization  appears in the 
context of elliptic curves and was formulated by S. Lang and H. Trotter in 1976 
\cite{LangTrotter}, as follows.

Let $E/\Q$ be an elliptic curve defined over $\Q$ and with arithmetic rank $\geq 1$. Let 
$a \in E(\Q)$ be a point of infinite order. For a prime $p$ of good reduction for $E$, let
$\overline{E}/\F_p$ be the reduction of $E$ modulo $p$ and $\overline{a}$ the reduction
of $a$ modulo $p$. Here, we are also excluding the primes dividing the denominators 
of the coordinates of $a$. 
Following standard notation, we  let $|\overline{E}(\F_p)| = p + 1 -a_p$ denote the number of 
$\F_p$-rational points of $\overline{E}$. The elliptic curve analogue of Artin's problem is to
determine the density of the primes $p$ for which 
$$
\overline{E}(\F_p) = \langle \overline{a} \rangle.
$$ 
In this case, we say that $a$ is a {\it{primitive point of $E$ modulo $p$}}.
In 1976, Lang and Trotter conjectured that this density exists; moreover, they conjectured
that there exists a constant $C_E(a) \geq 0$, depending on $E$ and $a$ such that as $x \rightarrow \infty$, 
$$
\#\{p \leq x: \overline{E}(\F_p) = \langle \overline{a} \rangle\} \sim C_E(a) \frac{x}{\log x}. 
$$

In most cases, this is still an open question. The only result known is due to Rajiv Gupta and
Ram Murty \cite{GM1} who proved that if $E/\Q$ has Complex Multiplication (CM) by the full ring
of integers of an imaginary quadratic field $K$, then under GRH, 
$$
\#\{p \leq x: a_p \neq 0, \overline{E}(\F_p) = \langle \overline{a} \rangle\}
=
C_E(a) \pi(x) + O\left(\frac{x \log \log x}{(\log x)^2}\right).
$$
Additionally, they proved that if 2 and 3 are inert in $K$, or $K = \Q(\sqrt{-11})$, then 
$C_E(a) > 0$; therefore, under GRH, as $x \rightarrow \infty$,
$$
\#\{p \leq x: a_p \neq 0, \overline{E}(\F_p) = \langle \overline{a} \rangle\} \gg \frac{x}{\log x}.
$$

The method used by Gupta and Ram Murty is an elliptic curve adaptation of Hooley's method.
It requires the assumption of GRH, even though it treats the case of
an elliptic curve with CM, which in many instances allows for unconditional results. 
This attests to the difficulty of the elliptic curve analogue of Artin's primitive root conjecture.
No results are known in the case that $E/\Q$ is without CM.

It is a remarkable piece of history that the above result of Gupta and Ram Murty led to the current
best result on the original Artin primitive root conjecture.  Indeed, in their paper, Gupta and
Ram Murty  also considered variations of the elliptic curve analogue of Artin's problem:
instead of working with only one point $a \in E(\Q)$ of infinite order, they also worked
with a free subgroup $\Gamma \leq E(\Q)$ of rank $ > 1$ and considered the
question of counting the primes $p \leq x$ for which 
$\overline{E}(\F_p)= \langle \Gamma (\mmod p) \rangle$. Subsequently,   they considered 
a similar variation of the classical Artin problem which led to the important results
by  Gupta and Ram Murty  \cite{GM} and  by Heath-Brown \cite{HB}.

Another interesting feature of the elliptic curve analogue of Artin's problem 
is that, in contrast with the classical situation for which the group $\F_p^{\ast}$
 is always cyclic,  when we require that $\overline{E}(\F_p) = \langle \overline{a} \rangle$
 we are making two assumptions: first, $\overline{E}(\F_p)$ is cyclic,  and second,
 it is generated by $\overline{a}$. 
 
In general, $\overline{E}(\F_p)$ is the product of two cyclic groups. Therefore it is natural
to consider the first question of finding the density of primes $p$ for which $\overline{E}(\F_p)$
is cyclic.  This question was first studied by Borosh, Moreno and Porta \cite{BMP} who conjectured that
there are infinitely many such primes. A few years later \cite{Serre77} Serre proved this conjecture under GRH,
by using an elliptic curve analogue of Hooley's method. More precisely, by considering the
field extensions $\Q(E[k])$ obtained by adjoining to $\Q$ the $x$ and $y$ coordinates
of the $k$-division points of $E$, he showed that under GRH,
\begin{equation}\label{cyclic}
\#\{p \leq x: \overline{E}(\F_p) \; \mbox{is cyclic}\}
=
C_E \pi(x) + O\left(\frac{x \log \log x}{(\log x)^2}\right),
\end{equation}
where
$$
C_E = \sum_{k \geq 1} \frac{\mu(k)}{[\Q(E[k]) : \Q]}.
$$
In the case that $E$ is with CM,
this result was made unconditional by Ram Murty \cite{Murty1}, and later using a new simpler proof by
A.C. Cojocaru \cite{Co2}.  Ram Murty did not provide an error term, Cojocaru provided the error term
$O(x/((\log x)\log \log x))$. Recently this error term has been considerably sharpened by
A. Akbary and Kumar Murty \cite{akbary1}.
In the case that $E$ is without CM,
the GRH assumption was relaxed
to a `quasi 3/4-GRH'  assumption by A.C. Cojocaru \cite{Co1}.
Unconditionally, it is known by Rajiv Gupta and Ram Murty \cite{GM2} that, for any
elliptic curve $E/\Q$ with $\Q(E[2]) \neq \Q$,  the number of primes $p \leq x$
such that $\overline{E}(\F_p)$ is cyclic is bounded from below by a constant times
$x/(\log x)^2$. Thus, unconditionally, for any such $E/\Q$,
 there are infinitely many primes $p$ for which $\overline{E}(\F_p)$ is cyclic. 

Recently, work of Banks and Shparlinski \cite{BaSh}, combined with work of N. Jones \cite{Jo},
 shows that, on average over a sufficiently large 
two-parameter family of elliptic curves $E$ over $\Q$, 
the cyclicity asymptotic formula $\#\{p \leq x: \overline{E}(\F_p) \; \mbox{is cyclic}\} \sim  C_E \pi(x)$
holds without any additional hypothesis. However, such an average result does not imply that
the formula holds for any curve $E$.

The division fields $\Q(E[k])$ are natural generalizations of the cyclotomic  fields $\Q(\zeta_k)$,
hence it is natural to expect that Hooley's method (which requires the use of the
Kummer extensions $\Q(\zeta_k, a^{1/k})$)
would be amenable to generalization to elliptic curves. However, the non-abelian nature and the size
of the extensions $\Q(E[k])/\Q$ lead to many obstacles in the elliptic curve case that cannot be
overcome as in the classical case. For a nice exposition of such difficulties (such as Brun-Titchmarsh
type results for $\Q(E[k])/\Q$) put in the context of elliptic analogues of other classical 
analytic problems, the reader is refered to Kowalski \cite{Kowalski}. 

For the Lang-Trotter conjecture on primitive points,  the fields 
$\Q(E[k], k^{-1} a)$ that occur  are direct generalizations  of $\Q(\zeta_k, a^{1/k})$. In this case, the
additional complications which arise because of the large size of the conjugacy classes 
involved could be overcome only if $E$ is with CM, under GRH. This is the content
of the work by Gupta and Ram Murty mentioned above \cite{GM1}. 

The positivity of the cyclicity constant $C_E$ was explained  by Serre in \cite{Serre77} and 
proved rigorously
by Cojocaru and Ram Murty in \cite{CojMur}: $C_E > 0$ iff $\Q(E[2]) \neq \Q$. In 
\cite{CojMur}, Cojocaru and Ram Murty
also show that the asymptotic formula (\ref{cyclic}) can be proven with substantially
better remainder terms: $O\left(x^{5/6} (\log x)^{2/3}\right)$ in the non-CM case and 
$O\left(x^{3/4} (\log x)^{1/2}\right)$
in the CM case, provided that the full strength of GRH holds.  These results, together
with the aforementioned ones,   suggest that there is an important
divergence between the cyclicity conjecture for elliptic curves and the Artin primitive 
root conjecture, as the latter does not seem to be  amenable to all the  approaches
of the first one. 

Similar problems have been considered in the context of elliptic curves over 
$\F_q[T]$ in works by Clark and Kuwata \cite{ClKu}, Cojocaru and T\'oth \cite{CoTo},
and Hall and Voloch \cite{HaVo}.  

A different but  related problem was also considered by  
 Chen  and Ju  in \cite{Chenju}.  Given an elliptic curve $E$ over a finite field $\mathbb F_p$ they considered the 
finite field extension $\mathbb F_p(E[ l])$ obtained by adjoining the $ l$-division points of $E$ to 
$\mathbb F_p$ with $ l$ a prime, which has degree $\le l^2 -1$. 
Chen and Ju investigated how often this degree is actually equal to $ l^2 -1$. Under
GRH they showed that when $E$ is not supersingular over $\mathbb F_p$, the set of primes 
$ l$ for which $[\mathbb F_p(E[ l])\colon \mathbb F_p]= l^2 -1$ has a positive density which can be given explicitly in 
terms of $$C_2 = {1 \over 4} \prod_{q> 2}\left( 1 - {2 \over q(q-1)} \right)=0.1337767531541596833\ldots,$$ where the 
product is over the odd primes $p$. They reduce
this problem to a variation of the Artin primitive root conjecture,  which 
subsequently can be dealt with by the methods of Hooley \cite{Hooley1}.\\
\indent We leave it as an exercise for the reader to show that $C_2=T/\pi^2$, with
$$T=2\prod_{q>2}\left( 1 - {1 \over (q-1)^2} \right),$$
the twin prime constant.

\subsection{Even order} 

\label{subsectioneven}
\noindent Let $g\ge 2$ be an integer. As we saw in the introduction, for
a prime $p$, the period of the decimal expansion of $1/p$ in base $g$ is maximal
iff ord$_p(g)=p-1$. 
In 1969 Krishnamurty \cite{krishna} wondered how often the
period in these decimal expansions are even. This can be studied
by algebraic methods similar to those sketched above.\\
\indent We say that a prime $p$ {\it divides a sequence} of integers, if it
divides at least one non-zero term of the sequence.
 Notice
that
$${\rm ord}_p(g)~{\rm is~even}\iff g^x\equiv -1({\rm mod~}p)
~{\rm has~a~solution}\iff p|\{g^k+1\}_{k=0}^{\infty}.$$
There is a unique $j\ge 1$ such that
$p\equiv 1+2^j ({\rm mod~}2^{j+1})$. Now ord$_p(g)$ is even
iff
$$g^{p-1\over 2^j}\not\equiv 1({\rm mod~}p).$$
It is more natural to consider the primes $p$ such that
ord$_p(g)$ is odd, that is the set
$\{p:2\nmid {\rm ord}_p(g)\}$. Let
$$P_j=\{p:p\equiv 1+2^j ({\rm mod~}2^{j+1}),~2\nmid {\rm ord}_p(g)
\}.$$
Note that
$$P_j=\{p:p\equiv 1+2^j ({\rm mod~}2^{j+1}),
~g^{p-1\over 2^j}\equiv 1({\rm mod~}p)\}.$$
Now $P_j$ consists of the primes $p$ that split completely
in $\mathbb Q(\zeta_{2^j},g^{1/2^j})$ but do not split completely
in $\mathbb Q(\zeta_{2^{j+1}},g^{1/2^j})$. Since
$\{p:2\nmid {\rm ord}_p(g)\}=\cup_{j=1}^{\infty}P_j$, and
the $P_j$ are disjoint, one expects that the natural
density of the set
$\{p:2\nmid {\rm ord}_p(g)\}$ equals $\sum_{j=1}^{\infty}
\delta(P_j)$, where $\delta(P_j)=\lim_{x\rightarrow \infty}
P_j(x)/\pi(x)$.
Now
$$\sum_{j=1}^{\infty}\delta(P_j)=\sum_{j=1}^{\infty}\left(
{1\over [\Q(\zeta_{2^j},g^{1/2^j}):\Q]}
-{1\over [\Q(\zeta_{2^{j+1}},g^{1/2^j}):\Q]}\right).$$
Taking $g=2$ one expects by computing the field degrees involved that the density of prime divisors
of the sequence $\{2^k+1\}_{k=0}^{\infty}$ equals
$$\delta=({1\over 2}-{1\over 4})+({1\over 8}-{1\over 8})
+\sum_{j\ge 3}({4\over 4^j}-{2\over 4^j})={7\over 24}.$$
With natural density replaced by Dirichlet density, this
result had been proved in 1966 by Hasse \cite{Hasse2}. Earlier, involving
fewer technical complications, he had considered the Dirichlet density
of primes $p$ such that an odd prime $l$ divides {\rm ord}$_p(a)$ \cite{Hasse1}.
The natural density of primitive divisors of sequences of the form $\{a^k+b^k\}_{k=1}^{\infty}$ with 
$a$ and
$b$ integers was computed independently along similar lines
by Ballot \cite{Ballot}, Odoni \cite{Odoni} and Wiertelak \cite{Wiertelak}.
For the   sequence $\{a^k+b^k\}_{k=1}^{\infty}$ with $a$ and $b$ canonical, the density is
$2/3$. In particular if $g=10$, one finds that the primes $p$ such that the period of the decimal expansion of $1/p$ is
even have density $2/3$.\\
\indent Note that this result, unlike Hooley's, is unconditional.
Both involve infinite towers of field extensions, but the one appearing
in this problem is so sparse (in the sense that the degrees of the
fields quickly increase) that the profusion of error terms that
occur in the Artin primitive root conjecture does not arise and
one can afford working with the relatively large error terms 
in $\pi_K(x)$ that are conditionally known. A probabilistic model for ord$_p(g)$ to be
even was found by Moree \cite{Model}. For a related heuristic approach see Ballot \cite{Ballot4}.\\
\indent Early authors working on divisors of sequences of the form $\{a^k+b^k\}_{k=1}^{\infty}$ were
interested in showing that primes in certain arithmetic progressions occur either all as divisors or
non-divisors, see, e.g., \cite{brauer, schinzel, sier}. Fermat was also interested in this and made
some false claims in one of his letters to Mersenne \cite{mosu, schinzel}. Moree and Sury \cite{mosu}
computed the density of primes $p$ such that $p\equiv c({\rm mod~}d)$ and $p|a^k+b^k$ for some $k\ge 1$ 
in case both $a$ and $b$ are positive integers.\\
\indent Midy (see Leavitt \cite{leav}) proved in 1836 that if $1/p$ has even period $2d$ (that is, 
ord$_p(10)$ is even), then writing 
$${1\over p}=0.(UV)(UV)\ldots,$$ where $U,V$ are blocks of $d$ digits each, one has $U+V=10^d-1$. 
Example: $10^3\equiv -1 ({\rm mod~}7)$, ${1\over 7}=0.{\underline{142857}}142857\ldots$, and $142+857=999$. The result of Midy
was generalized by Gupta and Sury \cite{GS} to the case where the period is $ld$, for some
recent results in this direction see Lewittes \cite{lewi} and H.W. Martin \cite{hmartin}.\\
\indent Instead of asking when $2|{\rm ord}_p(g)$, one can ask
for $m|{\rm ord}_p(g)$ or for $(m,{\rm ord}_p(g))=1$ an so forth. These
questions with sharp error terms are dealt with in various papers
of Wiertelak \cite{weird, wierook2, wiedewiet2, wierook, wiedewiet1}. Again it turns out the corresponding density exists and
is a rational number. Wiertelak \cite{wierook2} gave a complicated expression for this
density for the `$m|{\rm ord}_p(g)$' problem that was substantially simplified by 
Pappalardi \cite{Pappa}. Pappalardi's
expression was simplified by Moree \cite{Modelerd}, who based himself
on a simple to prove, yet previously unknown, identity for the associated counting
function (\cite[Proposition 1]{Modelerd}). For a number field generalization of some of these results see
Perucca \cite{Perucca}.\\

\subsection{Order in a prescribed arithmetic progression}
The problem of determining for how many primes $p\le x$ one has ord$_p(g)\equiv a({\rm mod~}d)$, with
$a$ and $d$ prescribed integers,  was first considered by Chinen and Murata \cite{CM} in case $d=4$. 
By 
a different method their results for $d=4$ are reproved in Moree \cite{Moreejapi1} for more general $g$. 
This gives also a rather more compact formulation of their results. 
In the
same paper also the case $d=3$ was dealt with.\\ 
\indent In this generality the problem turns out to be much harder than the case where $d|a$.
Under an appropriate generalization of RH it turns out that $\delta_g(a,d)$, the density of primes $p\le x$ such
that  ord$_p(g)\equiv a({\rm mod~}d)$ exists (see Chinen and Murata \cite{CM3} in case 
$d$ is a prime power, \cite{CM4, Moreejapi2} in the general case). Let $N_g(a,d)(x)$ be
the associated counting function.
The proof of the existence of $\delta_g(a,d)$  by Moree \cite{Moreejapi1} starts with noting the identity
$$N_g(a,d)(x)=\sum_{t=1}^{\infty}\#\{p\le x:r_p(g)=t,~p\equiv 1+ta({\rm mod~}dt)\},$$
where $r_p(g)=(p-1)/{\rm ord}_p(g)$ denotes the residual index. The individual terms
can be dealt with by a variation of Hooley's method, where now one also needs to keep track
of the dependence on $t$. The more novel aspects of the work rest in making the
resulting formula for $\delta_g(a,d)$ more explicit, which requires algebraic number theory.
Chinen and Murata \cite{CM4} express $\delta_g(a,d)$ as a six fold sum, in
\cite{Moreejapi1} it is expressed as a double sum,
Ziegler \cite{Ziegler} generalized the author's approach to $N_g(a,d)(x)$ to the case where instead
of primes one considers prime ideals, see Problem 10 below for more details.\\
\indent The much simpler quantity 
$\delta(a,d)$, the average density of elements of order $\equiv a({\rm mod~}d)$ in a field of prime
characteristic also exists \cite{Mofield}. It turns out that under the latter RH variant for almost
all $g$ with $|g|\le x$ one has $\delta_g(a,d)=\delta(a,d)$. Moreover, if $\delta_g(a,d)\ne \delta(a,d)$, then
$|\delta_g(a,d)-\delta(a,d)|$ tends to be small.
Nevertheless, it is possible to find specific $g$ (usually appropriately chosen
high powers of a small basis number), for which $\delta_g(a,d)$ shows highly 
non-canonical behaviour. For precise formulations and proofs see \cite{Moreejapi3}.\\ 
\indent Despite
the subtle arithmetic behaviour of $\delta_g(a,d)$, there are some surprises. For example, if $g>0$ then
always $\delta_g(1,4)\le \delta_g(3,4)$. The rule of thumb that a rational density in an Artin type problem
indicates that the associated tower of field extension is rather sparse, seems to be incorrect in this
context.\\
\indent The analogous problem of determining for how many primes $p\le x$ the residual index is in
a prescribed congruence class modulo $d$ was first studied by Pappalardi \cite{Pappa-1}. This turns
out to be a much easier problem.\\
\indent For a more detailed survey of the material in this section, see Moree \cite{elipie}. 
The main results in \cite{Moreejapi1, Moreejapi2} are
also obtained in papers by Chinen and Murata \cite{CM, CM3, CM4} by different methods (the comparison
of $\delta_g(a,d)$ with $\delta(a,d)$ which is the theme of Moree \cite{Moreejapi3}, has 
not (yet) been considered by Chinen and Murata).

\subsection{Divisors of second order recurrences}
The sequences  $\{a^k+b^k\}_{k=1}^{\infty}$ are special cases of
sequences satisfying a second order linear recurrence and one can
wonder about prime divisors of such recurrences in general, cf. Ballot \cite{Ballot}, or
about divisors in general. A good general comprehensive survey on results
on recurrences can be found in the book \cite{recur}. A lot of information on second order recurrences
can also be found in the book by Ribenboim \cite{riben}.\\
\indent A general linear second
order recurrence $R$ has the form $u_{n+2}=eu_{n+1}+fu_n$, with $ef\ne 0$. Once the values
of $u_0$ and $u_1$ are given, it is uniquely determined. The
{\it characteristic polynomial} associated to the recurrence is
defined as $f(X)=X^2-eX-f$. Stevenhagen \cite{Stevenfuture} shows that the set of prime divisors, $\pi_R$, of $R$ is
most conveniently described in terms of two parameters: the {\it root quotient} $r$ of the characteristic
polynomial and the {\it initial quotient} $q$ of the sequence. In this characterization the integer
sequence $R$ no longer plays a role. It greatly clarifies the group structure introduced by
Laxton \cite{Laxton1, Laxton2} on the set of equivalence classes of sequences satisfying a given second
order recurrence. Stevenhagen introduces the notion of an equivalence class $[q,r]$ of a second order
recurrence with root quotient $r\ne \pm 1$ and a sequence group $S(r)$. 
\subsubsection{Prime divisors of torsion second order recurrences}
A second order recurrence with
root quotient $r$ and initial quotient $q$ is said to be {\it torsion} if $[q,r]$ is an element of finite
order in $S(r)$. For example, the sequences $\{a^k+b^k\}_{k=1}^{\infty}$ are torsion.
Another case occurs when $K$ is a quadratic number field, $\epsilon$ its fundamental unit, 
${\tilde \epsilon}$ its conjugate, and $L_n={\epsilon}^n+{\overline \epsilon}^n$ (this is a generalized Lucas
sequence). Stevenhagen proved that a second order torsion sequence has a density of prime divisors
that is positive and rational. In the case of $L_n$ above the result is due
to Moree and Stevenhagen \cite{MStorsion}, who extended a result of Moree \cite{Moreetorsion}, who on his
turn generalized Lagarias' result \cite{Laga} that 
the sequence of Lucas numbers $L_n$ (defined by $L_0 = 2$, $L_1 = 1$ and $L_{n+1} = L_n + L_{n-1}$) has natural density $2/3$. 
\subsubsection{Prime divisors of non-torsion second order recurrences} 
If the characteristic polynomial of the recurrence is reducible over $\Q$,
then the recurrence has the form $R=\{ca^k-db^k\}_{k=0}^{\infty}$, and finding
the prime divisors amounts to finding the prime divisors of the sequence
with general term
$(a/b)^k-d/c=\alpha^k-\beta$. Now $p$ is a prime divisor of $R$, with
at most finitely many exceptions, iff
ord$_p(\beta)|{\rm ord}_p(\alpha)$. 
This can be phrased differently by saying that, with at most finitely
many exceptions, $p$ divides $R$ iff the subgroup mod $p$
generated by $\beta$, $\langle \beta \rangle$, is contained in the
subgroup $\langle \alpha \rangle$. One can then base a two variable
Artin conjecture on this, of which the qualitative form says that the
sequence $R$ should have infinitely many prime divisors. This is not difficult to prove, 
see P\'olya \cite{poll}. 
(Schinzel \cite[pp. 909-911]{Selecta} proved that if $a$ and $b$ are
rational integers with $a>0$ and $b\ne a^e$, then 
there are infinitely many primes $p$ that do {\it not} divide the sequence $\{a^k-b\}_{k=1}^{\infty}$.)
To investigate the quantitative form, we put
$S(x)=\sum_{p\le x,~p|R}1$.
Then
$$S(x)=\sum_{p\le x}\sum_{w|p-1\atop
{\rm ord}_p(\beta)|{\rm ord}_p(\alpha)=(p-1)/w}1
=\sum_{w\le x-1}\sum_{{p\le x\atop p\equiv 1({\rm mod~}w)}
\atop {\rm ord}_p(\beta)|{\rm ord}_p(\alpha)=(p-1)/w}1.$$
Now $p\equiv 1({\rm mod~}w)$ and
ord$_p(\beta)|{\rm ord}_p(\alpha)=(p-1)/w$ iff 
$p$ splits completely in
$\mathbb Q(\alpha^{1/w},\beta^{1/w},\zeta_w)$
and $p$ does not split completely in any of the fields
$\mathbb Q(\alpha^{1/qw},\beta^{1/w},\zeta_{qw})$ with $q$ a prime.
Denote the degree of the latter field by $n(qw)$. Then,
as $x\rightarrow \infty$, the limit $S(x)/\pi(x)$ tends
to
\begin{equation}
\label{Stephensdouble}
\sum_{w=1}^{\infty}\sum_{k=1}^{\infty}
{\mu(k)\over n(kw)},
\end{equation}
provided a sufficiently
strong generalisation of RH is true. The sum 
(\ref{Stephensdouble}) equals a rational number times
the {\it Stephens constant}
$$S=\prod_q\left(1-{q\over q^3-1}\right)=0.57595 99688 929\ldots $$
This result is due to Moree and Stevenhagen \cite{MS}
and improves on earlier work by Stephens \cite{Stephens3}. 
Stephens evaluation of (\ref{Stephensdouble})
(who restricted himself to the case where both $\alpha$ and $\beta$ are integers) is not always correct
and was corrected in \cite{MS}. Using the average-character-sum method 
\cite{LMS} a rather easier reproof can be given \cite{Most}. The Stephens constant, just like $A$, can also be evaluated with high numerical precision by
expressing it in terms of zeta values at integer arguments \cite{Monumi}.\\   
\indent A naive heuristic argument would lead us to expect that the
density equals the average density (over
the primes) of pairs $(u,v)\in \mathbb F_p^*\times
\mathbb F_p^*$
such that $u\in \langle v\rangle$. It is not difficult to see that
the density of the pairs $(u,v)\in \mathbb F_p^*\times
\mathbb F_p^*$
such that $u\in \langle v\rangle$ equals
$${1\over (p-1)^2}\sum_{d|p-1}d\varphi(d).$$
Thus the naive heuristic would give
$$\lim_{x\rightarrow \infty}{1\over \pi(x)}\sum_{1<p\le x}{1\over (p-1)^2}
\sum_{d|p-1}d\varphi(d)$$ as density. The above limit was evaluated by Luca \cite{lucka} to
be the Stephens constant, thus answering an open problem in von zur Gathen et al.\,\cite{veel}.
Luca's proof is quite similar to the evaluation of $\sum_{p\le x}\varphi(p-1)/(p-1)$ 
given in \S 7. Another interpretation of Stephens' constant, see 
Korolev \cite{Korolev}, is as the average value of
$${1\over \varphi(n)}\sum_{k=1\atop (k,n)=1}^n{1\over (k-1,n)}.$$
\indent Note that the problem considered here
is a variant of Artin's primitive root conjecture that is more
complicated than the original one, whereas the variant discussed
in \S \ref{subsectioneven} is easier, although both deal with second order recurrences. 
The variants discussed in \S \ref{subsectioneven} are torsion sequences. For every torsion sequence the prime
density exists and is a rational number. It can be explicitly determined. Here the case when the 
characteristic polynomial is reducible over
$\mathbb Q$ is easier, than when it is irreducible. For a non-torsion sequence assuming GRH the density can
be shown to exist and also be explicitly determined. Here again the case where the characteristic polynomial is
irreducible is the most difficult. An example is given by the sequence  
$a_n$ defined by $a_0 = 3$, $a_1 = 1$ and $a_{n+1} = a_n + a_{n-1}$. Lagarias \cite{Laga} raised the challenge
here of determining the density of prime divisors.
Moree and Stevenhagen \cite{mos} showed that under an appropriate generalisation of
RH  this density equals 
$${1573727\over 1569610}S=0.577470679956\ldots$$
Indeed, for every non-torsion sequence assuming GRH the density exists, can be explicitly evaluated and is a rational multiple
of the Stephens' constant $S$.\\
\indent For the convenience of the reader we give some references for the various cases:\\
-Torsion, reducible char. polynomial: \cite{Ballot, Hasse1, Hasse2, Odoni, Wiertelak}\\
-Torsion, irreducible char. polynomial: \cite{Laga, Moreetorsion, MStorsion, Stevenfuture}\\
-Non-torsion, reducible: \cite{Stephens3, MS}\\
-Non-torsion, irreducbile: \cite{mos}.\\
\indent Recently in addition to the techniques used in the above articles, also the theory of dynamical systems has been 
brought to bear on the issue of prime divisors of sequences, see, e.g., \cite{gott, Rafe, Rafe2}.\\
\subsubsection{General divisors}
One can also try to count the number of divisors $\le x$ of a given sequence. 
In case the sequence is $S_{p,1}$, where $S_{a,b}$ denotes the
sequence $\{a^k+b^k\}_{k=1}^{\infty}$, this plays a role in the analysis of the average
behaviour of Ulmer's rank formula, see \cite{PoS}.
The divisor counting problem also has some
applications outside number theory. For example, Pless et al.\,\cite[Theorem 3]{PSQ} have shown that
a nontrivial cyclic self-dual code of length $n$ exists iff $n$ does not divide the sequence 
$S_{2,1}$. Kanwar and L\'opez-Permouth \cite[Theorem 4.5]{KL} proved that if $p$ is a prime 
not dividing $n$
and $m$ is even, then nontrivial self-dual cyclic $\mathbb Z_{p^m}$-codes exist iff $n$ does not
divide $S_{p,1}$.\\
\indent In Moree \cite{PSQ} it is shown that for any $0<\epsilon\le 1/24$ and $r$ the smallest
integer such that ${1\over 3}(1/2^{r+3})\le \epsilon$, there exist positive constants $c_1,\ldots,c_r$
such that $G(x)$, the number of divisors $n$ of the sequence $S_{2,1}$ that are $\le x$, 
satisfies
$$G(x)={x\over \log x}\Big(c_1\log^{1/3}x+c_2\log^{7/24}x+\sum_{k=0}^r c_{3+k}\log^{2^{-k-3}/3}x
+O(\log^{\epsilon}x)\Big).$$
Similar results for divisors of the sequence $S_{a,b}$, with $a$ and $b$ integers
are derived in Moree \cite{mora}, and for Lucas numbers in Moree \cite{mora2}\\
\indent The divisibility of 
certain sequences is related to the level (Stufe) of a field $F$, $s(F)$, which is the smallest integer $s$ (if it
exists) such that $-1=\alpha_1^2+\ldots+\alpha_s^2$ with $\alpha_i$ in $F$.
In case $-1$ cannot be written as a sum of squares from $K$ we put $s(K)=\infty$.
Pfister \cite{Pfister} (cf. \cite[pp. 203-208]{GrosswaldBoek}) proved that in case $s(F)$ is finite we have $s(F)=2^j$ for some $j\ge 0$.
Hilbert (cf. \cite[XI, Theorem 1.4]{lammetje}) proved that if $F$ is an algebraic number field, then $s(F)\le 4$. It
follows that $s(F)\in \{1,2,4\}$ in this case. Note that $s(F)=1$ iff $\sqrt{-1}\in F$.\\
\indent Let us put $K_n=\mathbb Q(\zeta_n)$. If $4|n$, then $s(K_n)=1$. If
$n$ is odd, then clearly $s(K_{2n})=s(K_n)$ since $K_n=K_{2n}$. Thus we may assume
that $n$ is odd. P. Chowla \cite{C} proved that $s(K_p)=2$ when $p\equiv 3({\rm mod~}8)$
is a prime. In later unpublished papers John H. Smith and P. Chowla proved
independently that $s(K_p)=2$ also when $p\equiv 5({\rm mod~}8)$. In 1970, P. Chowla and S. Chowla \cite{CC} 
proved that $s(K_p)=4$ when $p\equiv 7({\rm mod~}8)$.
Fein et al.\,\cite{FGS} proved that for an odd prime $p$ we have
$s(K_p)=2$ iff $p|S_{2,1}$ (and so $s(K_p)=4$ iff $p\nmid S_{2,1}$). 
Apparently unaware of Hasse's work \cite{Hasse1}, they gave their own proof of the
fact that the Dirichlet density of prime divisors of $S_{2,1}$ is 17/24.
Moree \cite[Theorem 7]{mora} 
gave an asymptotic for the number of integers $m\le x$ such that $s(K_m)=4$.\\
\indent Recently Nassirou \cite{nassi} considered the level of ${\mathbb Q}_p(\zeta_n)$
with $p$ odd, where $\mathbb Q_p$ denotes the $p$-adic field. Since
$s({\mathbb Q}_p)=1$ when $p\equiv 1({\rm mod~}4)$, we may assume that 
$p\equiv 3({\rm mod~}4)$. Let $q\ne p$ be an odd prime. The results of Nassirou
imply that $s({\mathbb Q}_p(\zeta_q))=1$ iff $q|S_{p,1}$ and 
$s({\mathbb Q}_p(\zeta_q))=2$ iff $q\nmid S_{p,1}$.\\
\indent In the above we considered a fixed sequence $a^k+b^k$. In the context of supersingular
Fermat varieties, a reverse problem comes up. Here one fixes an integer $m$ and asks for
the density $\delta(m)$ of primes $p$ such that $m$ divides the sequence 
$\{p^k+1\}_{k=1}^{\infty}$. For this Yui \cite{Yui} and independently Waterhouse 
\cite{Water} gave an explicit formula. 
It is known that the Jacobian variety of a Fermat curve $X^m+Y^m=Z^m$ considered in
characteristic $p$, is supersingular (i.e. isogenous to a product of supersingular
elliptic curves), iff $m$ divides the sequence $\{p^n+1\}_{n=1}^{\infty}$, 
see, e.g., \cite{Shioda, Water, Yui}.
Using the explicit formula of Yui and Waterhouse, Schwarz and Waterhouse \cite{SWat}
gave an asymptotic formula for $\sum_{2<m\le x}\delta(m)$. This was improved on by
Nowak \cite{Nowak}, who proved that for every fixed $K$ there exist positive constants
$A_0,A_1,\ldots,A_K$ such that,
$$\sum_{2<m\le x}\delta(m)=x\sum_{k=0}^K A_k(\log x)^{2^{-k}/3-1}+O(x(\log x)^{2^{-K}/4-1}),$$
where the implied constant may depend on the integer $K$ and $x$ tends to infinity. Note
the similarity of this formula with that given above for $G(x)$.

\subsection{Lenstra's work}
\noindent The
most far reaching generalisation of
the Artin primitive root conjecture was
considered by Lenstra \cite{Lenstra}, in the
context of his research on Euclidean number
fields: 
Let $K$ be a global field (that is a finite
extension of $\mathbb Q$ or a function field in one variable over
a finite field) and $F$ a finite normal extension of $K$.  Let
$C$ be a subset of Gal$(F/K)$ which is stable under conjugation and let
$d$ be a positive integer (coprime to the characteristic of $K$ in the case
of a function field).  Consider a finitely generated subgroup $W$ of $K^*$
which has, modulo torsion, rank $r\geqq 1$, and let $M$ be the set of
prime ideals $\mathfrak P$ of $K$ satisfying
\begin{enumerate}
\item the Artin symbol $({\mathfrak P},F/K)\subseteq C$
\item the normalized exponential valuation attached to $\mathfrak P$ satisfies
ord$_{\mathfrak P}(w)=0$ for all $w\in W$.
\item if $\psi:W\to {\bar K}_{\mathfrak P}^*$ is the natural map, then
the index of $\psi(W)$ in ${\bar K}_{\mathfrak P}^*$ divides $d$.
\end{enumerate}
Here ${\bar K}_{\mathfrak P}^*$ denotes the multiplicative group of
${\bar K}_{\mathfrak P}$, the residue class field at $\mathfrak P$.
Lenstra conjectured that $M$ has a density.  He also obtained
necessary and sufficient conditions for this density to be nonzero under
a sufficiently strong generalisation of RH. His main result implies
Theorem \ref{lenstrarekenkundigerij} for example. Lenstra applied his
results
to show that under the usual
RH proviso certain principal ideal rings are Euclidean (see the next section).

\subsection{Artin's conjecture and Euclidean domains {\it (written with Hester Graves)}}
A Euclidean algorithm for an integral domain $R$ is a map $\phi:R\backslash\{0\}\rightarrow W$, a well-ordered set, such that for all $a,b\in R$ with $b\ne 0$, there exist $q,r\in R$ with $a=qb+r$ and $\phi(r)<\phi(b)$. From every Euclidean algorithm, one can create another Euclidean algorithm $\phi:R\backslash\{0\}\rightarrow W$, a well-ordered set, such that 
$\phi(0)=0$ and
for all $a,b\in R$ with $b\ne 0$, there exist $q,r\in R$ with $a=qb+r$ and $\phi(r)<\phi(b)$ or $r=0$. We will concern ourselves with this second type of Euclidean algorithm.  \\
\indent Nagata showed that one can replace ``$W$, a well-ordered set '' with ``$W$ an ordered set with minimum condition'' in the above definition.  He was nonetheless able to show that every such Euclidean algorithm implied the existence of a Euclidean algorithm that mapped to a well-ordered set.  In other words, his definition expanded the set of possible Euclidean algorithms, but not the set of possible Euclidean rings.\\ 
\indent If an integral domain has a Euclidean algorithm, then all of its ideals are principal. 
The converse of this is not true in general.  Not all principal rings of integers of imaginary quadratic fields are Euclidean.  The nine quadratic imaginary number fields with class number one are  $\mathbb Q(\sqrt{-d})$ with 
 $d=1,2,3,7, 11, 19,43,67$ and $163$; only the rings of integers of $\mathbb Q(\sqrt{-d})$ with $d= 1,2,3,7,$ or $11$ are Euclidean.\\
 \indent An integral domain $R$ equipped with a Euclidean algorithm is called a Euclidean domain. For algebraic number fields the
Euclidean algorithm most often studied is the absolute value of the norm. Such fields and their rings of integers are called
norm-Euclidean. 
The classification of rings of algebraic integers which are Euclidean (not necessarily 
for the norm function) is a major unsolved problem.
Dedekind showed in Supplement XI to Dirichlet's {\it Vorlesungen \"uber Zahlentheorie} \cite{Diri}
that $\mathbb Q(\sqrt{d})$ is norm-Euclidean for $d=-11,-7,-3,-2,-1,2,3,5,13$.  Now it's known that the
full list of norm-Euclidean quadratic fields also includes $d=6,7,11,17,19,21,29,33,37,41,57$ and $73$.
Lemmermeyer \cite{Lemmermeyer-1} determined all real cubic norm-Euclidean fields with discriminant
less than 4692. Clark \cite{cclark} gave two examples of cubic fields which are Euclidean but not
norm-Euclidean. Further examples were found by Cavallar and Lemmermeyer \cite{cavlem}.\\
\indent Motzkin \cite{Motzkin} constructed the so called minimal Euclidean algorithm. Given a nonempty
collection of Euclidean algorithms $\phi_{j}$ on $R$, the map defined by
$\phi(r)=\min_{j}\phi_j(r)$ is easily seen to also be a Euclidean algorithm on $R$. If the minimum
is taken over all the Euclidean algorithms of the ring $R$, the resulting algorithm $\phi$ is
called the minimal algorithm. 
Let $A_0=\{0\}$, and define $A_j$ for $j\ge 1$ by
$$A_j-A_{j-1}=\{a\in  R:{\rm~each~residue~class~of~}R/(a){\rm~contains~an~element~of~}A_{j-1}\}.$$
Note that $A_1-A_0$ is the set of units of $R$.
If
\begin{equation}
\label{euclidie}
R=\cup_{j=0}^{\infty}A_j,
\end{equation}
then $R$ is Euclidean, with $\phi(a)=\min\{j:a\in A_j\}$. Motzkin \cite{Motzkin} proved that the
condition (\ref{euclidie}) is also necessary if $R$ is to be Euclidean.
As an example, if $R =\mathbb Z$, then $A_1=\{-1,0,1\}$. This set consists of three consecutive
integers and thus it provides representatives for the classes mod 2 and mod 3, so that
$A_2=\{-3,-2,-1,0,1,2,3\}$. Here we have 7 consecutive integers so that $A_3$ is the interval
$[-7,+7]$ consisting of 15 consecutive integers. It can be easily shown that if $\theta(n)$ denotes
the number of binary digits of $|n|$, then the minimal algorithm on $\mathbb Z$ is given by $\theta$.
Using the sets $A_n$ we can easily
see, for example, why $\mathbb Q(\sqrt{-19})$ is not Euclidean for any algorithm. This follows because
$\pm 1$ are its only units and every non-trivial ideal has norm at least 4 so the construction
of the sets $A_n$ stops at $A_1$.  If the reader wants to learn more about Euclidean rings,  Lenstra \cite{Leneuclid} wrote an introductory account. Samuel \cite{samuel} wrote a more advanced survey on Euclidean rings.  In this survey, he studied the function $\theta$ defined below in Theorem \ref{Lenstratheta} and noticed that the condition on the value of $n_{\mathfrak p}$ (defined in order to study $\theta$) was similar to Artin's conjecture.  He may have been the first to notice the connection between Euclidean rings and Artin's conjecture.  In addition, Samuel asked questions that shaped the study of Euclidean rings.  For example, he asked whether there were a real quadratic number field such that its ring of integers were Euclidean but not norm-Euclidean; he suggested the  ring $\Z[\sqrt{14}]$ as a candidate. \\
\indent Two years later, assuming GRH, Weinberger \cite{Wein1} was able to use Hooley's work on Artin's conjecture in the Euclidean ring situation and 
established the following result.
\begin{theorem} 
\label{blubber}
{\rm (GRH)}.
Let $K$ be an algebraic number field whose ring of integers both is a principal
ideal domain and has infinitely many 
units. Then the ring of integers of $K$ is Euclidean.
\end{theorem}
\indent At first sight GRH seems to have little to do
with the statement of the result. It comes from the
fact that the following variant of
Artin's primitive root conjecture is needed for the proof of Theorem \ref{blubber} and
can presently be only proved assuming GRH. 
\begin{theorem} 
\label{bups}
{\rm (GRH)}. Let $K$ be an algebraic number field whose ring of integers $R$ both is a principal
ideal domain and has infinitely many units. Let $\epsilon$ be a fundamental unit of $R$, i.e.
$\epsilon\ne \epsilon_1^n$ for all units $\epsilon_1$ and all integers $n>1$. Let $\mathfrak p$ be a
prime ideal of $K$ and let $I$ be the multiplicative group of ideals of $K$ prime to $\mathfrak p$ and
let $H=\{\mathfrak a\in I:\mathfrak a\equiv (1)~ ({\rm mod~}\mathfrak p)\}$. Then every ideal class of $I/H$
contains infinitely many prime ideals $\mathfrak q$ for which $\epsilon$ is a primitive root
mod $\mathfrak q$, i.e. the class of $\epsilon$ generates the (cyclic) group of units in the finite field $R/\mathfrak q$.
\end{theorem}
The proof of this result is a variant of Hooley's proof of the original Artin conjecture.\\
\indent The proof that Theorem \ref{bups} implies Theorem \ref{blubber} consists of two steps. We first define irreducibility; an element $b$ is irreducible if the ideal $(b)$ is prime.  Assuming GRH, one uses Theorem
\ref{bups} to show that every irreducible  of $R$ is in $A_3$. If one chooses an irreducible element $p$, then every ideal class of $I/H$ has infinitely many prime ideals $\mathfrak q = (q)$ such that $(R/ \mathfrak q)^{\times} = \langle [ \epsilon ] \rangle.$ This implies the $q$ are in $A_2$ and therefore, for all $a \in R-(p)$, there exists an irreducible $q \in A_2$ such that $a \equiv q \pmod p.$ We conclude that $p$ is in $A_3.$ \\
\indent This means that every element of $R -A_1$ can be written as a product of irreducibles in $A_2 - A_1$ and in $A_3 - A_2.$  If $n_2$ is the number of these irreducibles in $A_2 - A_1$ and $n_3$ is the number of these irreducibles in $A_3 - A_2$, then we define the height of $b$ to be $\rm ht(b) =$  $2n_2 + 3n_3.$  One then shows that if $b \in R-A_1$, then $b \in A_{\rm ht(b)}$.  Thus the sets $A_j$ exhaust $R$ and
so $R$ is Euclidean.\\
\indent Let $K$ be a global field, and let $S$ be a non-empty set of prime divisors of $K$ containing
the set $S_{\infty}$ of archimedean prime divisors of $K$. The ring of $S$-integers of $K$ is 
$R_S=\{x\in K:\nu_{\mathfrak p}(x)\ge 0~{\rm ~for~all~primes~}p\not\in S\}$. Lenstra \cite{Lenstra} generalized the result of Weinberger
above as follows.
\begin{theorem}
\label{Lenstratheta} Suppose that $R_S$ is a principal ideal ring, and that $\# S\ge 2$. Further, if
$K$ is a number field, assume that for every squarefree integer $n$ and every finite subset
$S'\subset S$ the Dedekind zeta function $\zeta_{K(\zeta_n,{R_S^*}^{1/n})}(s)$ satisfies the Riemann
Hypothesis. Then $R_S$ is Euclidean, and its minimal algorithm $\theta$ is given by
$$\theta(x)=\sum_{\mathfrak p\not\in S}\nu_{\mathfrak p}(x)n_{\mathfrak p}~(x\in R_S,~x\ne 0),$$
where the sum is over all primes $\mathfrak p$ of $K$ which are not in $S$ and
$$n_{\mathfrak p}=
\begin{cases}
1 & {\rm if~the~natural~map~}R_S^*\rightarrow {\overline K_{\mathfrak p}}^*{\rm ~is~surjective};\cr
2 & {\rm else}.
\end{cases}
$$
\end{theorem}
The function field case of this result is due to Queen \cite{queen}. Lenstra also points out
that the Euclidean algorithm given by Weinberger is not the minimal one. The assumption $\# S\ge 2$
is only an apparent restriction in the latter result, as in case $\# S=1$ Lenstra gives the
complete list of $R_S$, both for the number and function field situations, that are principal and those that are Euclidean.\\
\indent Under some further assumptions, the GRH condition
is not needed. 
Gupta, Ram Murty, and Kumar Murty \cite{GMM} established the following result.
\begin{theorem} Let $K$ be a number field  Galois over $\mathbb Q$ and let $g$ be the gcd of the 
numbers $N_{K:\mathbb Q}(\mathfrak P) -1$ for $\mathfrak P\in S-S_{\infty}$.  If \\
\indent (a) $|S|\geq\max(5,2[K:\mathbb Q] -3)$, and \\ 
\indent (b) $K$ has a real embedding or $\zeta_g\in K$,\\
then the ring  of S-integers $R_S$ is principal iff it is Euclidean.
\end{theorem}
The proof combines ideas from a paper by Gupta and Ram Murty on Artin's conjecture
\cite{GM} and from a paper of Ram Murty and Kumar Murty \cite{murtymurty}
on an analogue of the Bombieri-Vinogradov theorem.\\
\indent Harper and Ram Murty published further results on Euclidean rings that did not assume GRH.  In his thesis \cite{Harper}, Harper adapted Motzkin's construction of the sets $A_j$ defined above in the number field situation.  He decided to define $B_0 = \{ 0 \cup R^* \}$ and $B_i = B_{i-1} \cup \{ \textnormal{irreducible elements } p : R^* \twoheadrightarrow (R/(p))^* \}.$ Using Dirichlet's theorem on primes in arithmetic progressions, he showed  that if $R$ is principal and $\cup B_i$ contains every irreducible element of $R$, then $R$ is Euclidean.  He then further refined the sets by augmenting $B_0$ by an `admissible set of primes.'   A set of primes $\{\pi_1,\dots,\pi_s\}$ is called admissible if, for all integers 
$\beta$ composed of only these primes, every coprime residue class modulo $\beta$ can be represented 
by a unit. He then used the large sieve to show that if $|\{ b \in B_1: \textnormal{Nm}(b) \leq x \} | \gg \frac{x}{\log^2 x}$, then $\cup B_i$ contains all of the irreducible elements of $R$.   In particular, if a real quadratic field has two admissible primes, then its ring of integers is principal iff it is Euclidean.   He found two such primes in the ring $\Z[\sqrt{14}]$ and used them to show that $\Z[\sqrt{14}]$ is Euclidean, thereby answering Samuel's question.  He then used the same techniques to show that if the ring of integers of a cyclotomic field, $\Z[\zeta_d]$, is principal, then it is Euclidean.\\
\indent Harper and Ram Murty  \cite{HM} generalized this work by proving that if $K$ is a finite Galois extension of $\Q$ with unit rank greater than 3, and if the ring of integers of $K$ is principal, then there exists a large enough set of admissible primes, implying that the ring of integers is a Euclidean domain. The truth of this result assuming GRH was previously established by
Weinberger (Theorem \ref{blubber}).\\
\indent Petersen and Ram Murty extended Ram Murty's earlier work  with Harper by removing the Galois condition in certain cases.  They showed that if $K$ is a number field with ring of integers $R$, $\mathrm{rank}(R^*) > 3$ and $M$ is a subfield of $K$ such that $K$ is a Galois extension of $M$ of degree $> 3$,  then $R$ is a principal ideal domain iff $R$ is a Euclidean ring \cite{EuclidMurtyPetersen}. It is worth noting that Petersen and Ram Murty did other work involving Artin's conjecture.  They showed that if $1 \leq a <q$, with $(a,q)=1$, and if $K$ is a finite Galois extension of $\Q$ such that $\mathrm{rank}(R^*) >3$ and $\zeta_q \notin K$, then there are infinitely many prime ideals $\mathfrak{p}$ of first degree in $R$ such that $|\mathrm{Nm}(\mathfrak{p})| \equiv a \pmod{q}$ and $R^* \twoheadrightarrow (R / \mathfrak p)^{\times}$ \cite{MurtyPetersen}.  They used this to show that if $K$ is a finite Galois extension of $\Q$, $\mathrm{rank}(R^*) >3$ and $\sqrt{-1} \notin K$, then there are infinitely many maximal $h_K$-cusped subgroups of $\mathrm{PSL}_2(R)$ \cite{MurtyPetersen}. Previously, Petersen \cite{Petersenalone} proved such a result for all $K$ with positive unit rank and $\sqrt{-1}\not \in K$,  assuming the GRH.   As all of these subgroups are congruence subgroups, this result is a contrast to situation when $K$ is $\Q$ or an imaginary quadratic, where there are many maximal $h_K$-cusped subgroups, but few of these are congruence subgroups. (cf. \cite{Petersson} and \cite{OnePetersen}.)   The geometry is discussed in  \cite{MurtyPetersen, Petersenalone, OnePetersen}.\\
\indent Harper's proof that $\mathbb Z[\sqrt{14}]$ is Euclidean, but not norm-Euclidean, answered a question 
that had inspired a great deal of mathematical activity. The earliest known example 
of a quadratic Euclidean field that is not norm-Euclidean is
$\mathbb Q(\sqrt{69})$ (\cite{Clark, Niklasch}).  Since Samuel asked the question,  much has been done to find a function
other than the norm
which makes $\mathbb Z[\sqrt{14}]$ Euclidean (see \cite{Lemmermeyer-1}). Nagata was unable to prove that $\mathbb Z[\sqrt{14}]$ is Euclidean.  Instead, he introduced the concept of a `pairwise algorithm' and showed that $\mathbb Z[\sqrt{14}]$ has such an 
algorithm (\cite{Nagatapair}).
Clark and Ram Murty \cite{CMurty} proved that ${\mathbb Z}[\sqrt{14},1/p]$ is Euclidean for 
$p=1298852237$. On
combining their techniques with earlier work of Gupta, Ram Murty and 
Kumar Murty \cite{GMM}, Harper \cite{Harper} observed that it can be shown that 
${\mathbb Z}[\sqrt{14},1/p]$ is Euclidean for any $p$. \\
\indent  In 1995, Lemmermeyer \cite{Lemmermeyer-1} wrote
a survey of known results and open questions on Euclidean and non-Euclidean algebraic number fields.   It was updated in 2004, but in 2007, Narkiewicz proved that there are at most two real quadratic fields such that their ring of integers is principal but not Euclidean \cite{Narkeuclid}.  He also showed that there is at most one normal cubic number field with class number one that is not Euclidean.  He proved this by combinining Harper and Ram Murty's work with his earlier results (Narkiewicz \cite{Nark1}) on units in residue classes.   There have been no other major advances in the study of Euclidean rings since then, but there has been work on Euclidean ideals and Euclidean systems.\\
\indent In 1979, Lenstra \cite{lenst} introduced the notion of a Euclidean ideal class which generalizes that of
a Euclidean ring.  As the existence of a Euclidean algorithm for a domain implies trivial class group, the existence of a Euclidean ideal in a Dedekind
domain implies cyclic class group.  Lenstra was inspired by norm-Euclidean rings of integers.  Suppose that $K$ is a number field and that $\O_K$ is its
ring of integers.  
Let $\Nm$ denote the field norm.
If, for all $a,b \in \O_K$, $b \neq 0$, there exists some $q,r \in \O_K$ such that 
$$a = qb + r,  \textnormal{ where } r = 0 \textnormal { or  }|\Nm (r)| < |\Nm (b)|,$$
then 
$$|\Nm \left ( \frac{a}{b} - q \right )| = |\Nm \left ( \frac{r}{b} \right )| < 1.$$
In other words, $R$ is norm-Euclidean iff, for all $x \in K$, there exists some $y \in \O_K$ such that 
$$|\Nm (x-y)| < 1 = \Nm (R).$$
Lenstra then asked what would happen if $R$ were replaced by some fractional ideal $C$.  If, for all $x \in K$, there exists some $y \in C$ such that 
$$|\Nm(x-y)| < |\Nm(C)|,$$
then $C$ is a norm-Euclidean ideal. 
The only quadratic number fields with a norm-Euclidean ideal are $\Q(\sqrt{d})$, with $d= -1,-3,-5,-7,-11,-15,-20, 3,5,6,7,8,10$ $11,13,15,17,19,21,$ $29,33,37,41,57,73,$ and $85$ (\cite{lenst},\cite{Graves/Ramsey}).  Note that only for  $d= -15, -5,10,15,$ and $85$, does $\Q(\sqrt{d})$ have a non-principal norm-Euclidean ideal.\\
\indent As mentioned above, there are rings that are Euclidean but not norm-Euclidean.  Lenstra generalized norm-Euclidean ideals as follows.  For the rest of the section, we define $E$ to be the set of all inverses of integral ideals of a Dedekind domain $R$. Suppose that $C$ is a fractional ideal of $R$. If there exists a function $\psi: E \longrightarrow W$, $W$ a well-ordered set, such that for all $I \in E$ and all $x \in IC \setminus C$, there exists some $y \in C$ such that $$\psi((x +y )^{-1} I C) < \psi(I).$$  We say $\psi$ is a $\textit{Euclidean algorithm for  } C$ and $C$ is a \textit{Euclidean ideal}.\\   
\indent If we define $\mathbb{I}$ to be the set of all integral ideals of $R$, we can rewrite this definition. A fractional ideal $C$ is  $\textit{Euclidean}$ if there exists a function $\psi: \mathbb{I} \longrightarrow W$, $W$ a well-ordered set, such that for all integral ideals $I$ and all $x \in I^{-1} C \setminus C$, there exists some $y \in C$ such that $$\psi((x -y ) I C^{-1}) < \psi(I).$$  \\
\indent If $C$ is a Euclidean ideal, then so is every ideal in its ideal class.  If $C$ is a Euclidean ideal, then its ideal class $[C]$ is called a \textit{Euclidean ideal class}.  One can check that if $R$ is a Euclidean ring, then $R$ is a Euclidean ideal and $[R]$ is a Euclidean ideal class. If $C$ is a Euclidean ideal, then its ideal class generates the class group of $R$, implying that the class group is cyclic.  Lenstra then proceeded to prove a generalization of Weinberger's result for Euclidean ideals.  \\
\indent Lenstra showed \cite{lenst} by generalizing Weinberger's work on Euclidean rings, that 
assuming GRH if $K$ is a number field with cyclic class group and with $\O_K$ having infinitely many units, then $[C]$ generates the class group of $K$ iff $C$ is a Euclidean ideal.  To be more precise,  the assumption is that for every squarefree $n$, the $\zeta $ function associated to $K((R^*)^\frac{1}{n})$ satisfies the RH.\\
\indent As in the Euclidean ring situation, there has been work to prove this without assuming the GRH.  In her thesis \cite{Graves}, Graves generalized Motzkin's reformulation of the Euclidean algorithm and Harper's work on Euclidean rings to the Euclidean ideal situation.  One can define analogous constructions for the $A_i$'s and $B_i$'s.  One can then, using a variation of the large sieve, prove Harper's theoretical results for the Euclidean ideal situation.  They imply that if $K$ is a number field such that $|\O_K^{\times}| = \infty$,  if $[C]$ generates $\mathrm{Cl}_K$, and if
$$\left | \left \{
 \textnormal{prime ideals } 
\mathfrak p 
 \left | \textnormal{Nm}(\mathfrak p) \leq x,[\mathfrak p] = [C], R^* \twoheadrightarrow (R / \mathfrak p)^{\times} \right. \right \} \right | \gg \frac{x}{\log^2 x},$$ then $C$ is a Euclidean ideal.   In a later paper \cite{Gravesexample}, she used this along with a result of Narkiewicz \cite{Nark1}, to show that $\Q(\sqrt{15}, \sqrt{35})$ has a Euclidean ideal.\\
 \indent Recently, Graves and Ram Murty used the above growth results to prove Lenstra's result in some cases without assuming the GRH.  More precisely, if $K$ is a finite Galois extension of $\Q$ with ring of integers $R$, $\mathrm{rank}(R^*) >3$, and cyclic class group $\mathrm{Cl}_K$, such that its Hilbert class field $H(K)$ has an abelian Galois group over $\Q$, then $[C]$ generates $\mathrm{Cl}_K$ 
iff $[C]$ is a Euclidean ideal class \cite{GravesMurty}.  They then used this to show that $\Q(\sqrt{5}, \sqrt{21}, \sqrt{22})$ has 
a non-principal Euclidean ideal class that is not norm-Euclidean \cite{GravesMurty}.\\

\subsection{Miscellaneous}
\noindent We mention some other related results without delving deeply into the mathematics involved.\\

\noindent {\bf 1}) {\it How many primes $p$ are there with $p-1$ squarefree?}
There are two lines of attack.
One is \`a la Artin's primitive root conjecture.
In this approach one notes that $p-1$ is squarefree iff $p$
does not split completely in any of the
cyclotomic fields $\mathbb Q(\zeta_{q^2})$, where $q$ runs over all primes.
By inclusion-exclusion this gives then the conjectural
natural density
$\delta=\sum_{n=1}^{\infty}\mu(n)/[\mathbb Q(\zeta_{n^2}):\mathbb Q]$. The latter sum is easily seen to equal $A$.
An approach to this problem in this spirit was
made by Knobloch \cite{Hasse, knoflook}, a student of Hasse. It turns out, however, that
a quite elementary proof is possible. It was found by Mirsky \cite{Mirsky}. This uses the trivial identity
$\sum_{d^2|n}\mu(d)=1$ if $n$ is not squarefree and zero otherwise.
Note that the number of primes $p\le x$ such
that $p-1$ is squarefree equals
$\sum_{p\le x}\sum_{a^2|p-1}\mu(a)$. Now swap the order of summation.
The inner sum is a prime counting sum and can be estimated by the
Siegel-Walfisz theorem (\ref{siwa}) for $a$ small enough. The contribution from all large $a$ we
can trivially estimate. Carrying out this procedure gives us the
same sum for the density $\delta$ as found in the earlier approach.  
(Mirsky gives the details of another proof and then for the relevant
 one for our purposes says it can be done similarly. This is worked
 out, however, in more detail in, e.g., Moree and Hommersom \cite[pp. 17--20]{Huib}.) 
The same approach was recently followed 
by Broughan and Zhou \cite{BZhou} to show that the density of primes $p$ such that $p-1=2^em$ with $m$
odd and squarefree, equals $2A$. Clary and
Fabrykowski \cite{clear} determined the density of primes $p\equiv l({\rm mod~}k)$ such that $ap+b$ is squarefree, where
$(l,k)=1$ and $a>0$. It is an explicitly defined rational multiple of the Artin constant.\\
\indent In view of the analogy between the groups $\F_p^{*}$ and $E(\F_p)$, where $E/\F_p$ is
an elliptic curve, it is natural to also consider the squarefreeness of the elements of the  sequence
$\#E(\F_p) = p+1-a_p$, as $p$ goes to infinity. This problem may be interpreted as a stronger
version of that of calculating the number of primes $p$ for which the group $E(\F_p)$ is cyclic,
where $E$ is a global elliptic curve defined over $\Q$ and then reduced modulo primes $p$.
The problem of how often $p+1-a_p$ is squarefree was first considered by A.C. Cojocaru
in her PhD thesis \cite{Co3}, where she proved an unconditional asymptotic formula for these primes
if $E/\Q$ is with CM (see also \cite{Co5}) and a conditional asymptotic formula if $E/\Q$ is without CM.\\
\indent The above topic is related to the question of how many ways, $R(n)$ say, there are one 
can write a given number $n$ as a sum of a prime and a squarefree number. 
Let $B>1$ be fixed. In 1936 Walfisz \cite{Walvis} (see also 
\cite[pp. 386-387]{mova}) showed that
$R(n)=c(n){\rm Li}(n)+O(n\log^{-B}n)$, where
$$c(n)=A\prod_{q|n}\Big(1+{1\over q^2-q-1}\Big).$$
The appearance of $A$ in this context seems to be unrelated to the multiplicative order. It comes from
the first sum below, which happens to be equal to the second sum arising in primitive root theory:
$$c(n)=\sum_{d=1,~(d,n)=1}^{\infty}{\mu(d)\over \varphi(d^2)}=
\sum_{d=1,~(d,n)=1}^{\infty}{\mu(d)\over d\varphi(d)}.$$
For some material related to the result of Walfisz, see \cite{Ester, Page}.\\

\noindent {\bf 2}) {\it Higher rank Artin}. Note that $g$ is a primitive root mod $p$ if the subgroup generated
by $\bar g$, the reduction of $g$ mod $p$, equals $(\Z/p\Z)^*$.
If $g_1,\ldots,g_r$ are rational numbers, one can wonder about \\
a) the density of primes $p\le x$ such that
$\langle {\bar g_1},\ldots,{\bar g_r} \rangle$ is the full multiplicative
group $(\Z/p\Z)^*$;\\
b) the density
of primes $p$ such that each $g_i$, $1\le i\le r$, is a primitive root modulo $p$.\\
a) This was studied by Cangelmi and Pappalardi, see \cite{Cangelmi, Paplepel}.
The analogue of $A$ in this setting is the
$r$-rank Artin constant $$A_r=\prod_q \Big(1-{1\over q^{r}(q-1)}\Big).$$
Trivially $A_1=A$.\\
\indent Let $\kappa(p)$  denote the smallest integer $r$ such that the first $r$ primes
generate $(\mathbb Z/p\mathbb Z)^*$. Brown and Zassenhaus conjectured in \cite{BZ} that
$\kappa(p)\le [\log p]$ for almost all primes $p$ and that $\kappa(p)>\log p$ for infinitely many primes $p$. The latter
part of the assertion is true by the result of Graham and Ringrose \cite{graring} that the least quadratic non-residue 
exceeds $c (\log p )\log \log \log p$ for infinitely many primes $p$. \\
\indent If $U$ denotes the number of primes $p\le x$ for which $g(p)\ge T$, then, by (\ref{pappaheefthemnodig})
we have $UT\ll \sum_{p\le x}g(p)\ll \pi(x)\log^2 x(\log \log x)^4$. For any $\epsilon>0$, we choose
$T=\log^2 x(\log \log x)^{4+\epsilon/2}$ so that $U=o(\pi(x))$ and since $g(p)\le T$ is a product of primes
$\le T$, we infer that for almost  all primes $p$ we have the estimate
$\kappa(p)\le \log^2x(\log \log x)^{4+\epsilon/2}
\le (\log p)^2(\log \log p)^{4+\epsilon}$. 
Pappalardi \cite{paminimal} improved this by showing that $\kappa(p)=O(\log^2p/\log\log p)$ for
almost all primes $p\leq x$. Under GRH he showed in another paper \cite{Paplepel} that the
Brown-Zassenhaus conjecture is true. Moreover, he shows that
there is a positive absolute constant $B$ such that for every
divergent function $y=y(x)\le \log x/(4 \log \log x)$ and for all primes $p\le x$ with at most 
$O(\pi(x)B^{-y(x)})$ exceptions we have $\kappa(p)\le [y(p)]$. In 1993, Konyagin and Pomerance \cite{KoPo} and Pappalardi \cite{pathesis} independently proved that for
all $\epsilon>0$  and for all primes $p\le x$, $G(p)\le x^{\epsilon}$ with at most $O_{\epsilon}(1)$ exceptions.\\
b) This was studied by K. Matthews \cite{Matthews}, 
who determined the density under GRH (formula (1.4)). (This formula contains a typo, one
has to replace $q>2$ by $q\ge 2$.) 
The analogue of $A$ in this setting is 
$$A_r'=\prod_q \Big(1-{1\over q-1}\Big(1-\Big(1-{1\over q}\Big)^n\Big)\Big).$$
Trivially $A_1'=A$. Pappalardi \cite{panew} studied  the infinitude of this set of primes $p$ assuming 
Schinzel's Conjecture H. Schinzel \cite{schinzelrecent}, relying heavily on results from 
Matthews paper, determined assuming GRH the density of primes
$p$ such that two prescribed disjoint sets of odd primes $A$ and $B$ have the property that
$(b/p)=1$ for every $b\in B$ and every $a\in A$ is a primitive root modulo $p$. He used this to determine
the density $D(p_k,p_l)$ of primes $p$ such that the $k$th-prime $p_k$ is for $p$ the least quadratic non-residue and
$p_l$ the least prime primitive root modulo $p$.\\ 
\indent The shortest route known to calculating the densities is provided by the average-character-sum method, see
Moree and Stevenhagen \cite{Most}. In their setup Schinzel's result follows very easily.\\
\indent Again, this higher rank analogue of Artin's conjecture can be formulated in the context of
reductions of an elliptic curve $E/\Q$ with arithmetic rank $> 1$, as well as that of an elliptic curve $E/\F_q[T]$. 
As mentioned before, the rational situation was studied by Gupta and Ram Murty in 
\cite{GM1, GM2}, and the 
function field one by Hall and Voloch in \cite{HaVo}.\\

\noindent {\bf 3}) {\it Near-primitive roots}. Fix an integer $t$. For the primes $p\equiv 1({\rm mod~}t)$ we
can ask for the subset of them such that $\langle \bar g \rangle$ generates a
subgroup of index $t$ in $(\Z/p\Z)^*$, see \cite{cunningham, Lehmers, Lenstra, LMS, moeller, MoreeNP, 
Mojapan, Golombearly, MGolomb, murat-1, Wagstaff}. 
We let $N_{g,t}(x)$ denote
the corresponding counting function.
Assuming GRH the corresponding density, $\delta(g,t)$, can be shown to be
$$\sum_{n=1}^{\infty}{\mu(n)\over [\mathbb Q(\zeta_{nt},g^{1/nt}):\mathbb Q]}.$$ 
The first person to do so seems to have been H. M\"oller \cite{moeller}, followed a few years later
by Wagstaff \cite{Wagstaff}.
Since now $nt$ can be non-squarefree making the degree evaluation technically
more involved, the resulting answer is less elegant than in the case $t=1$. For $g>0$ or $g$ not
minus a square, the resulting density was given in Euler product form by Wagstaff \cite{Wagstaff}. For general
$g$ the density has been given in Euler product form by Moree \cite{MGolomb}. It turns out that the
average-character-sum method \cite{LMS} leads far more efficiently to this result. The latter method leads
also very directly to an iff result for $\delta(g,t)$ to vanish. Lenstra \cite{Lenstra} was the
first to give these vanishing conditions, but he did this without proof. Since Wagstaff did not work out the
Euler product form of the density for every $g$, his result did not make such a proof possible. In the
paper by Moree cited above, \cite{MGolomb}, a proof along the lines of Wagstaff's approach can be found.\\
\indent Moree introduced a function $w_{g,t}(p)\in \{0,1,2\}$ for which he proved 
(see \cite{MoreeNP}, for a rather easier reproof see \cite{Mojapan}) under GRH
that 
$$N_{g,t}(x)=(h,t)\sum_{p\le x,~p\equiv 1({\rm mod~}t)}w_{g,t}(p){\varphi((p-1)/t)\over p-1}
+O\Big({x\log \log x\over \log^2 x}\Big).$$
This function $w_{g,t}(p)$ has the property that, under GRH, $w_{g,t}(p)=0$ for
all primes $p$ sufficiently large iff $N_{g,t}$ is finite. Since the definition of
$w_{g,t}(p)$ involves nothing more than the Legendre symbol, it is then not difficult
to arrive at the vanishing conditions.\\
\indent To sum up, there are three different approaches to determine exactly when $\delta(g,t)$ vanishes:\\
1) completing Wagstaff's work and bringing $\delta(g,t)$ into Euler product form \cite{MGolomb};\\
2) using asymptotically exact heuristics for $N_{g,t}(x)$ \cite{MoreeNP, Mojapan};\\
3) using the average-character-sum method \cite{LMS}.\\
\indent There is no doubt though, that at present approach 3 yields the most elegant and short derivation of the vanishing conditions.\\ 
\indent Unaware of the above results Solomon Golomb made in 2004 \cite{CM} the following near-primitive root
conjecture (a review of near-primitive roots from the perspective of Golomb's conjecture is given in Moree \cite{Golombearly}. 
This perspective was abandoned in Moree \cite{MGolomb}, the final version).
\begin{conjecture}
For every squarefree integer $g>1$, and for every positive integer $t$, the 
set $N_{g,t}$ is infinite. Moreover, the density of such
primes equals a constant (expressible in terms of $g$ and $t$) times
the corresponding density for the case $t=1$ (Artin's conjecture).
\end{conjecture}
In a 2008 paper Franc and Ram Murty \cite{FM} made some progress towards
establishing this conjecture. In particular they prove the conjecture in 
case $g$ is even and $t$ is odd, assuming GRH. In general though, this conjecture
is false, since in case $g\equiv 1({\rm mod~}4)$, $t$ is odd and $g|t$, $N_{g,t}$ is
finite. To see this note that in this case we have $({g\over p})=1$ for the primes
$p\equiv 1({\rm mod~}t)$ by the law of quadratic reciprocity and thus $r_p(g)$ must
be even, contradicting the assumption $2\nmid t$.\\
\indent The author proposes the following conjecture (which on GRH can be shown to be true).
\begin{conjecture}
For every squarefree integer $g>1$, and for every positive integer $t$, the 
set $N_{g,t}$ is finite iff $g\equiv 1({\rm mod~}4)$, $t$ is odd and $g|t$. In case $N(g,t)$ is
infinite it has a natural density that equals a positive constant (expressible in terms of $g$ and $t$) times
the corresponding density for the case $t=1$.
\end{conjecture}
\indent Golomb apparently made his conjecture having in mind the problem of estimating the number of primes
in the set $S$, defined as follows. Let $\Phi_n(x)$ denote the $n$-th cyclotomic polynomial. Let $S$ be the set of primes $p$ such that
if $f(x)$ is any irreducible factor of $\Phi_p(x)$ over $\mathbb F_2$, then $f(x)$ does not divide
any trinomial. Using the explicit evaluation of $\delta(g,t)$ and 
results from Golomb and Lee \cite{gollee}, one can then deduce that, on GRH, the set $S$ has a subset
of density $>0.95$, see \cite{Golombearly}.\\
\indent Felix and Ram Murty \cite{Felix2, Felix1} are interested in the problem of proving estimates of the type
\begin{equation}
\label{adamfelix}
\sum_{p\le x}f(r_p(g))\sim A_f(g)\pi(x),
\end{equation}
The case $f(x)=\log x$ and $g=2$ was first studied by Bach et al.\,\cite{BLSW}. Fomenko \cite{Fomenko} showed that 
(\ref{adamfelix}) holds true in case $f(x)=\log x$, assuming GRH and Conjecture A of Hooley \cite[p. 112]{Hoolyboek}. Felix
and Ram Murty \cite{Felix1} managed to prove this result without assuming Conjecture A in case
$f(x)=\log^{\alpha}x$, with $0<\alpha<1$. In fact they can deal with a wider class
of functions including some arithmetic functions, e.g., $\omega$ and $\Omega$. In a 
sequel Felix \cite{Felix2} establishes (\ref{adamfelix})
with $r_p(g)$ replaced by its higher rank analogue. When the rank is 
at least two, a much more precise version of (\ref{adamfelix})
can be established (under GRH), namely with error term $O(x^{\theta})$ and $\theta<1$.\\
\indent Ram Murty and Simon Wong \cite{MurWong} showed that at least one of the following holds true:\\
1) $N_{2,2}$ is infinite.\\
2) There exist infinitely many primes $p$ such that $2^p-1$ is composite.\\
In the same spirit, they announce that a variation of their argument gives that if $2$ is not a primitive
root for infinitely many primes $p$, then $(2^p+1)/3$ is composite for infinitely many primes $p$.\\

\noindent {\bf 4}) {\it Primitive $\lambda$-roots}. The definition of primitive roots was extended by
Carmichael to
that of {\it primitive $\lambda$-roots} for composite moduli $n$, which
are integers with the maximal order modulo $n$ (which equals the
exponent, $\lambda(n)$, of the group $(\mathbb Z/n\mathbb Z)^*$). 
We have $\lambda(n)|\varphi(n)$. It is an exercise in elementary number theory to describe those $n$ for which
$\lambda(n)=\varphi(n)/2$, see, e.g., Lee et al.\,\cite{leeal}.
Let $N_g(x)$ be the
number of natural numbers up to $x$ for which $g$ is a primitive
$\lambda$-root. One might wonder then, in
analogy with Artin's conjecture, whether $N_g(x)\sim B(g)x$
for some positive constant $B(g)$ depending on $g$, perhaps with
some exceptional values of $g$. Li \cite{Li}
provides results suggesting that $N_g(x)/x$
does not typically tend to a limit.
Let $R(n)$ denote the number of residues modulo $n$ which are primitive $\lambda$-roots for $n$. 
One has $R(n)\ge \varphi (\varphi(n))$. M\"uller and Schlage-Puchta \cite{slag} showed that 
the set of $n$ with $R(n)=\varphi (\varphi(n))$ has density zero. 
Li \cite{Li}
showed that $(1/x)\sum_{n\le x}R(n)/n$ does not tend to a limit as $x$ tends to infinity. In particular,
one has
$$\lim_{x\rightarrow \infty} \sup {1\over x}\sum_{n\le x}{R(n)\over n}>0,~
\lim_{x\rightarrow \infty} \inf {1\over x}\sum_{n\le x}{R(n)\over n}=0.$$
Let $\mathcal E$ denote the set of integers $g$ which are a power higher than the first power or a square
times a member of $\{\pm 1,\pm 2\}$. Li \cite{Li-1} showed that for $g\in \mathcal E$ we have
$N_g(x)=o(x)$, and that for every integer $g$ one has $\lim \inf_{x\rightarrow \infty} N_g(x)/x=0$. Li
and Pomerance \cite{LiP2} showed that, under GRH, for each integer $g$ not in $\mathcal E$ one has $\lim \sup_{x\rightarrow \infty} N_g(x)/x>0$. Furthermore, Li \cite{Li2} showed that, for $y\ge \exp((\log x)^{3/4})$, one
has
$${1\over y}\sum_{1\le g\le y}N_g(y)\sim \sum_{n\le x}{R(n)\over n},~{\rm~as~}x\rightarrow \infty.$$
The range for $y$ was extended by Li and Pomerance \cite{LiP3}.
See also \cite{slis, Li2, Li3, slag} for further results and \cite{LiP} for a survey on primitive
roots, with special focus on primitive $\lambda$-roots.\\

\noindent {\bf 5}) {\it Squarefreeness}. Pappalardi \cite{Pappa} studied the problem of determining how many primes
$p\le x$ there are such that ord$_p(g)$ is squarefree (and similarly ord$_m(g)$, where
$m$ ranges over the integers $\le x$ for which the latter order is defined). Using a
simple characterisation of those $n$ for which Carmichael's function $\lambda(n)$ is
squarefree together with Mirsky's result from \S 9.7.1, Pappalardi, Saidak and Shparlinski \cite{PSS} have
shown that the number of $n\le x$ such that $\lambda(n)$ is squarefree is asymptotically
equal to $cx\log^{A-1}x$, for some constant $c>0$. Note that Artin's constant appears in
an unusual position in the formula! There are of course various generalisations possible,
e.g., to `roughly' squarefree values \cite{BL} and $k$th powerfree values \cite{BP}.\\ 

\noindent {\bf 6}) {\it Goldstein's conjecture and follow up}. 
Larry Goldstein formulated a conjecture in 1968 which implies Artin's
conjecture (announcement \cite{goldannounce}, detailed version \cite{gold-1}).
For each prime $q$, let $L_q$
be a Galois extension over $\Q$.  For each integer $k$ let $L_k$ be the 
compositum of those $L_p$ for which $p$ divides $k$, and 
let $n(k)$ be the degree of $L_k$ and $d_k$ its discriminant; then the set of primes $p$ that split completely in none of the 
$L_q$ equals
\begin{equation}
\label{murdi}
\sum_{k=1}^{\infty}{\mu(k)\over n(k)}
\end{equation} 
if this series converges absolutely. This conjecture is known to be false, even for certain abelian 
$L_q$ (Weinberg \cite{Wein}). In 1973, Goldstein \cite{gold} proved a special case of his 
conjecture (assuming GRH) for families of abelian 
$L_q$ satisfying certain conditions on their discriminants, degrees, and splitting of primes.
An
unconditional version of his result was established by Ram Murty \cite{MurtyGoldstein}. 
\begin{theorem}\label{thm:L_q}
Suppose that for $q$ sufficiently large, the extensions $L_q$ are abelian
over $\mathbb Q$ and
\begin{enumerate}
\item $\log |d_k|=O(n(k)\log k)$
\item if a prime $p$ splits completely in $L_q$,then for $q$ sufficiently large $p\ge f_q$, where $f_q$ is the conductor of $L_q$;
\item $\sum_{k\ge y}{1\over n(k)}=o({1\over \log y})$, as $y\rightarrow \infty$, 
\end{enumerate}
then the set of primes which do not split completely in
any $L_q$ has a density $\delta$ given by (\ref{murdi}).
\end{theorem}
In 1983, assuming GRH Ram Murty established a variant of Goldstein's conjecture in a lattice theoretic
setting and applied it to elliptic curves \cite{Mu1}.\\

\noindent {\bf 7}) {\it Brizolis' conjecture}. Brizolis conjectured that for every prime $p>3$ there exist an $x$ and a primitive root
$g({\rm mod~}p)$ such that $x\equiv g^x({\rm mod~}p)$. Let $N(p)$ denote the number of $1\le x\le p-1$ such that $x$ is a 
primitive root mod $p$ and $x$ and $p-1$ are coprime. It is easy to see that if $N(p)\ge 1$, then Brizolis' conjecture is 
true for the prime $p$. 
(Let $y$ be the multiplicative inverse of $x$ modulo $p-1$ and consider $g=x^y$. Then
$g^x\equiv x^{xy}\equiv x({\rm mod~}p)$ and $g$ is a primitive root for $p$.) 
Hausman \cite{Hausman} showed that 
$N(p)\ge 1$ for every sufficiently large prime $p$. Let $\omega(n)$ denote 
the number of distinct prime divisors of $n$ and 
 $d(n)$ the number of divisors of $n$. Cobeli and Zaharescu \cite{CZ} and, independently Wen Peng Zhang \cite{Zhang} 
using character sums and the Weil bounds proved that 
 $$N(p)={\varphi(p-1)^2\over p-1}+O(\sqrt{p}~ 4^{\omega(p-1)}\log p),$$ and used this to show that 
 Brizolis' conjecture is true for every prime $p>10^{50}$. The full Brizolis' conjecture was 
 settled by Campbell \cite{Campbell}. 
In a later paper \cite{LPSound} she and her coauthors using a numerically explicit ``smoothened'' version of the
P\'olya-Vinogradov inequality showed that for each prime $p>3$, there is a primitive root $g$ for $p$ in 
$[1,p-1]$ that is coprime to $p-1$. This led to a new proof of Brizolis conjecture requiring far less
computer calculation.\\
\indent Note that if $x\equiv g^x({\rm mod~}p)$, 
 then $x$ is 1-cycle under taking the discrete logarithm with respect to $g$. For small $k$, $k$-cycles of the
 discrete logarithm problem are studied in Holden and Moree \cite{HM2} (for more 
experimental work in this direction see \cite{Holden, HoldenLindle}), see also Glebsky and Shparlinski \cite{GSh}.
Let $F(p)$
denote the number of pairs $(g,h)$ satisfying $g^h\equiv h({\rm mod~}p)$ with $1\le g,h\le p-1$.
Let $G(p)$ denote the number of pairs $(g,h)$ satisfying $g^h\equiv h({\rm mod~}p)$ with $1\le g,h\le p-1$, where
in addition we require $g$ to be a primitive root modulo $p$.
Holden and Moree \cite[Conjecture 8.3]{HM2} conjectured that
$$\sum_{p\le x}{F(p)\over p-1}=(1+o(1))\pi(x){\rm ~and~}\sum_{p\le x}{G(p)\over p-1}=(A+o(1))\pi(x).$$
This conjecture with quantitative error estimates has meanwhile been established
by Bourgain et al.\,\cite{BKS}. A related conjecture by Holden and Moree \cite{HM2} states that one should have $F(p)=(1+o(1))p$. 
Bourgain et al.\,\cite{BKS} showed that $F(p)=p+O(p^{4/5+\epsilon})$ for a set of primes $p$ of relative density 1
and in a later paper, \cite{BKS2}, $F(p)\ge p+O(p^{3/4+\epsilon})$ and $F(p)=O(p)$.\\ 
\indent The above questions are all modulo $p$, we can also ask them modulo prime powers. A preliminary approach to these
questions was made by Holden and Robinson \cite{Holdrob} using $p$-adic methods, primarily Hensel's lemma and $p$-adic interpolation.\\

\noindent {\bf 8}) {\it $\omega(n)=\omega(n+1)$ and Artin}. Erd\H{o}s conjectured that there are infinitely many integers $n$ such that
$\Omega(n)=\Omega(n+1)$, where $\Omega(n)$ denotes the total number of prime divisors of $n$ and
a similar result for $\omega(n)$. The former conjecture was 
proved in 1984 by Heath-Brown \cite{opeen}, but the latter was only proved by Schlage-Puchta in
2003 \cite{Schlage}. 
Kan in 2004 \cite{Kan}, unaware of the latter
result, proved the weaker result that not both the $\omega(n)=\omega(n+1)$ conjecture
and the qualitative
form of Artin's primitive conjecture can be false.\\
\indent Heath-Brown actually proved that
there are infinitely many integers $n$ such that $d(n)=d(n+1)$ (thus solving a conjecture of Erd\H{o}s
and Mirsky) and remarks that his method can similarly prove that $\Omega(n)=\Omega(n+1)$ infinitely often, 
where $d(n)$ denotes the number of positive divisors of $n$.\\

\noindent {\bf 9}) {\it Artin over number fields}. Over arbitrary number fields, there are two obvious ways in which Artin's conjecture can be generalized. Fix a number 
field $K$ with ring of integers ${\mathcal{O}}_K$, and a nonzero element $\alpha\in 
{\mathcal{O}}_K$, which is not a root of unity. 
We expect the following generalization to hold: the set of 
primes $\mathfrak p$ of $K$ for which $\alpha$ generates the cyclic group $(\mathcal{O}_K/\mathfrak{p}\mathcal{O}_K)^*$ has a density 
inside the set of all primes of $K$. Moreover, the situation is highly similar to the rational case, as the set of 
primes $\mathfrak p$ of $K$ for which $(\mathcal{O}_K/\mathfrak p\mathcal{O}_K)^*$ is isomorphic to 
$(\mathbb Z/p\mathbb Z)^*$ has density 1. In 1975 Cooke and Weinberger \cite{CW} 
proved that this generalization of 
Artin's conjecture is indeed true and the density is given by the 
infinite sum in (\ref{classicalhooley}) with $\mathbb Q$ replaced by $K$, 
if an appropriate generalization of RH holds.\\
\indent Roskam \cite{Ros2} considers a generalization of Artin's conjecture in a ``rational" direction: the set of rational 
primes $p$ for which the order of $\alpha$ in $(\mathcal{O}_K/p\mathcal{O}_K)^*$ is equal to the exponent of 
this group has a 
density inside the set of all rational primes. He proves that, assuming that an appropriate generalization of RH 
is true, this conjecture holds for quadratic fields $K$. A less general result in this direction he established
in \cite{Ros1} (but with a rather easier proof). Kitaoka and various of his pupils 
did a lot of work on Artin's conjecture for units (rather than
general elements) in number fields, cf. \cite{CKY, I, Kat-1, Kat, kit, kit1, kit2, kit3} (this being
only a partial list of references).\\

\noindent {\bf 10}) {\it Artin for matrices}. Let $A$ be a hyperbolic matrix in SL$_2(\mathbb Z)$, that is a matrix with absolute trace greater 
than 2. We may define ord$_n(A)$ to be the order of $A$ in ${\rm SL}_{2}(\mathbb Z/n\mathbb Z)$. It can be shown
that for a hyperbolic matrix in SL$(2,\mathbb Z)$ there is a density-one sequence of integers $n$ such that
its order mod $n$ is slightly larger than $\sqrt{n}$. This is a crucial ingredient in an ergodicity theorem 
regarding eigenstates of quantized linear maps $A$ of the torus (`cat maps'), see Kurlberg and Rudnick 
\cite{KR}. Under an appropriate generalisation of RH it can be shown \cite{K1} that for almost all $n$ we have
ord$_n(A)\gg n^{1-\epsilon}$ for any $\epsilon>0$. Call a matrix $A$ exceptional if it is of finite order or if it is diagonalizable and a 
power $A^r$ of $A$ has all eigenvalues equal to powers of a single rational integer $>1$, or equal to powers 
of a single unit $\not=±1$ of a real quadratic field. Corvaja et al.\,\cite{CRZ} established that if $A$ is not 
exceptional, then the quotient ${\rm ord}_n(A)/\log n$ tends to infinity with $N$, 
with their proof ultimately relying on the Schmidt subspace theorem 
from Diophantine approximation. Their result is best possible since, for any exceptional matrix $A$, there 
exist some $c>0$ and arbitrarily large integers $n$ for which ${\rm ord}_n(A)<c\log n$.\\ 
\indent Since the eigenvalues of $A$ are units in a quadratic number
field, this problem is closely related to that discussed in part 9.\\

\noindent {\bf 11)} {\it Artin for K-theory of number fields}. (Written by W. Gajda.) 
For basic references concerning K-theory the reader should consult 
\cite{Ga1} and \cite{Ga2}. For a number field $K$ we denote by 
${\O}_K$ the ring of algebraic integers. The class group of ${\O}_K$ 
is isomorphic to the quotient $K_0({\O}_K)/\Z,$ where $K_0(R)$ denotes 
the Grothendieck of a commutative ring $R$ with unity. The group of units 
$({\O}_K)^{\star}$ coincides with the group $K_1({\O}_K)$ defined by 
Bass in the sixties. In 1972 Quillen introduced groups $K_n(R),$ for all 
integers $n\geq 0,$ {\sl the} so called {\it higher algebraic K-theory groups of $R.$} In the case 
$R={\O}_K$ one can treat Quillen's groups  as higher dimensional analogs of the class group of 
${\O}_K$ - for $n$ even, and of the group of units $({\O}_K)^{\star}$
 - for $n$ odd. The groups $K_n(R)$ enjoy many very useful properties. In particular,  
they depend functorially on the ring $R.$ For a fixed $n,$ the reduction at 
a prime ideal $\mathfrak p$ of ${\O}_K$ induces a group homomorphism: 
$K_n({\O}_K)\longrightarrow K_n({\O}_K/\mathfrak p),$ which for $n=1$ 
can be easily identified with the map 
$({\O}_K)^{\star}\longrightarrow ({\O}_K/\mathfrak p)^{\star}.$ Due to 
a classical result of Quillen $K_n({\O}_K)$ are finitely generated 
abelian groups. On the other side, $K_n({\O}_K/\mathfrak p)$ vanishes, for $n{>}0$ and even,  
and is a finite cyclic group, if $n{>}0$ and odd. The ranks of the groups $K_n({\O}_K)$ were determined 
by Borel. In particular, if ${\O}_K{=}\Z,$ then the group $K_n({\O}_K)$ is of rank 
1, if $n{=}0$ or $n{=}4m+1$ for an integer $m{>}0,$ and is finite, otherwise. We expect that the following 
analog of the Artin conjecture holds for the K-groups of the integers: {\it 
for a fixed $n{=}4m+1$ and $m{>}0,$ the set of primes $p$ for which the homomorphism 
$K_n(\Z)\longrightarrow K_n(\F_p)$ is onto, has a natural positive density.} Note that for $n$ 
as above, due to another classical result of Quillen $K_n(\F_p)\simeq \Z/(p^{2m{+}1}{-}1),$ and conjecturally 
$K_n(\Z)\simeq \Z,$ at least up to 2-torsion. 
For some results supporting this conjecture and the discussion of the relation of the reduction 
map on the K-groups to deep questions in number theory, such as the 
Kummer-Vandiver conjecture, we refer the reader to the survey paper 
\cite{Ga2}.\\

\noindent {\bf 12)} {\it Griffin's dream}. R. Griffin (a mathematical amateur) thought in 1957 that the decimal expansions of $1/p$ should have
period length $p-1$ for all primes of the form $10X^2+7$. The first 16 primes $p$ of the form $10X^2+7$ have
indeed decimal period $p-1$, but this is not true for $p=7297$, the 17th such prime. D. Lehmer \cite{larti} found in 1963 that
326 is a primitive root for the first 206 primes of the form $326X^2+3$. More impressive examples in the
same spirit can be given using recent results on prime producing quadratics by Moree \cite{Moreegriffin} 
and K. Scholten \cite{Scholten}. Y. Gallot holds the
record in which 206 is being replaced by 38639. 
Scholten \cite{Scholten} was the first to consider prime producing cubics and he found one having the property that the first 10011 primes $p$ it 
produces are such that that 11045 is a primitive root modulo $p$.\\

\noindent {\bf 13}) {\it Artin and binomial coefficients}. Interest in prime divisors of binomial coefficients dates back at least
to Chebyshev. It is not difficult to see that unless $k$ is fairly close to one of the ends of the 
range $0\le k\le n$, then the coefficient ${n\choose k}$ contains many prime factors and many of these occur with 
high multiplicity. Erd\H{o}s made a number of conjectures quantifying these observations, one well-known conjecture being 
that the middle coefficient ${2n\choose n}$ is not squarefree for any $n>4$. Granville and Ramar\'e \cite{GR}
proved this conjecture in 1996. The corresponding problem for
$q$-binomial and, more generally, $q$-multinomial  coefficients is much more complicated. Sander \cite{Sander} 
showed how Artin's primitive root conjecture implies the respective results.\\

\noindent {\bf 14)} {\it Ulmer's rank result}. Let $\mathbb F_q$ be the finite field of $q$ elements of prime characteristic $p$. 
If $(n,q)=1$, then it can be shown that $X^n-1$ factors into
$$i_q(n)=\sum_{r|n}{\varphi(r)\over {\rm ord}_q(r)}$$
distinct irreducible factors over $\mathbb F_q$. 
The quantity $i_q(n)$ appears in various mathematical
settings, for example in Ulmer's study of elliptic curves with large rank over function
fields \cite{Ulmer}, the study of the cycle structure
of repeated exponentation modulo a prime $p$ \cite{choushpa, Sha}, 
see Moree and Sol\'e \cite{Moso} for further examples). Using the Cauchy-Frobenius 
identity coupled with elementary group-theoretical considerations, Deaconescu \cite{diaconie} has shown
that 
$$i_q(n)={1\over {\rm ord}_q(n)}\sum_{d|{\rm ord}_q(n)}\varphi\Big({{\rm ord}_q(n)\over d}\Big)(n,q^d-1).$$
Ulmer considers the parametric
family of curves $E_d:~y^2+xy=x^3-t^d$ over the function field $\mathbb F_q(t)$, where $d$ is a positive
integer. Ulmer \cite[Theorem 9.2]{Ulmer} showed that if $d$ is a divisor of the sequence $\{p^k+1\}$, then the
rank $R_q(d)$ of $E_d$ over $\mathbb F_q(t)$ is given by $R_q(d)=i_q(d)-\epsilon_q(d)$, with
$\epsilon_q(d)$ an explicit correction term that always satisfies $0\le \epsilon_q(d)\le 4$.\\ 
\indent Using techniques from the study of Artin primitive root
type problems (for example a version of the Chebotarev density theorem due to Lagarias 
and Odlyzko), Pomerance and Shparlinski \cite{PoS} showed that the typical rank in Ulmer's family $R_q(d)$ tends
to infinity as $d\rightarrow \infty$. 
In particular, they proved that except for $o_{p,\epsilon}(x)$ values of $d\le x$, one has
$R_q(d)\ge (\log d)^{(1/3-\epsilon)\log \log \log d}$. Further they established the existence of an absolute constant
$\alpha>1/2$ such that for all finite fields $\mathbb F_q$ and all sufficiently large values of $x$ (depending
only on the characteristic $p$ of $\mathbb F_q$), one has
$x^{-1}\sum_{d\le x}R_q(d)\ge x^{\alpha}$.
This shows that the family $E_d$ is quite special, as
Brumer \cite{brumer} has shown that the average rank for all elliptic curves over a function field of positive
characteristic is asymptotically bounded above by $2.3$.\\

\noindent {\bf 15)} {\it Primitive roots compared with
quadratic non-residues}. Let $n$ be a power $p^h$ or $2p^h$, with $p$ an odd prime. 
Write $$PR(n)=\{g\in(\mathbb Z/n\mathbb Z)^\times,\; g {\rm ~is~a~primitive~root~modulo~} n\},$$
$$QNR(n)=\{a\in(\mathbb Z/n\mathbb Z)^\times,\;a{\rm ~is~a~quadratic~non{\rm -}residue~modulo~} n\}$$
for the set of primitive roots and quadratic non-residues modulo $n$, respectively. 
Obviously $PR(n)\subseteq QNR(n)$. 
By methods from elementary number theory and combinatorics Gun
et al.\,\cite{gun} classified those $n$ for which $$\#(PR(n))=\#(QNR(n))-2^r.$$
Modulo $p$ the number of nonquadratic residues which are not primitive roots is
obviously $g(p):=(p-1)/2-\varphi(p-1)$. The values assumed by $g$ were
investigated by Luca and Walsh \cite{LW} and by Robbins \cite{Robbins}. More recently, Gun et al.\,\cite{gun2}
applied character sum estimates to prove results on consecutive quadratic non-residues
modulo $p$ that are not primitive roots. They showed for example that given
a fixed $0<\epsilon<1/2$ and any positive integer $N$, then for all primes
$p\ge \exp((2\epsilon^{-1})^{8N})$ satisfying $\varphi(p-1)/(p-1)\le 1/2-\epsilon$,
we can find $N$ consecutive quadratic non-residues
modulo $p$ that are not primitive roots. 
More recently Luca et al.\,\cite{lucca2} showed that in fact the same property holds starting with
significantly smaller primes. 
These results are analogous to earlier
results by A. Brauer \cite{braueroud} (quadratic residues, quadratic non-residues)
and Szalay \cite{Szalay} (primitive roots). \\

\noindent {\bf 16)} {\it The Crandall-Collatz $qx+1$ problem}. For positive odd numbers $q$ and $m$, let
$C_q(m)$ denote the largest odd factor of $qm+1$. The sequences of iterates $C_q(m),~C_q(C_q(m)),\ldots$ consists
of the odd numbers in the orbit of $m$ under the $qx+1$ function. The $3x+1$ conjecture (unproven) asserts that for $q=3$ such an orbit always contains the number 1, cf. Lagarias \cite{3x+1} or the book \cite{3x+1book}. Crandall conjectured that
for every odd number $q>3$ there always exists an integer $m$ such that 1 does not occur in the orbit
of $m$ under the $qx+1$ function. A number $q$ for which such an $m$ exists we call a {\it Crandall number}.
Crandall \cite{Crandall} showed that $5,181$ and $1093$ are Crandall numbers. The proofs for 
$5$ and $181$ involve cycles: for $q=5$, we have the cycle $13,33,83,13$, and for $q=181$, we have the
cycle $27,611,27$. The proof for $1093$ depends on the fact that it is a Wieferich prime, that is a prime
$p$ such that $2^{p-1}\equiv 1({\rm mod~}p^2)$. It is conjectured that there are infinitely many, but
presently only two of them are known: 1093 and 3511. We define a positive odd integer $q$ to be 
a {\it Wieferich number} if $q$ and $(2^{{\rm ord}_q(2)}-1)/q$ are not coprime. These definitions are consistent in the sense that a Wieferich prime is a Wieferich number, and a prime Wieferich number is a
Wieferich prime. Franco and Pomerance \cite{FP} showed that every Wieferich number is a Crandall number, and
that the Wieferich numbers have relative density 1 in the odd numbers. The appearance of ord$_q(2)$ in
this setting is trivial: the odd number $r$ is in the range of $C_q$ iff $2^jr\equiv 1({\rm mod~}q)$ for some
integer $j$, so that the residue classes modulo $q$ in the range of $C_q$ are precisely those in the subgroup
of $(\mathbb Z/q\mathbb Z)^*$ generated by 2 modulo $q$ (see also \cite{KZ}). Also in some weaker versions of the $3x+1$ problem, 
the order of $2$ modulo $q$ makes its appearance, see, e.g., \cite{Caraiani}.\\
\indent As a byproduct Franco and Pomerance \cite{FP} showed that $N_m(x)=o(x)$ (see \S 9.7.17 for the definition of $N_m(x)$).\\

\noindent {\bf 17)} {\it Counting integers with odd order}. Let $m\ge 1$ be an arbitrary integer and let $N_m(x)$ count the number of odd integers $u\le x$
such that $m\nmid {\rm ord}_u(2)$. In case $m=q$ is an odd prime H. M\"uller \cite{HeM1} proved that
$${x\over \log^{1/(q-1)}x} \ll N_q(x)\ll {x\over \log^{1/q}x}.$$
This was improved in Moree \cite{pmuller2} to
$$N_q(x)=c_q{x\over \log^{q/(q^2-1)}x}\Big(1+O\Big({(\log \log x)^5\over \log x}\Big)\Big),$$
with $c_q>0$ a constant. The latter result was subsequently generalized by M\"uller \cite{HeM2} to
the case where $m=q^n$. Subsequently this result was generalized to 
arbitrary $m$ see Moree \cite{pmuller1}. For a nice survey of the material discussed under \S 9.7.16 and \S 9.7.17, 
see M\"uller \cite{HeM3}.\\

\noindent {\bf 18)} {\it Pseudorandom number generators}. Of importance in computer science are pseudorandom number generators, an example being
the power generator, where the sequence is given by $u_{i+1}\equiv u_i^e({\rm mod~}n)$. This
generator is periodic, and for it to behave pseudorandom, the period must be large. Some light on this
can be shed by techniques used in the study of Artin's primitive root conjecture; see, e.g., \cite{FPS, KP, MP}. Goubin, Mauduit and S\'ark\H{o}zy \cite{GMS} studied the pseudorandomness
of sequences $e_1,\ldots,e_N$ with $e_n=({f(n)\over p})$ if $p$ and $f(n)$ are coprime, and
$e_n=1$ otherwise, where $f(X)\in \mathbb F_p[X]$ has positive degree. In their
considerations the order of $2({\rm mod~}p)$ plays an important role (the larger it is
the better).\\ 
\indent The analogous  problem of studying
pseudoprime reductions of an elliptic curve over $\Q$ was considered by Cojocaru, Luca
and Shparlinski in \cite{CoLuSh}. Some of their results were recently improved on by David
and Wu \cite{DavidWu}.\\ 

\noindent {\bf 19)}  {\it Non-primitive roots}. Let $E(M,N)$ denote the number of integers in the interval $[M+1,M+N]$ which are not
primitive roots modulo $p$ for any odd prime $p\le \sqrt{N}$. Gallagher, in the fundamental paper \cite{Ga} in which
he introduced what is now called Gallagher's large sieve, gave the estimate
$E(M,N)=O(\sqrt{N}\log N)$. This was improved by Vaughan \cite{va} to $E(M,N)=O(\sqrt{N}\log^c N)$ for some
explicit $c=0.37\ldots$. Note that the result is quite sharp as $E(M,N)$ includes the squares
in the interval and so trivially $E(0,N)\gg \sqrt{N}$.
A number field analogue of Gallagher's result was established by Hinz \cite{hinz1a}. 
In a later paper Hinz \cite{hinz2} obtained an analogue of Vaughan's result with an explicit value of $c<1$ depending on 
the field.\\

\noindent {\bf 20)} {\it Romanoff's result and follow up}. In Bilharz's \cite{Bilharz} work on the Artin primitive root
conjecture over $\mathbb F_q[T]$ the density 
\begin{equation}
\label{bilhoi}
\sum_{m=1,~q\nmid m}^{\infty}{\mu(m)\over m~{\rm ord}_m(q)}
\end{equation}
appeared. The absolute convergence of this sum was proven by Romanoff \cite{Romanoff}. A simpler proof
was given by Erd\H{o}s and Tur\'an \cite{ET}. However, a much simpler proof of this was given by 
Erd\H{o}s \cite{romer}.
Romanoff's result
was strengthened and generalized by Landau \cite{Landau} and more recently by Ram Murty et
al.\,\cite{MRS}, who also established analogous results over number fields and for abelian varieties.\\
\indent Romanoff proved that the integers of the form $p+2^k$ have positive lower density (for 
a nice reproof see Montgomery and Vaughan \cite[pp. 98-99]{mova}).
He raised
the following question: does there exist an arithmetic progression consisting only of odd numbers, no
term of which is of the form $p+2^k$? Erd\H{o}s \cite{romer-1} found such an arithmetic progression by considering integers
which are congruent to $172677$ modulo $5592405(=(2^{24}-1)/3)$. Thus the density of numbers of the form 
$p+2^k$ is less than $1/2$, the trivial bound obtained from the the odd integers.
Put
$${\underline d}=\lim \inf_{x\rightarrow\infty}{\#\{p+2^k\le x\}\over x/2}{\rm ~and~}
{\overline d}=\lim \sup_{x\rightarrow\infty}{\#\{p+2^k\le x\}\over x/2}.$$
Habsieger and Roblot \cite{HR} showed that ${\overline d}<0.9819$.  Pintz \cite{Pintz} showed
that ${\underline d}\ge 0.18734$. The first explicit lower bound for $\underline d$ was found by
Chen and Sun \cite{CS}.\\
\indent A related conjecture is due to Erd\H{o}s, who conjectured that $7,15,21,45,75$ and $105$ are the
only integers $n$ for which $n-2^k$ is a prime for every $k$ with $2\le 2^k<n$. Vaughan \cite{va} gave
an estimate for the number of integers $n\le x$ of this form \cite[(1.5)]{va}. Assuming that $2$ is a primitive
root infinitely often, it is not difficult (see Hooley \cite[pp. 113--115]{Hoolyboek}) to show that the number of such 
integers $n\le x$ is $O(x^{1-A+\epsilon})$. Narkiewicz \cite{nar-1} sharpened this estimate to
$O(x^{1-A/(\log 2)+\epsilon})$.\\
\indent Ram Murty and Srinivasan \cite{MSri} considered a sum akin to (\ref{bilhoi}). They showed that if
$$S_a(x):=\sum_{p\le x,~p\nmid a}{1\over {\rm ord}_p(a)}=O(x^{1/4}),$$ then Artin's conjecture holds true. Unconditionally they
showed that the latter sum is $\ll \sqrt{x}$. Erd\H{o}s and Ram Murty \cite{EM} later improved on this by showing there exists
$\delta>0$ such that $S_a(x)\ll \sqrt{x}\log^{-1-\delta}x$. Felix \cite{Felix} showed that if $x/\log x=o(y)$, then
$${1\over y}\sum_{a\le y}S_a(x)=\log x+O(\log \log x)+O({x\over y}),$$
showing that $S_a(x)$ equals $\log x$ on average. He established a similar result for $S_a(x)$, but with
ord$_n(a)$ replaced by $\varphi({\rm ord}_n(a))$, showing that this sum is $\zeta(2)\zeta(3)\log x/\zeta(6)$ on average. The
asymptotic he obtains involves a constant arising in the Titchmarsh divisor problem (cf. Problem 40).\\

\noindent {\bf 21)} {\it Erd\H{o}s-Kac type theorems for the multiplicative order}.
An arithmetic function is said to have a {\it normal order} $F(n)$, iff for any $\epsilon>0$,
for almost all $n\le x$, one has $(1-\epsilon)F(n)<f(n)<(1+\epsilon)F(n)$. In 1917 Hardy and Ramanujan 
\cite{HardyRamanujan} showed that $\omega(n)$ has normal order $\log \log n$. 
In 1940, Erd\H{o}s and Kac \cite{EK}
proved a more refined result: the quantity 
$${\omega(n)-\log \log n\over \sqrt{\log \log n}}$$
is normally distributed, thus the normal order $\log \log n$ serves as the mean, and
$\sqrt{\log \log n}$ as standard deviation. This result can be more precisely formulated as
follows:
$$\lim_{x\rightarrow \infty}\Big({1\over x}\cdot \#\Big\{n\le x:\alpha\le {\omega(n)-\log \log n\over 
\sqrt{\log \log n}}\le
\beta\Big\}\Big)=G(\alpha,\beta),$$
where $G(a,b):={1\over \sqrt{2\pi}}\int_a^b e^{-t^2/2}dt$ denotes the Gaussian normal distribution. 
This celebrated result marked the birth of a whole new branch of number theory: probabilistic number
theory. See Elliott \cite{Prob1, Prob2} and Tenenbaum \cite{tenenbaum} for books dealing with this topic.

Under GRH the following result was proved by Saidak in his PhD thesis and Ram Murty and Saidak \cite{MurS}. 
Let $a\ge 2$ be an integer. We have
$$\lim_{x\rightarrow \infty}\Big({1\over x}\cdot \#\Big\{n\le x:(a,n)=1, \alpha\le {\omega({\rm ord}_n(a))-{1\over 2}(\log
\log n)^2\over 
\sqrt{{1\over 3}(\log \log n)^3}}\le
\beta\Big\}\Big)=G(\alpha,\beta){\varphi(a)\over a}.$$
Later Li and Pomerance \cite{LiP2} gave a different proof, also requiring GRH.
Erd\H{o}s and Pomerance \cite{EP} had earlier conjectured the latter result and proved
that it is true with ord$_n(a)$ replaced by $\varphi(n)$.
For analogues in the setting of elliptic curves, see, e.g., Cojocaru \cite{Co4} and 
Liu \cite{Liu1, Liu2},
and in the context of Drinfeld modules see \cite{Co6, KuoL-1, KuoL, KuoL2}.\\
\indent In \cite{KuoL-1} Kuo and Liu give a short survey of Erd\H{o}s-Kac type theorems and sketch
the proof of the main result in their paper \cite{KuoL}.\\

\noindent {\bf 22)} {\it Gauss periods}. In connection with a problem involving Gauss periods, Gao et al.\,\cite{GGP} raised the question of
determining the density of primes $p\equiv 1({\rm mod~}k)$ such that ${p-1\over {\rm ord}_p(g)}$
and ${p-1\over k}$ are coprime. For a given $g$ and $k$, von zur Gathen and Pappalardi \cite{GP} showed that, 
under GRH, this density exists and computed an Euler product for it. Gauss periods are of interest to perform fast
arithmetic in finite fields \cite{fasta}.\\

\noindent {\bf 23)} {\it Golomb's conjectures}. In 1984, Golomb \cite{gol} made four conjectures, which if true 
could be applied to construct so-called Costas arrays which were first
considered by Costas \cite{cost} in attempting to construct sonar signal patterns. 
Independently, Gilbert \cite{Gilbert} also wrote about them in the same year, publishing what is now known as the logarithmic Welch method of constructing Costas arrays.
Costas arrays concern primitive elements in finite fields. For
our purposes it suffices to formulate the first three of these conjectures in the case where the finite
field is of prime order.\\
A) If $p$ is odd, then there are two not necessarily distinct primitive roots $g_1$ and $g_2$ such that
$g_1+g_2\equiv 1({\rm mod~}p)$.\\ 
B) If $p$ is odd, then there are two not necessarily distinct, primitive roots $g_1$ and $g_2$ such that
$g_1+g_2\equiv -1({\rm mod~}p)$.\\ 
C) For all $p$ large enough, every $c$ with $p\nmid c$ can be written as $c\equiv g_1+g_2$ with
$g_1$ and $g_2$ primitive roots modulo $p$.\\

\noindent After work of various authors, it is now known that these conjectures are true;
see S. Cohen and Mullen \cite{CoMu}. A more
general form of Conjecture C is also known to be true \cite{cohie2}. A related question is to estimate
the number of solutions $g_1,g_2$ for a fixed $c$ (representations
of $c$) asymptotically as $p\rightarrow \infty$ \cite{cohsum, Zhang2}.\\

\noindent {\bf 24)} {\it Fibonacci primitive roots}. Various authors considered primitive roots that in addition are required to satisfy
a polynomial equation. In this context, most frequently the {\it Fibonacci primitive roots} have been
considered: a primitive root $g$ modulo a prime $p$ is called a Fibonacci primitive root 
if $g^2\equiv g+1 ({\rm mod~}p)$. More
generally, if $m>1$, $A$ and $B$ are integers, then a primitive root $g$ 
modulo $m$ is called a generalized Lucas primitive root modulo $m$ with parameters $A$ 
and $B$ if $m$ divides $g^2-Ag+B$, see, e.g., Mollin \cite{Mollin}.\\
\indent Shanks \cite{Shanks} conjectured that the density of Fibonacci primitive roots equals
$27A/38$. Under GRH this was shown to be true by Sander \cite{Sander1} in 1990. However, this
result was already established by Lenstra in 1977 \cite[Theorem 8.1]{Lenstra}.\\

\noindent {\bf 25)} $N(H,p)$. For primes $p$ let $N(H,p)$ denote the number of primes $q \le H$ which are primitive 
roots modulo $p$. Heuristically one expects that  $N(H,p)$ should be well
approximated by ${\phi(p-1) \over p-1}\pi(H)$. Results in this direction were established by
Elliott \cite{Elliott} and, more recently,  Nongkynrih \cite{non}.\\

\noindent {\bf 26)} {\it Primitive roots and groups}. There are many papers in the group theory literature
based on the problem of recognizing a finite group by the set of its 
element orders. For a finite group the set of orders of all its elements is called a spectrum of 
this group. A finite group is called recognizable from its spectrum if all finite groups with the 
same spectrum are isomorphic to this group. 
 Lyuchido and Mogkhaddamfar \cite{LM}
proved that for an odd prime $p$ the projective 
special linear group $L_p(2)$ is recognizable from its spectrum if $2$ is a primitive root 
modulo $p$.\\ 
\indent Recall 
that ${\rm Aut}(C_n(p))$ is the symplectic group of 
$2n\times 2n$ matrices over $\mathbb F_p$. Thompson \cite{thom} proved that 
if $p$ is an odd prime which is a primitive root modulo 
the prime $l\geq5$, then ${\rm Aut}(C_{(l-1)/2}(p))$ is the Galois group of an extension of $\mathbb Q(t)$.\\
\indent Let $p$ and $r$ be odd primes. Ab\'ert and Babai \cite{AB} considered the wreath product
$W(r,p)$ of the cyclic groups $C_p$ and $C_r$. They showed that the maximum size of the
minimal generating set of $W(r,p)$ equals $2+(p-1)/{\rm ord}_p(g)$. They used this result
to answer a question of Lubotzky.\\
\indent Martin and Valette \cite{MV} considered the solvable Baumslag-Solitar group
$BS_n=\{s,b:aba^{-1}=b^n\}$ for $n\ge 2$. They showed for example that if the Artin conjecture holds for the
integer $n$, then the spectrum of the associated Markov operators contains the set $\{z\in \mathbb C:|z|=1/2\}$.\\

\noindent {\bf 27)} {\it Least primitive root modulo $p$}. Let $g(p)$ denote the least primitive root modulo $p$, 
and $h(p)$ be the least positive primitive root modulo $p^2$. 
Using Gauss sums Vinogradov \cite{vino2} 
(see also Landau \cite{Landauvor})
showed in 1930 that $g(p)< 2^{\omega(p-1)}p^{1/2}\log p$.
This was improved to $g(p)<2^{\omega(p-1)+1}\sqrt{p}$ in 1942 by
Hua \cite{Hua}. Erd\H{o}s \cite{romer-2} showed that $g(p)=O(\sqrt{p}\log^{17} p)$, and 
Erd\H{o}s and Shapiro \cite{ES} that  $g(p)=O(\omega(p-1)^c \sqrt{p})$, for some absolute constant $c$.
Burgess \cite{Burgess}
using his celebrated character sum estimates proved that
$g(p)=O(p^{1/4+\epsilon})$ and $h(p)=O(p^{1/2+\epsilon})$. 
The former estimate was achieved independently by Wang \cite{wangie}.
S. Cohen et al.\,\cite{COS} 
sharpened the latter estimate to $h(p)=O(p^{1/4+\epsilon})$.
Elliott and Murata \cite{pdta2} considered the least primitive root mod $2p^2$. They
also considered moments of $g(p)$, that is sums of the form
$\pi(x)^{-1}\sum_{p\le x}g(p)^{\delta}$ (in \cite{pdta3}); in particular they showed that
if $\delta<1/2$ assuming GRH the latter sum tends to a constant depending on $\delta$. 
Elliott \cite{Elliott} has shown that for
all but $O(X^{\epsilon})$ primes $p$ up to $X$, one has $g(p)=O(\log^{c_{\epsilon}}p)$, with
$c_{\epsilon}$ a constant depending on $\epsilon$. In the same paper he proved that
$g(p)<475\log^{1.6}p$ for infinitely many primes $p$. 
The latter result was superseded in 1984 by Gupta and Murty \cite{GM} who obtained the bound
$g(p)<2250$ for infinitely many primes $p$. Subsequently this bound was reduced to $g(p)\le 7$
for infinitely many primes $p$ by Heath-Brown \cite{HB}.
Under GRH, Wang \cite{wangie} has
shown that $g(p)\ll (\log^2 p )\omega(p-1)^6$. 
Utilizing a combinatorial sieve due to Iwaniec, Shoup \cite{shoup} improved this to
$g(p)\ll \omega(p-1)^4(\log (\omega(p-1))+1)^4 \log^2 p$. 
Towards a lower estimate for $g(p)$, it
is known that infinitely often $g(p)$ can exceed $c(\log p)(\log \log \log p)$ for 
a certain $c>0$ \cite{graring}. 
This improves dramatically on a result of
Pillai \cite{pil}, who showed that unconditionally $g(p)=\Omega(\log \log p)$.
Assuming GRH one has $g(p)>c(\log p)(\log \log p)$ for infinitely many primes $p$ (see Montgomery \cite[Theorem 13.5]{monties}).
An
often studied quantity in the literature is $n_2(p)$, which is defined as the least quadratic
non-residue modulo $p$. 
For example in 1952, Ankeny \cite{Ankeny} proved the important result that $n_2(p)=O(\log^2p)$ on GRH. This consequence
of GRH has the important corolloray of makng Miller's primality test polynomial time (cf. Bach \cite{bachdiss}).
Since clearly $g(p)\ge n_2(p)$, any $\Omega$-result for $n_2(p)$ will
imply the same result for $g(p)$.\\ 
\indent Burgess and Elliott \cite{BE} have shown that
\begin{equation}
\label{pappaheefthemnodig}
{1\over \pi(x)}\sum_{p\le x}g(p)\ll \log^2 x(\log \log x)^4,
\end{equation}
improving on an earlier bound $x\log^A x$ (with an unspecified $A$) of Burgess \cite{burgessconf}.
Cohen et al.\,\cite{COS} showed that in (\ref{pappaheefthemnodig}) one can replace $g(p)$ by $h(p)$, thus
improving on an earlier result of Burgess \cite{burger} to the effect that $\pi(x)^{-1}\sum_{p\le x}h(p)\ll \log^3 x(\log \log x)^6$.
Murata \cite{murat} has shown that  $\pi(x)^{-1}\sum_{p\le x}g(p)\ll \log x(\log \log x)^7$, assuming
the GRH, meanwhile the exponent 7 has been replaced
by any number $>4$ by Elliott and Murata \cite{pdta3}.\\
\indent Let $D(i)$ be the density of primes $p$ such that $g(p)=i$. Elliott and Murata \cite{pdta3} 
proved under GRH that $D(i)=\sum_{M}(-1)^{|M|-1}A_M$, where $M$ runs over all the subsets of the
set $\{1,2,\ldots,i\}$ containing $i$, with $A_M$ the density (under GRH) of primes $p$ such that
each $a_i\in M$ is a primitive root modulo $p$, as computed by K. Matthews \cite{Matthews}. This formula
can be transformed using finite difference methods, into a form which allows one to prove that $D(i)>0$
if $i$ is not a power; see Buttsworth \cite{but1, but2}.
Paszkiewicz and Schinzel \cite{PaS} carried out extensive numerical work concerning $D(i)$ for
$i\le 32$. Paszkiewicz did extensive computations to study whether or not $g(p)=h(p)$. Litver and Yudina
had computed in 1971 that $g(p)=h(p)$ for $p\in [2,2^{32}]$ with the single exception of $p=40487$.
Paszkiewicz \cite{PaS2} showed that in the interval $[2^{32},10^{12}]$ there is one further such prime: $p=6692367337$.\\
\indent Let $d>1$ be an integer. Linnik proved in 1944 the celebrated result that there is a constant
$c$ such that every reduced residue class modulo $d$  contains a prime not exceeding $d^c$. Subsequent
proofs and improvements use, like Linnik, zero-free regions for $L$-series, zero density estimates for
these and a quantitative version of the Deuring-Heilbronn phenomenon (for an introduction, see Stopple \cite{Stopple}). 
This also applies to Heath-Brown's paper \cite{HB2} in which he
established $c=11/2$. Recently T. Xylouris in his Diplomarbeit (Bonn, 2009), established $c=5.2$, see \cite{Xylouris}.
Elliott \cite{pdta1} gave a rather different proof in which as a `waystation' (as
he called it), he established the result that the least prime primitive root, $g^*(p)$, to a prime modulus $p$ does not
exceed $p^c$. 
In another paper Elliott \cite{Ellio} proved that $g^*(p)< 475 (\log p)^{8/5}$ infinitely often.
Shoup \cite{shoup}, assuming GRH, gave a much sharper upper bound. G. Martin \cite{martin} showed
that for all but $O(Y^{\epsilon})$ primes $p\le Y$ we have $g^*(p)\ll \log^{C(\epsilon)}p$. Assuming GRH, 
Paszkiewicz and Schinzel \cite{PaS1} derived
from the work of K.R. Matthews \cite{Matthews} a formula for the density $E(q)$ of primes $p$ for which $g^*(p)=q$. They computed $E(q)$ for
the first 25 primes $q$. Bach \cite{Bachsearch} conjectured that 
$$\lim_{p\rightarrow \infty} \sup {g^*(p)\over (\log p)(\log \log p)^2}=e^{\gamma}.$$
Murata \cite{murat} showed that, under GRH, for any positive constant $D$ we have
$$|\{p\le x: g^*(p)\ge D(\log p) \log \log p\}|\ll_D {\pi(x)\over \log \log x}$$
and $g^*(p)=O(p^{\epsilon})$ for almost all primes $p$.\\
\indent A polynomial analogue of $g(p)$ also exists. There one can prove, as did Hsu \cite{Hsu2}, that there
exists a $c>0$ such that if $n\ge c\cdot{\rm deg}(P)/{\log_q {\rm deg}(P)}$, one can find at least one positive
irreducible of degree $n$ which is a primitive root modulo $P$. Davenport \cite{daven} had proven much earlier that
if $q$ (taken to be prime) is large enough with respect to deg$(P)$, then there even exists a linear 
polynomial which is a primitive root.\\
\indent Various authors \cite{hinz0, hinz1, hinz1b, hinz2, hinz3, WB, WG} considered the 
least primitive root in number fields. 
For example, Hinz \cite{hinz-1} generalized the P\'olya-Vinogradov inequality to arbitrary algebraic number
fields $K$. As an application he estimated the totally positive primitive root $\nu$ modulo a prime ideal 
$\mathfrak{p}$ of least norm and showed that it satisfies $N(\nu)\ll N(\mathfrak{p})^{1/2+\epsilon}$. In a later
paper \cite{hinz0} he improved this to $N(\nu)\ll N(\mathfrak{p})^{1/4+\epsilon}$.\\
\indent Dieulefait and Urroz \cite{dieu} used results on $g(p)$ in
order to study the malleability of RSA moduli. 
An encryption algorithm is malleable if it is possible for an adversary to transform 
a ciphertext into another ciphertext which decrypts to a related plaintext. That is, 
given an encryption of a plaintext $m$, it is possible to generate another ciphertext 
which decrypts to $f(m)$, for a known function $f$, without necessarily knowing or learning $m$.\\
\indent These results on $g(p)$ immediately give an efficient search procedure for primitive roots, that is, lead to 
the construction a small set $S$ with at least one element that is a primitive root modulo $p$. Bach \cite{Bachsearch} gave a different
search procedure, where the set $S$ has size $O(\log^4 p/(\log \log p)^3)$ and can also be constructed in polynomial time,
assuming GRH. However, the set tends to have larger elements. (In writing this section I made grateful use
of \S 2.2.4.3 in Narkiewicz's book \cite{Narkboek}.)\\

\noindent {\bf 28)} {\it Inverse primitive roots}. Zhang \cite{Zhang0} considered the distance between a primitive root $g$ with
$1\le g\le p-1$ and its inverse, ${\bar g}$ (which is also a primitive root), with $1\le {\bar g}\le p-1$. 
For about 75\% of the primitive roots $g$ modulo $p$ in $[1,p]$, one has that $|g-{\bar g}|\leq p/2$.\\

\noindent {\bf 29)} {\it Consecutive primitive roots}. Vegh, see, e.g., \cite{Vegh}, wrote a series of papers on consecutive primitive roots
and primitive roots in arithmetic progression. 
Mozingo \cite{moz}, using elementary arguments, showed that the only positive 
integers $>1$ having all their primitive roots consecutive are $2,3,4,5$ and $6$.\\

\noindent {\bf 30)} {\it Average multiplicative order}. J. von zur Gathen et al.\,\cite{veel} defined $u(n)$ to be the 
average multiplicative order
of the elements of $(\mathbb Z/n\mathbb Z)^*$. Note that $u(n)\le \lambda(n)$. They proved various
results comparing $u(n)$ with $\lambda(n)$, see also \cite{lucka, lucca, Mofield}.\\

\noindent {\bf 31)} {\it Smooth orders}. An integer is said to be ``$y$-smooth" (or ``$y$-friable") if none of 
its prime factors exceed $y$. The distribution of smooth integers is by now very well understood; for
a survey up to 1993 see Hildebrand and Tenenbaum \cite{HT}. 
A more recent survey is due to Granville \cite{granulaat}.
Banks et al.\,\cite{BFPS} established
some results on the smoothness of the order function, thus improving upon earlier results by
Pomerance and Shparlinski \cite{posh}.\\

\noindent {\bf 32)} {\it Distribution of $\varphi(p-1)/(p-1)$}. Trivially $\varphi(p-1)/(p-1)\le 1/2$ if $p$ is
odd. For any real numbers $x\ge 2$ and
$u$, let 
$$D(x,u)={1\over \pi(x)}\sum_{p\le x\atop \varphi(p-1)/(p-1)\le x}1.$$
Elliott \cite{pdta-1} proved that the limit ${\rm lim}_{x\rightarrow \infty}D(x,u)=D(u)$ exists for all
real numbers $u$. The function $D(u)$ turns out to be continuous and is strictly increasing on the interval
$[0,1/2]$. Schoenberg \cite{schoen}, see
also Kac \cite{Kac}, had earlier considered the distribution problem for $\varphi(n)/n$ and proved the
existence of the analogue of $D(u)$ and its continuity as a function of $u$. Li \cite{slis} considered
the analogue of $D(x,u)$ with $R(n)/\varphi(n)$ as function, with $R(n)$ the number of primitive $\lambda$-roots
in $[1,n]$. He showed for example that in this case there exists $u_0>0$ such that for every $u$ in $(0,u_0)$
the function $D(x,u)$ does not have a limit as $x\rightarrow \infty$.\\
\indent In the stochastic model proposed by Esseen \cite{Esseen} the sum
$$\sum_{p\le x}\Big({\varphi(p-1)\over p-1}\Big)^k$$
with $k\ge 1$ an integer played an important role. Without proof Esseen stated an asymptotic for this sum, namely
that the latter sum is asymptotic to $A_r'{\rm Li}(x)$ with $A_r'$ as in \S 9.7.2.
A proof was given much earlier by Vaughan \cite[Lemma 4.4]{va}. An easier proof
was given by Holden and Moree \cite{HM2} (following a suggestion of Pomerance). \\
\indent If in the latter sum we sum over integers $n\le x$ we have (see \cite[p. 42]{mova})
$$\sum_{n\le x}\Big({\varphi(n)\over n}\Big)^k=x\prod_q\Big(1-{1\over q}(1-(1-{1\over q})^k\Big)+O(x^{\epsilon}).$$
Vaughan \cite{va} studied the average relative density of
primitive roots mod $p$ versus non-primitive roots mod $p$,
$${1\over \pi(x)}\sum_{2<p\le x}{\varphi(p-1)\over p-1-\varphi(p-1)},$$
by relating it to the above `Esseen sums'.\\

\noindent {\bf 33)} {\it Cubic reciprocity primitive root criteria}. Chebyshev's criterium we mentioned in the introduction 
only involves quadratic reciprocity. Fueter \cite{Fueter} formulated some 
criteria which also make use of cubic reciprocity. He proved for
example that if $p$ is a prime of the form $p=1+2^{2n}3^{2m+3}$, with $n>0$ and $m\ge 0$, then 5 is a
primitive root modulo $p$ iff $m\equiv n ({\rm mod~}2)$.\\

\noindent {\bf 34)} {\it Coding theory}. 
In coding theory sometimes the multiplicative order arises.
See, e.g., \cite{BJ, DMS, gofu, guru, KL, Lenstracode, Moso, PSQ, RB, Rodier0, Rodier} and
cf. \S 9.4.2. For example Lenstra
applied his methods from \cite{Lenstra} to yield existence theorems, assuming GRH, for perfect, one-error-correcting
arithmetical codes, cf. \cite{gofu, Lenstracode}.\\

\noindent {\bf 35)} {\it Hooley-Heath-Brown Hybrids}. There are various results which are hybrid between Hooley's result and
Heath-Brown's unconditional result in the sense that under an assumption weaker than GRH a conclusion weaker
than Hooley's is established, but strongr than Heath-Brown's. For example Nongkynrih \cite{non2} proved that if for every squarefree 
integer $k$ the Dedekind
zeta function $\zeta_K(s)$ of the
field $K=\mathbb Q(\zeta_k,2^{1/k})$ is zero free in Re$(s)>1-e^{-12A/5}-\delta$, then 2 is a primitive root
for a positive proportion of primes. For other examples, see Cojocaru \cite{Co1, Co3} and Ram Murty \cite{Murty}.\\

\noindent {\bf 36)} {\it Polynomials representing primitive roots}. Let
$f$ be a polynomial over a finite field. Is there a primitive root in its image?
Madden \cite{mad} proved that for almost all finite fields $k$, every squarefree 
polynomial of degree $l$ over $k$ represents a primitive root in $k$. 
Let $f_1(x),\ldots,f_r(x)$ be polynomials modulo $p$, which are squarefree and
relatively prime in pairs. Carlitz \cite{car} gave an asymptotic formula for the number
of $x({\rm mod~}p)$ for which each of $f_1(x),\ldots,f_r(x)$ is a primitive root.
The best
results to date in this direction are due to S. Cohen \cite{cohenbest}.\\

\noindent {\bf 37)} {\it Permutations and primitive roots}. Consider the 
permutation 
$$(1,2)(1,2,3)\ldots(1,2,\ldots,m)$$
(multiplication of cycles from right to left), 
then (for a proof see \cite{Aul}) the product is a single cycle containing all of $1,2,\ldots,m$ iff $2m+1$ is 
a prime number having 2 as primitive root.\\
\indent To a primitive root $g$ modulo $p$ we associate the permutation $\sigma_g$ of $X$ defined
by $\sigma_g(x)\equiv g^x({\rm mod~}p)$. More precisely, $\sigma_g(x)=y$, the unique element in $\{1,2,\ldots,p-1\}$ 
satisfying $y\equiv g^x({\rm mod~}p)$. For example, if $p=7$, then $\sigma_5=(1~5~3~6)(2~4)$. 
Lewittes and Kolyvagin \cite{LKol} determined the parity $s(\sigma_g)$ of the permutation $\sigma_g$.
Assume $p>3$ and put $w=((p-1)/2)!$. Then
$$
s(\sigma_g)=
\begin{cases}
-wg^{(p-1)/4}({\rm mod~}2) & {\rm ~if~}p\equiv 1({\rm mod~}4);\cr
-w({\rm mod~}2) & {\rm ~if~}p\equiv 3({\rm mod~}4).
\end{cases}
$$
In case $p\equiv 1({\rm mod~}4)$, using Wilson's theorem one sees that $\pm w$ are precisely
the two roots of $X^2+1\equiv 0({\rm mod~}p)$. On the other hand, $\pm g^{(p-1)/4}$ are also roots
modulo $p$. The above results allows one to equate these roots.\\
\indent In case $p\equiv 3({\rm mod~}4)$, the  congruence can be rewritten as
$s(\sigma_g)\equiv h_{\mathbb Q(\sqrt{-p})}({\rm mod~}4)$, using Mordell's \cite{mordell} result
that $w\equiv (-1)^a ({\rm mod~}p)$, with
$a=(1+h_{\mathbb Q(\sqrt{-p})})/2$.\\

\noindent {\bf 38)} {\it Perfect card shuffling}. Consider 
an even number of cards numbered $0$ to $2k-1$, with 0 the top card. Then
a perfect card shuffle takes a card in position $i$ and sends it to $2i({\rm mod~}2k-1)$.
If the cards continue to be perfectly shuffled, the deck returns to the order it was in before
the first shuffle. For example, if one takes a deck of 52 cards, cut it in half, and
perfectly shuffles it (with the bottom card staying at the bottom), then it returns
after 8 times to its starting order. Thus Shuf$(52)=8$. In general let
Shuf$(2k)$ be the least number of perfect shuffles that will return a deck of $2k$ cards to its
starting order. It is not difficult so show that Shuf$(2k)={\rm ord}_{2k-1}(2)$, cf. \cite{card}.\\
\indent To wit, she that can card shuffle can do algebraic number theory, as the
level (see \S 9.4.3 for the definition) of the cyclotomic number field $\mathbb Q(\zeta_{2k-1})$ can be determined using the
parity of Shuf$(2k)$. This level is 4 if Shuf$(2k)$ is odd, and 2 otherwise.\\

\noindent {\bf 39)} {\it Parabolic generators of $\Gamma(p)$}. Frasch \cite{fra} has
shown that $\Gamma(p)$, the principal congruence subgroup of the modular group
$\Gamma$ modulo a prime $p$, can be generated as a free group by the generators
$S^p=({1\atop 0}{p\atop 1})$ and $(p-1)p(p+1)/12$ other generators. 
Recall that $g(p)$ is the smallest primitive root modulo $p$ (cf.~\S 9.7.27).
Grosswald
\cite{Grosswald-1} has shown that if $g(p)<\sqrt{p}-2$, then these $(p-1)p(p+1)/12$
generators
can be chosen in such a way that they are all parabolic. In a later paper, 
\cite{Grosswald}, he showed that for all $p$ sufficiently large we have
$g(p)<\sqrt{p}-2$.\\

\noindent {\bf 40)} {\it Pseudopowers}. An $x$-pseudopower to base $g$ is a positive integer which
is not a power of $g$, yet is so modulo $p$ for all primes $p\le x$. Let $q_g(x)$ denote the least
$x$-pseudopower to base $g$. Bach et al.\,\cite{BLSW} proved that if RH holds for Dedekind zeta functions, then
there is a constant $A>0$, depending only on $g$ such that $q_g(x)\ge \exp(A\sqrt{x}\log^{-2}x)$.
Estimating $q_g(x)$  is closely related to determining the behaviour of the order on average. 
Using upper bounds
of Baker and Harman \cite{baker} for the Brun-Titchmarsh inequality on average (see also
the book by Harman \cite{harman}) and new bounds on
exponential sums, Konyagin et al.\,\cite{KPS} showed that $q_g(x)\le e^{0.88715x}$ for 
$x$ sufficiently large and $|g|\le x$. This improved upon the upper bound
$q_g(x)\le e^{(1+o(1))x}$ due to Bach et al.\,\cite{BLSW}. In the course of the proof the authors
showed that $\prod_{p\le x,~p\nmid g}{\rm ord}_p(g)\ge e^{0.58045x}$ for $x$ sufficiently
large and for $g$ an integer with $2\le |g|\le x$. Recently the result by Konyagin et al. was improved
by Bourgain et al. \cite{BKPS} who showed that $q_g(x)\le e^{0.86092x}$ for $x$ sufficiently large and $|g|\le x$.\\ 

\noindent {\bf 41)} {\it Ducci sequences}. The Ducci-sequence generated by
$X:=(x_1,\ldots,x_n)\in \mathbb Z^n$ is the sequence $(X,D(X),D^2(X),\ldots)$ where
$D:\mathbb Z^n\rightarrow \mathbb Z^n$ is defined by
$D(x_1,\ldots,x_n)=(|x_1-x_2|,|x_2-x_3|,\ldots,|x_n-x_1|)$. Every Ducci sequence
$(X,D(X),D^2(X),\ldots)$ gives rise to a cycle; there are integers $i$ and $j$
with $0\le i<j$ with $D^i(X)=D^j(X)$. When $i$ and $j$ are as small as possible, we say
that the Ducci-sequence has period $(j-i)$. Breuer et al.\,\cite{Breuer} studied links
between Ducci-sequences, primitive polynomials and Artin's primitive root conjecture.
They proved for example that if $p$ is a prime and $2$ is a primitive root modulo $p$, 
then Ducci-sequences of length $p$ have only one period other than $1$.\\
\indent Breuer \cite{Breuer-1} introduced Ducci sequences over abelian groups. They are
sequences $U,D(u),D^2(u),\ldots \in G^n$, where $D(u_1,u_2,\ldots,u_n)=(u_1+u_2,u_2+u_3,\ldots,u_n+u_1)$.
In \cite{Breuer2} he proved that in case $G=(\mathbb Z/p^t\mathbb Z)$ and $p\nmid n$, the
maximal period of a Ducci sequences, $P_{p^t}(n)$, divides $p^{t-1}(p^{{\rm ord}_p(n)}-1)$ and
$p^{t-1}n(p^{{\rm ord}_p(n)/2}-1)$ if ord$_p(n)$ is even. Furthermore, he showed that
$P_{p^t}(2)={\rm ord}_{p^t}(2)$.\\

\noindent {\bf 42)} {\it Divisible trinomials}. Golomb and Lee \cite{gollee} considered irreducible polynomials
which divide trinomials in $\mathbb F_2$. In their considerations ord$_p(2)$ played
an important role. They proved for example that
if $\Phi_p(x)={x^p-1\over x-1}=f_1(x)f_2(x)\cdots f_r(x)$ is a product of
$r$ irreducible polynomials, $r$ is even and $p>7^{r/2}$, then the $f_i(x)$'s 
divide no trinomials. Note that $r=(p-1)/{\rm ord}_p(2)$. In this context Golomb made a conjecture
on near-primitive roots that subsequently turned out to be false; see \cite{MGolomb}.\\

\noindent {\bf 43)} {\it Sidelnikov primes}. In the theory of pseudo-random sequences, there is some
interest in finding prime powers $q$, such that $q-1$ has a large prime divisor $r$ such that 2 is a primitive
root modulo $r$. Let $P(x,y)$ be the number of primes $p\in (x,2x]$ such that $p-1$ has a divisor $r\ge y$
satisfying ord$_r(2)=r-1$. Friedlander and Shparlinski \cite{FS} showed that under GRH for any fixed
$1/2\le \alpha<17/32$, 
$$P(x,x^{\alpha})\ge \big(A\log(17/32\alpha)/100+o(1)){x\over \log x}.$$

\noindent {\bf 44)} {\it Diffie-Helllman triples}. Let $g$ be a primitive root modulo $p$. Canetti et al.~\cite{CFS} proved that the triples $(g^x,g^y,g^{xy})$, $x,y=1,\ldots,p-1$, are uniformly distributed modulo $p$ in the sense of H. Weyl. In cryptography it is often assumed that the triples $(g^x,g^y,g^{xy})$ cannot be
distinguished from totally random triples in feasible computation time. The above result implies that this
is in any case true for a constant fraction of the most significant bits, and for a constant fraction of the least significant bits.\\

\noindent {\bf 45)} {\it Carmichael numbers in number rings}. If $a^n\not\equiv a({\rm mod~}n)$, for some
integer $a$, then $n$ is composite. This is a rudimentary Compositeness Test. However, there are some
$n$ for which this test fails, no matter how we pick $a$. A {\it Carmichael number} is a composite integer
$n$ such that $a^n\equiv a({\rm mod~}n)$ for all integers $a$. The smallest such number is 561. It
turns out that there are infinitely many Carmichael numbers, see \cite{AGP}. Steele \cite{Steele} defined
a {\it Carmichael ideal} to be a composite ideal $\mathfrak n$ such for $\alpha \in \O_K$, we have
$\alpha^{N(\mathfrak n)}\equiv \alpha({\rm mod~}\mathfrak n)$. Using Heath-Brown's result that all
primes, with the possible exception of at most two are primitive roots for infinitely many primes $p$,
Steele \cite{Steele} showed that if the integer $n$ is the product of at least three distinct primes, then
there exist infinitely many cyclotomic extensions $\mathbb Q(\zeta_ q)$, such that $n$ is not
Carmichael in $\mathbb Q(\zeta_q)$.\\

\noindent {\bf 46)} {\it Common primitive root}. Let $\prod_{i=1}^r p_i^{e_i}$ be the canonical prime factorisation
of an integer $n$. Finizio and Lewis \cite{FL} called an integer $m$ a {\it common primitive root} if
$m$ is a primitive root modulo $p_i^{e_i}$ for $1\le i\le r$.\\

\noindent {\bf 47)} {\it Mapping the discrete logarithm}. A functional graph is a directed
graph such that each vertex must have exactly one edge directed out from it. If $S$ is
a finite set and $f:S\rightarrow S$ a mapping, then one can associate a functional graph
to the mapping $f$ by interpreting each element in $S$ as a vertex. The edges are defined
such that an edge $(a,b)$ is in the graph iff $f(a)=b$. Cloutier and Holden \cite{cloutier}
considered the functional graph associated to the mapping $x\rightarrow g^x$ modulo $p$, with
$p$ a prime and $S=\{1,2,\ldots,p-1\}$.  They computed various statistics of this
graph and compared them with those of a random mapping.\\

\noindent {\bf 48)} {\it Hilbert $p$-class field towers}. Furuta \cite{furuta} and
B. Schmidt \cite{schmidt} proved results relating the multiplicative order to sufficiency
criteria for the existence of Hilbert $p$-class field towers of $\mathbb Q(\zeta_m)$.
For example, Schmidt proved that if $m=kq$, where $k$ is an integer, and $q$ is a prime with
$q\equiv 1({\rm mod~}p)$, $(k,q)=1$ and $q^n\not\equiv -1({\rm mod~}k)$, for
$n=1,2,\ldots$, and also such that $\varphi(k)/{\rm ord}_k(q)\ge 8p+12$, then $\mathbb Q(\zeta_m)$ has
an infinite Hilbert $p$-class field tower. Shparlinski \cite{shparclass} used the latter result to show
that for infinitely many $m$ there is an infinite Hilbert $p$-class field tower over
$\mathbb Q(\zeta_m)$ for some $p\ge m^{0.3385+o(1)}$. Furthermore, he used the
result of Furuta \cite{furuta}  to show that for almost all positive integers $m$,
$\mathbb Q(\zeta_m)$ has an infinite Hilbert $p$-class field with high rank Galois
groups at each step, simultaneously for all primes $p$ to size up to about
$(\log \log m)^{1+\epsilon}$. As a consequence, Shparlinski inferred that for all $m\le x$, with
$x$ large enough,
except possibly $O(x(\log\log x)^{-0.08})$ exceptions, the class number of $\mathbb Q(\zeta_m)$
is divisible by all primes $p\le \log \log x/(10 \log \log \log x)$.\\

\noindent{\bf 49)} {\it Cunningham chains}. Let $p_1,\ldots,p_k$ be a chain of primes such that
for $2\le j\le k$ we have $p_j=2p_{j-1}+1$; that is we have
\begin{equation}
\label{ketting}
p_j=2^{j-1}p_1+2^{j-1}-1=2^{j-1}(p_1+1)-1.
\end{equation}
Such a chain of primes is known as a Cunningham Chain. A basic example is $2,5,11,23,47$, while the
longest one known is the Cunningham Chain with $k=16$ and $p_1=810433818265726529159$, discovered by
Carmody and Jobling in 2002. Let $k(p)$ be the length of the longest Cunningham Chain starting from $p$. One has $k(p)=1$ for
all primes $p\le x$, except for $O(x\log^{-3}x)$ of them. (This follows by a standard upper bound for
Sophie Germain primes coming from sieve theory.) By the first equality in (\ref{ketting})  
and Fermat's little theorem we find
$k(p)\le {\rm ord}_p(2)$. Conditonally one can do far better than this. Assuming that
${\mathcal{P}}(2)(x)\sim A\pi(x)$, which by Hooley's result is true under GRH, Ford et al. \cite{FKL} 
have shown that
$$\lim_{p\rightarrow \infty} \sup {k(p)\over \log p} \le {1\over A}.$$

\noindent{\bf 50)} {\it Dynamical systems and Artin}. Let $k$ be an ${\bold A}$-field. Denote by $P(k)$ the set of all places of $k$. Let $P_{\infty}(k)$ denote the set of Archimedean places if $k$ is a number field or the set of infinite places if $k$ is a function field of transcendence degree one over a finite field. 
Suppose one is given an element $\xi \in k^{\ast}$ξ, and any set $S \subset P(k) \subseteq P_{\infty}(k)$ such that $\xi$ is integral 
for all $w \notin S\cup P_{\infty}(k)$. The dual group of the $S$-integers $R_S$, denoted by $X$, is a compact abelian group. 
Let $\alpha \colon X \to X$ be a continuous group endomorphism which is the dual of the 
monomorphism $\widehat{\alpha}\colon R_S \to R_S$ given by $\widehat{\alpha}(x) = \xi x$. Dynamical systems of the 
form $(X, \alpha) = (X^{(k,S)}, \alpha^{(k,S,\xi)})$ are called {\it $S$-integer dynamical systems}. These were
introduced in Chothi et al.\,\cite{CEW}. The authors show, for example, that
if the qualitative version of Artin's conjecture holds, then there exist examples with uniformly 
distributed periodic points. 

In \cite{WardFF} (building on an earlier paper \cite{Ward})
Ward considers a family of isometric extensions of the full shift on $p$ symbols (with $p$ a prime) 
parametrized by a probability space. Using Heath-Brown's work \cite{HB} on the
Artin conjecture, he shows that for all but two primes $p$ the set of limit points of the
growth rate of periodic points is infinite almost surely. This shows in particular that
the dynamical zeta function is not algebraic almost surely.\\

\noindent {\bf 51)} {\it Converse theorems and Artin}. Associated to a newform $f(z)$ is a Dirichlet series $L_f(s)$ with functional equation and Euler product. Hecke
showed that if the Dirichlet series $F(s)$ has a functional equation of a particular form, then $F(s)=L_f(s)$ for some
holomorphic newform  $f(z)$ on $\Gamma(1)$. Weil extended this result to $\Gamma_0(N)$ under an assumption on the twists
of $F(s)$ by Dirichlet characters. Farmer and Wilson \cite{FW} take another approach by also making the assumption
that $L_f(s)$ has an Euler factor at the prime 2:
$$L(s,f)=(1-a_22^{-s}+2^{k-1-2s})\sum_{2\nmid n}{a_n\over n^s}.$$
Their main result has various conditions, one being a weak form of the Artin conjecture, described
in \S 10.36, that ensures that a
certain subset of matrices of $\Gamma_0(N)$ actually is a generating set for $\Gamma_0(N)$.\\

\noindent {\bf 52)} {\it Finding $M$ from a string of $g$-ary digits of $1/M$?} 
It has been shown by Blum et al. \cite{BBS} that given $k=2L+3$ successive digits of the
$g$-ary expansion of $1/M$, one can find $M$ in polynomial time $L^{O(1)}$. On
the other hand, under Artin's conjecture, it is also shown in \cite{BBS} that
$k=L-1$ digits are not enough to determine $M$ unambiguously. Konyagain and
Shparlinski \cite{KoSh} proved that for any $\epsilon>0$, given a string of
$k=[(3/37-\epsilon)\epsilon)L]$ consecutive $g$-ary integers, there are at
least $(1+o(1))\pi(g^L)$ prime numbers $p<g^L$ such that the $g$-ary expansion
of $1/p$ contains this string. Recently Bourgain et al. \cite{BKS} showed that
the ${3\over 37}=0.0810\ldots$ can be replaced by ${41\over 504}=0.0813\ldots$.\\

\noindent {\bf 53)} {\it Fermat quotients}. For a prime $p$ and and integer $u$ with $p\nmid u$, the Fermat quotient
is defined by the condition
$$q_p(u)\equiv {u^{p-1}-1\over p}({\rm mod~}p),~0\le q_p(u)\le p-1.$$
It has the property that $q_p(uv)\equiv q_p(u)+q_p(v)({\rm mod~}p)$.
It is expected that the map sending $u$ to $q_p(u)$ behaves very similarly to a random map on the 
set $\{0,\ldots,p-1 \}$. For results in this direction see \cite{Ostafe}. Let 
$l_p$ be the smallest positive integer $a$ such that $p\nmid q_p(a)$. It might be
true that $l_p\le 3$. Bourgain et al.~\cite{BFKS} have shown that
$l_p\le \log^{463/252+o(1)}p$ as $p$ tends to infinity. As a consequence, they derive a stronger
version of Lenstra's squarefree test. Shparlinski \cite{shparfermat} shows that for every $p$ there exists
$n\le p^{3/4+o(1)}$ such that $q_p(n)$ is a primitive root modulo $p$.\\

\noindent {\bf 54)} {\it Sums of squares of primitive roots modulo p}. 
Let $p\nmid c$ be an integer and let $M(p,c)$ be the number of solutions of the congruence 
$r^2+s^2\equiv c({\rm mod~}p)$ with $1\le r,s\le p-1$ primitive roots.  Zhang \cite{bao} shows that
$$M(p,c)={\phi^2(p-1)\over p}+O\Big({\phi^2(p-1)\over (p-1)^2}\sqrt{p}4^{\omega(p-1)}\Big).$$

\noindent {\bf 55)} {\it $G(n)$}. 
Let $G(n)$ be the least integer $G$ such that $\{y: 1\le y\le G~{\rm~and~}(y,n)=1\}$ generates 
$(\mathbb Z/n\mathbb Z)^*$. 
Any good upper bound for $G(n)$ would have implications
for deterministic primality testing. 
Assuming GRH, Montgomery \cite{monties} showed that $G(n)=O(\log^2 n)$, which
was made effective by Bach \cite{ebach} who showed that $G(n)\le 3\log^2 n$ for $n\ge 2$. As to unconditional
results: Bach and Huelsbergen \cite{BH} showed that $G(n)=O(\sqrt{n}(\log n)\log \log n )$, which
was sharpened by Burthe \cite{Bu1} to $G(n)\ll_{\epsilon} n^{\alpha+\epsilon}$, where $\alpha=1/(3\sqrt{e})=0.20217\ldots$.
The proof is based on the formula
$G(n)=\max\{\min\{a:\chi(a)\ne 0{\rm ~and~}\chi(a)\ne 1\}\}$,
where $\chi$ runs over the non-principal Dirichlet characters modulo $n$.
Norton \cite{norton} improved Burthe's bound to $G(n)\ll_{\epsilon} n^{\beta+\epsilon}$, where $\beta=1/(4\sqrt{e})=0.15163\ldots$. 
Bach and Huelsbergen \cite{BH}  conjectured that 
$$\limsup_{n\to\infty}{G(n)\over \log n\log\log n}=\frac1{\log 2}
{\rm ~and~}{1\over x}\sum_{n\le x}G(n)\sim (\log \log x )\log \log \log x.$$
Using zero density estimates Burthe \cite{Bu2} showed that $x^{-1}\sum_{n\le x}G(n)=O(\log^{97}x)$.\\
\indent Note that if $n$ is prime, $G(n)$ is bounded below by the least quadratic nonresidue of 
$n$, and above by the least primitive root of $n$. In this case $G(n)=\kappa(n)$, with $\kappa(n)$ defined
as in Section \S 9.7.2.\\ 

\noindent {\bf 56)} {\it Other Artin surveys}. For a much shorter introduction to Artin's primitive 
root conjecture 
see Ram Murty \cite{Murty1}. Li and Pomerance have written a survey with special emphasis on $\lambda$-roots, 
see \cite{LiP}. Cheng \cite{Cheng} wrote a survey on algorithms for finding 
primitive roots. Konyagin and Shparlinski \cite{KoSh} considered in their
book various questions related
to the distribution of integer powers $g^x$ for some integer $g>1$ modulo a prime $p$.

\section{Open problems}
\noindent {\bf 1)} Remove the GRH assumption that is made in many results in this area.\\

\noindent {\bf 2)} The two papers of Holden and Moree \cite{HM1, HM2} contain various open problems. 
Two of these have meanwhile been solved by Shparlinski; see \S 9.7.7.\\

\noindent {\bf 3)} Study the average behaviour of ranks in other families of elliptic curves than those
considered by Ulmer (see \S 9.6.13), for example those appearing in  Darmon \cite{darm}. \\

\noindent {\bf 4)} Let $g_1,\ldots,g_t$ denote the primitive roots modulo $p$.  The Dence brothers \cite{dence} considered
symmetric functions in the primitive roots of primes. They considered for example
$s_2(p)=\sum_{1\le i<j\le t}g_ig_j$ and show that this quantity, when considered modulo $p$, assumes
only values in $\{-1,0,1\}$. They wondered about the density 
$\delta(j)$ of primes for which, say $s_2(p)\equiv j({\rm mod~}p)$.
The author conjectures that $\delta(-1)=A/4$, $\delta(0)=1-A$ and $\delta(1)=3A/4$. 
In general one has the following problem: given a symmetric function in the primitive
roots, which values are assumed and with what frequency? For some partial progress see 
Moree and Hommersom \cite{Huib}.\\

\noindent {\bf 5)} Guy \cite[Section F9]{guy} formulated some unsolved problems regarding
primitive roots.\\

\noindent {\bf 6)} Find a lower bound, $b(p)$, say, such that
ord$_p(2)>b(p)$ for almost all $p$ (i.e. for all but $o(x/\!\log x)$ primes $p \le x$), as $x$ tends
to infinity. Erd\H{o}s \cite{E} proved that one can take $b(p)=p^{1/2-\delta}$. This was later improved
by Erd\H{o}s and Ram Murty \cite{EM} to $b(p)=p^{1/2+\epsilon(p)}$ for any function $\epsilon(p) \rightarrow 0$. They also proved that under the assumption of GRH for the Dedekind zeta function of $\mathbb{Q}(\zeta_d, 2^{1/d})$ one can take $b(p)=p/f(p)$ where $f(x)$ is any function that tends to infinity as $x\rightarrow \infty$.  For related results see Indlekofer and Timofeev \cite{indle} and Pappalardi \cite{papipo}. For a subgroup $\Gamma$ of rank $r\geq 1$ of $\mathbb{Q}^*$, Erd\H{o}s and Ram Murty \cite{EM} and Pappalardi \cite{papipo} considered the problem of finding  a lower bound $b(p)$, for almost all $p$,  such that $|\Gamma_p |>b(p)$, where $\Gamma_p$ is the reduction mod $p$ of $\Gamma$. Again they proved that under GRH one can take $b(p)=p/f(p)$ where $f(x)$ tends to infinity. See Elsholtz \cite{els} for a related result. Unconditionally Pappalardi proved that for almost all primes $p$ one can take
$b(p)=p^{r/r+1} \exp(\log^\tau{p})$ for some $\tau\simeq 0.15$.
The analogous problems for  subgroups of the group of rational points of CM elliptic curves were considered by Akbary \cite{great}, and Akbary and Kumar Murty \cite{akbary3}. For an elliptic curve $E$ defined over $\mathbb{Q}$ and a subgroup $\Gamma$ of rank $r$ of the group of rational points of $E$,
  Akbary, Ghioca, and Kumar Murty \cite{akbary2.5}  proved that, under GRH, one can can take $b(p)=p/f(p)$ as long as $r>18$ (for non-CM curves) or $r>  5$ (for CM curves). The problem remains open for lower ranks.
Duke \cite{Duke} considered the problem of finding a lower bound, for almost all $p$, for the exponent of reduction modulo $p$ of an elliptic curve $E$ and proved (under GRH) that for almost all primes $p$, the group of rational points of and $E$ mod $p$ contains a cyclic subgroup of order at least $p/f(p)$. A related dynamical interpretation for the size of reduction mod $p$ of the orbit of a point on a quasiprojective variety under the action of a monoid of endomorphisms of the variety proposed by Akbary and Ghioca \cite{AG}.
Related questions in a more general context of algebraic groups
and abelian varieties were considered by C. Matthews \cite{cmatthews}. He mentioned further applications
to nilpotent groups and to manifolds due to Milnor, Tits and Wolf.\\
\indent In the paper of Erd\H{o}s \cite{E} mentioned above, he conjectured that if $c<1$ one
can take $b(p)=p^c$. If one could take $c=0.8$ this would imply a conjecture of Ska\l ba \cite{skalba} to the extent
that almost all primes occur as prime divisors of the numbers $2^a+2^b+1$. Under GRH it follows by a result of
Pappalardi \cite[Theorem 2.3]{papipo} that we indeed can take $c=0.8$.
In the same paper Ska\l ba conjectured
that there are infinitely many primes $p$ such that $p$ does not divide any number of the form $2^a+2^b+1$. It
is easy to see that this set includes all Mersenne primes. Ska\l ba showed that if $\Omega(2^m-1)<\log m/\log 3$, then
there exists a prime divisor $q$ of $2^m-1$ such that $q$ is not a divisor of $2^a+2^b+1$ for any positive integers
$a$ and $b$. Luca and Stanica \cite{LucaSt} conjecture that $\omega((a^n-1)/(a-1))\ge (1+o(1))(\log n) (\log \log n )$ for almost all
integers $n$ and give heuristic arguments to support this. Their conjecture, if true, implies that the integers $m$ satisfying
Ska\l ba's condition $\Omega(2^m-1)<\log m/\log 3$ are not typical. Ska\l ba \cite{skalba} conjectures that there are
infinitely many primes $q$ such that $q$ does not divide any number of the form $2^a+2^b+1$.\\

\noindent {\bf 7)} Roskam \cite{Ros1a} raised the following questions and demonstrated their applications
to the study of divisors of general linear recurrences: let $f\in \mathbb Z[X]$ be an irreducible monic
polynomial, $\alpha$ a root of $f$, and define $K=\mathbb Q(\alpha)$ with ring of integers $\O_K$. Furthermore,
let $T$ be the set of primes that are inert in $K/\mathbb Q$, and define for each $h\ge 1$ the set $T_h$ of primes
$p$ in $T$ for which the subgroup $\langle {\overline \alpha}\rangle\subset ({\O}_K/p{\O}_K)^*$ has
index $h$. Is it true for generic $\alpha$ that $T_h$ has a natural density $\delta(T_h)$?
Is it true that $\sum_{h=1}^{\infty}\delta(T_h)=\delta(T)$? Earlier Brown and Zassenhaus had made
a similar conjecture \cite{BZ}.\\

\noindent {\bf 8)} Let $q$ be an odd given prime power and $m$ a natural number. One can wonder whether
there exists an extension ${\mathbb F}_{q^n}$ of ${\mathbb F}_{q^m}$ such that
${\mathbb F}_{q^n}$ has an optimal normal basis over $\mathbb F_q$. A number
theoretic question that arises in this context, see
\cite{CGPT},  is whether there is a prime $p$
such that $p\equiv 1({\rm mod~}m)$, and $q$ is a primitive root modulo $p$, and if
yes to provide a small upper bound for the smallest such prime. The latter part
of the question seems very difficult (it is already a very hard
problem to find a bound for the smallest prime $p$ such that
$p\equiv 1({\rm mod~}m)$, see, e.g., Heath-Brown \cite{HB2}), but it should be possible
to shed some light on the first part of the question.\\

\noindent {\bf 9)} For a given odd prime $l$ and prime power $q$, Ballot 
\cite{Ballot1}
has computed the Dirichlet
density of primes $P$ of $\mathbb F_q[X]$ such that $l$ divides the order of $X({\rm mod~}P)$,
by an elementary method using neither Kummer theory, nor the Chebotarev Density Theorem. 
In an earlier paper Ballot \cite{Ballot2} considered this problem for the case $l=2$.
Can
one likewise compute the Dirichlet density of primes $P$ such that the order of $X({\rm mod~}P)$ is
in a prescribed arithmetic progression? Indeed, in this setting one can wonder whether the Dirichlet
density is the `good density' to consider. For a lively discussion of this problem see
Ballot \cite{Ballot3}.\\

\noindent {\bf 10)} Let $K$ be a number field, $F/K$ a Galois extension and $C$ a union of
conjugacy classes of Gal$(F/K)$. Furthermore, let $\alpha\ne 0$ be an algebraic
integer in a number field $K$ which is not a root of unity.
In \cite{Ziegler}, Ziegler was interested in the
set $P_{\alpha}(K,F,C,a,d)$ of primes ${\mathfrak p}$  of $K$ such that
$({\mathfrak p}, F/K)\in C$, ${\mathfrak p}\nmid (\alpha)$, and
ord$_{\mathfrak p}(\alpha)\equiv a({\rm mod~}d)$, where 
$({\mathfrak p}, F/K)$ denotes the Frobenius automorphism of ${\mathfrak p}$, 
and ord$_{\mathfrak p}(\alpha)$ denotes the order of the algebraic
integer  $\alpha $ in ${\O}/{\mathfrak p}{\O}$, ${\O}$ being
the ring of integers of $K$. Ziegler \cite{Ziegler} showed assuming GRH
that
this set has a natural density that can be given as a double sum
involving field degrees and Galois intersection coefficients. In the
case $K=\mathbb Q$ his result simplifies to results proven by Moree  \cite{Moreejapi2}.
In the sequel to the latter paper, \cite{Moreejapi3}, Moree showed that there is a `generic' density,
and that almost always the density equals the generic density. Do similar
results hold for the more general case considered by Ziegler as well?\\

\noindent {\bf 11)} Moree \cite{Model} gave an asymptotically exact heuristic for the number of
primes $p\le x$ dividing a sequence of the form $\{a^k+b^k\}_{k=1}^{\infty}$. These numbers
satisfy a linear recurrence of degree 2. Is it possible to derive an asymptotically exact
heuristic for the number
of primes $p\le x$ dividing  a linear recurrence of degree 2 consisting of integers only?\\

\noindent {\bf 12)} Let $g$ be a primitive root modulo $p$. 
How large is the
smallest $N$ such that any residue class modulo $p$ is representable in the 
form $g^x-g^y\pmod p$ with $1\le x,y<N$? This question is attributed to Odlyzko, who 
conjectured that one can take $N$ as small as $p^{1/2+\epsilon}$, for any fixed $\epsilon$
and $p$ large enough in terms of $\epsilon$.  Rudnick and Zaharescu \cite{RZ} proved that
one can take $N=c_op^{3/4}\log p$, and this was sharpened by 
Cuauht\'emoc Garcia \cite{mexico} to $N=2^{5/4}p^{3/4}$. 
Using properties of Sidon sets the latter estimate was improved to $N=(\sqrt{2}+\epsilon)p^{3/4}$ by
Cilleruelo \cite{Javier} (with $p>p(\epsilon)$) and more recently by Cilleruelo and 
Zumalac\'arregui \cite{Javier2} to $N=\sqrt{2}p^{3/4}$ for every prime $p$.
For
related papers on the distribution of small powers of a primitive root see \cite{CGZ, mont, VZ}.\\

\noindent{\bf 13)} Grosswald \cite{Grosswald} wrote that an inspection of tabulated 
values of $g(p)$, the smallest positive primitive root modulo $p$, strongly suggests that
$g(p)<\sqrt{p}-2$ for all primes $p>409$. He proved that $g(p)<p^{0.499}$ if
$p-1\ge e^{e^{24}}$. Can one lower this lowerbound? (Andrew Booker, personal
communication, suggested the answer is yes and even substantially so.) For an
application of this bound see \S 9.7.39.\\

\noindent{\bf 14)} Cat maps are maps of the unit torus given by $2\times2$ matrices with integer entries, 
unit determinant, and real eigenvalues. The chaotic nature of these maps is traditionally depicted by
showing what the map does to the face of a cat (after several iterations).
The cat map models chaotic components of phase space in Hamiltonian systems. In the analysis of
its periodic orbits the order of the cat map modulo $p$ plays an important role. Keating \cite{Keating}
studied these aspects using heuristic arguments. Meanwhile a much more rigorous treatment using
for example the methods in Kurlberg's paper \cite{K1} is possible. As a working title for
this project the author proposes: Cat maps: the dogged approach. (See also \S 9.7.10 and \cite{Dyson, perci}.)\\ 

\noindent{\bf 15)} In connection with a question related to codes, Rodier \cite{Rodier0} was interested
in computing the density of primes $p\equiv 7({\rm mod~}8)$ such that the subgroup in $(\mathbb Z/p\mathbb Z)^*$ generated
by 2 is of index 2.  
He argued that assuming GRH the density is $A/2$. 
Recently, Peter Malicky independently in connection with work on the periodic orbits of a 
certain 2-dimensional dynamical system asked about the infinitude of this set of primes. 
More generally one might ask for an
explicit evaluation of the density of primes $p\equiv a({\rm mod~}f)$ such that the subgroup 
in $(\mathbb Z/p\mathbb Z)^*$ generated
by $g$ is of index $t$.\\

\noindent{\bf 16)}  Brauer \cite{brauer} showed that the infinitely
many primes $p\equiv 9({\rm mod~}16)$ which can be represented as
$65x^2+256xy+256y^2$ are all non-divisors of the sequence $\{2^k+1\}_{k=1}^{\infty}$. 
Can this be generalized?\\

\noindent {\bf 17)} {\it Transposition invariant words}. In the context of a problem
of words invariant under certain transpositions, Lepist\H{o} et al. \cite{LPS} prove
that there are infinitely many primes $p$ such that there is a primitive root mod $p$
that divides $p+1$. They conjectured that this set of primes has a density
$B\ge 0.65$ that should be exactly computable assuming GRH.\\

\noindent{\bf 18)} Nagata \cite{Nagata} raised a question on $\mathbb Z[\sqrt{14}]$. A positive answer would give a proof different from Harper's \cite{Harper} that $\mathbb Z[\sqrt{14}]$ is Euclidean.\\

\noindent{\bf 19)} Granville and Soundararajan \cite{GSound} studied a conjecture of Erd\H{o}s,
that every odd positive integer can be written as the sum of a squarefree number and
a power of 2. They showed that this is connected with the behaviour of ord$_{p^2}(2)$ and ord$_p(2)$ and made various conjectures
on these orders. They showed for example that if $\sum_{p^2|2^{p-1}-1}1/{\rm ord}_p(2)<\infty$, then there
exists an integer $k$ such that almost every odd integer can be written as the sum of a squarefree number plus
no more than $k$ distinct powers of 2. They conjectured that there exists a constant $\delta>0$ such
that, for any finite set of primes $\{p_1,p_2,\ldots,p_m\}$ and any choice of integers $a_1,a_2,\ldots,a_m$, the
proportion of positive integers which belongs to at least one of the congruence classes
$a_i({\rm mod~}{\rm ord}_{p_i^2}(2))$, is $<1-\delta$. They showed that $\sum_p 1/{\rm ord}_{p^2}(2)<\infty$ if
this conjecture is true. In case the latter sum is $<1$, they inferred that all but $O(x/\log x)$ of the
odd integers $n\le x$ can be written as a sum of a squarefree number and a power of 2. They also asked what
the true order of magnitude is of
$$\prod_{p\le x}\Big(1-{1\over {\rm ord}_p(2)}\Big).$$

\noindent{\bf 20)} Elliott and Murata \cite{pdta2} conjectured that the following limits all exist and are
finite:
$$\lim_{x\rightarrow \infty}{1\over \pi(x)}\sum_{p\le x}g(p)
{\rm ~and~}\lim_{x\rightarrow \infty}{1\over \pi(x)}\sum_{p\le x}g(2p^2),$$
and likewise with $g$ replaced by $g^*$ (recall that $g^*(p)$ is the smallest prime
primitive root modulo $p$). Bach \cite{Bachsearch} conjectured that 
$$\lim_{p\rightarrow \infty} \sup {g^*(p)\over (\log p)(\log \log p)^2}=e^{\gamma}.$$

\noindent {\bf 21)} {\it Schinzel-W\'ojcik problem}. Given any rational $a,b\in \mathbb Q\backslash \{0,\pm 1\}$, Schinzel and W\'ojcik
\cite{SW} proved that there exist infinitely many primes $p$ such that $\nu_p(a)=\nu_p(b)=0$
and ord$_p(a)={\rm ord}_p(b)$. The first condition ensures that the orders in the second one are
defined and is satisfied for all but finitely many primes. More generally, given
$a_1,\ldots,a_s\in \mathbb Q\backslash \{0,\pm 1\}$, the Schinzel-W\'ojcik problem is to determine
whether there exist infinitely many primes $p$ for which the order modulo $p$ of each
$a_1,\ldots,a_s$ is the same. Pappalardi and Susa \cite{PSu} proved assuming GRH that the primes with this
property have a density, and in the special case where each $a_i$ is a power of a fixed rational number,
they showed unconditionally that such a density is positive. In the case where all the $a_i$'s are prime, they
expressed the density as an infinite product.\\

\noindent {\bf 22)} Y. Gallot found a quadratic polynomial $f$ and an integer $g$ such that
the first (distinct) 38639 primes $p$ produced by $f$ are such that $g$ is a primitive root modulo $p$. Can this
be improved (cf. \S 9.7.11)?\\

\noindent {\bf 23)} K. Szymiczek asked the following question: Do there exist infinitely many natural numbers
$n$ such that $(2^n-3,3^n-2)=1$? In this direction Ska\l ba \cite{skalba2}, using ideas from a paper of
Erd\H{o}s \cite{E}, proved that the number of primes $p<x$ dividing both $2^n-3$ and $3^n-2$ is
$O(x\log^{-1.0243}x)$. Ballot and Luca \cite{Balluc} considered the number
of primes $p\le x$ such that $p$ divides both $a^n-b$ and $c^n-d$.\\

\noindent {\bf 24)} Does the unconditional bound as established by
Vinogradov, see (\ref{wijn}), hold in greater generality?\\

\noindent {\bf 25)} Let $q$ be an odd prime power. For each integer $a$ with $(a,q)=1$, is there a pair of primitive roots $g_1$ and 
$g_2$ modulo $q$ such that $a^2\equiv g_1g_2({\rm mod~}q)$?  Zhang \cite{Zhang1} proved that if $q$ is
an odd prime and $({a\over q})=1$, then the answer is yes.\\

\noindent {\bf 26)} Prove or disprove the conjecture of Bach et al. \cite{BLSW} to
the effect that $q_g(x)$, the least $x$-pseudopower to base $g$, should be
about $\exp(c_gx/\log x)$, where $c_g>0$ (see \S 9.7.40).\\

\noindent {\bf 27)} Let $G(n)$ be the smallest integer $k$ such that the 
primes $p\leq k$ generate the multiplicative group modulo $n$. Prove or disprove
the conjectures formulated by Bach and Huelsbergen \cite{BH}, see \S 9.7.55.\\ 

\noindent {\bf 28)} {\it V.I. Arnold's conjectures}. In a series of
papers, Arnold \cite{arnold1, arnold2, arnold3, arnold4, arnold5, arnold6, arnold7}
considered dynamical systems related to linear transformations in finite fields
and residue rings and made a series of conjectures. Of these many were
confirmed, several refuted, and Shparlinski \cite{shpar-1} observed that several required adjustments. 
For example, in
\cite{arnold5} Arnold considered the quantity
$$T_g(L)={1\over L}\sum_{l=1,~(g,l)=1}^L {\rm ord}_l(g),$$
where the sum is over the integers $l$ coprime with $g$ and suggested it grows
asymptotically like $c(g)L/\log L$, with $c(g)>0$ a constant. Under GRH 
Shparlinski showed, however, that
$$T_g(L)\ge {L\over \log L}\exp(C(g)(\log \log \log L)^{3/2}),$$
for some constant $C(g)>0$. Shparlinski observed that it should be possible
to show under GRH that 
$$T_g(L)\ge {L\over \log L}\exp((\log \log \log L)^{2+o(1)}).$$
Some of Shparlinski's other results in \cite{shpar-1} have recently been improved
by Chang \cite{chang}.\\
\indent Since ${\rm ord}_l(g)\le \lambda(l)$, with $\lambda$ the Carmichael function, and on GRH there is a set of
positive upper density of numbers coprime to $g$ such that 
equality holds as Li and Pomerance \cite{LiP2} showed, it is natural
to compare $T_g(L)$ with $L^{-1}\sum_{l=1}^L \lambda(l)$.
Erd\"os et al.\,\cite{EPS} proved that
\begin{equation}
\label{lollie}
{1\over x}\sum_{n\le x}\lambda(n)={x\over \log x}\exp\Big({B\log \log x\over \log \log \log x}(1+o(1)\Big),
\end{equation}
where 
$$B=e^{-\gamma}\prod_q\big(1-{1\over (q-1)^2(q+1)}\Big)=0.3453720641\ldots$$
\indent Kurlberg and Pomerance \cite{KPStephens} showed that under GRH actually 
$$T_g(x)={x\over \log x}\exp\Big({B\log \log x\over \log \log \log x}(1+o(1)\Big),$$
as $x$ tends to infinity, uniformly in $g$ with $1<|g|\le \log x$. On comparing this result with (\ref{lollie}) we see that the mean values of $\lambda(n)$
and ord$_n(g)$ are of a similar order of magnitude. 
Further, in the same paper Kurlberg and Pomerance established an asymptotic for
$${1\over \pi(x)}\sum_{p\le x}{\rm ord}_p(g).$$
Earlier Stephens \cite{Stephens3} showed that on GRH, the limit
$$\lim_{x\rightarrow \infty}{1\over \pi(x)}\sum_{p\le x}{{\rm ord}_p(g)\over p-1},$$
exists and equals the Stephens constant $S$ time a rational correction factor $c_g$ depending on $g$ 
(caution: his $c_g$ must be adjusted in certain cases as noted by Moree and Stevenhagen \cite{MS}).\\
 
\noindent {\bf 29)} (Erd\H{o}s, \cite{romer-1}). Is it true that $7,15,21,45,75$ and $105$ are the
only integers $n$ for which $n-2^k$ is a prime for every $k$ with $2\le 2^k<n$? (See \S 9.7.20)\\

\noindent {\bf 30)} Are there Crandall numbers that are not Wieferich numbers other than 5 and 181?
Are there infinitely many Crandall numbers that are not Wieferich numbers? (See \S 9.7.16 for definitions and 
references.)\\

\noindent {\bf 31)} Show that
$$\sum_{n\le x,~(n,a)=1}{1\over {\rm ord}_n(a)}=O(x^{1/4}).$$
Then Artin's primitive root conjecture follows by Ram Murty and Srinivasan \cite{MSri}, who
conjectured that actually $O(x^{\epsilon})$ should hold true. Felix \cite{Felix} showed that the above
sum is $\log x$ on average.\\

\noindent {\bf 32)} {\it Very odd sequences}. 
Let $V(x)$ denote the number of binary integers $n=pq\le 2x-1$ with $p<q$ such
that 
ord$_p(2)=(p-1)/2$, ord$_q(2)=(q-1)/2$ and $({\rm ord}_p(2),{\rm ord}_q(2))=1$.
The author is inclined to believe that $V(x)$ counts a positive fraction of all binary
integers $n\le 2x-1$, that is
$$V(x)\sim c_0{x\log \log x\over \log x},$$
for some positive constant $c_0$. If true, this has some consequence for the
theory of {\it very odd sequences} (introduced by J. Pelik\'an in 1973 \cite{pelikaan}). For a given
number $n$ fix integers $a_i$ with $a_i\in \{0,1\}$ for $1\le i\le n$. Put
$A_k=\sum_{i=1}^{n-k}a_ia_{i+k}$ for $0\le k\le n-1$. We say that
$a_1\ldots a_n$ is a {\it very odd sequence} if $A_k$ is odd for $1\le k\le n$. By $S(n)$
we denote the number of very odd sequences of length $n$ and by $N_k(x)$ the number
of integers $m\le x$ such that $S(m)=k$. Then if the above conjecture holds
true we have, as $x$ tends to infinity, 
$$N_{16}(x)\sim V(x)\sim c_0{x\log \log x\over \log x}.$$
See Moree and Sol\'e \cite{MS} for further information.\\

\noindent {\bf 33)} {\it Cherednik's heuristic}. The heuristic considerations sketched in this survey
begin with the behaviour of the sum $\sum_{p\le x}\varphi(p-1)/(p-1)$. Cherednik \cite{Cherednik}
proposed to work instead with 
$${\sum_{p\le x}\varphi(p-1)\over \sum_{p\le x}p-1}.$$
From Pillai's work \cite{pi0} it follows that the latter ratio has $A$ as a limit. One can try to redo the asymptotically
exact heuristics of the author in this spirit. Some numerical work, see 
\cite{Cherednik}, suggests that perhaps the Cherednik variation
of the heuristics tends to be numerically closer to the actual primitive root counting function than the standard heuristics. Try 
to investigate and understand this.\\

\noindent {\bf 34)} {\it Least primitive root (mod $p$) versus least primitive root (mod $p^2$)}. 
Paszkiewicz \cite{PaS2} conjectures that for most primes $p$ we have $g(p)=h(p)$. Quantify this: do there exist 
infinitely many primes $p$ for which $g(p)\ne h(p)$? (Presently only two are known, see also \S 9.7.27).\\

\noindent {\bf 35)} {\it Primitive roots in image of polynomial mapping}.  Let $f$ be a polynomial over a finite field. Is there a primitive root in its image? (See \S 9.7.36.) \\

\noindent {\bf 36)} {\it Farmer-Wilson variant of the Artin conjecture}. This conjecture
states that if $(d,bM)=1$, then there exist integers $k$ and $n$ such that 
$b\equiv 2^n({\rm mod~}d+kM)$.\\
\indent Note that this follows from Artin's conjecture, for if $p=d+kM$ is prime and $2$ a primitive root
modulo $p$ we are done. Since we do not require $d+kM$ to be prime, nor 2 to be a primitive root modulo $d+kM$,
the above conjecture is actually much weaker than Artin's. (See also \S 9.7.51.)\\  

\noindent {\bf 37)} Let $M\in \Bbb F_q[T]$ be a fixed polynomial and $k\ge 2$ be an integer. Determine the density of the set of all monic irreducible polynomials $P$ for which $P+M$ is a $k$-free polynomial. 
This is the polynomial variant of problems discussed
in \S 9.7.1. and was earlier considered by Wei-Chen Yao \cite{wcyao}, but
according to MR1841909 (2002d:11144) his paper contains mistakes.\\

\noindent {\bf 38)} Let $p\equiv 1({\rm mod~}n)$. Put $s_n(p)=\sum \{g/p\}$, where the sum is over the subgroup
of $(\mathbb Z/p\mathbb Z)^*$ of order $n$, and $\{~\}$ denotes the fractional part. Hadano et al.\,\cite{hada} raise 7 conjectures related to this quantity.\\ 

\noindent {\bf 39)} {\it Divisors of $a^{f(n)}-1$}. {}From Fermat's little theorem, we know that the set of primes which divide $a^n-1$ for some $n$ is precisely
the set of primes not dividing $a$. Ballot and Luca \cite{ballotluca} investigated what happens if we replace the exponent $n$ here
by a different polynomial expression in $n$. Let $a$ be an integer with $|a|>1$ and $f(X)\in \mathbb Q[X]$ a nonconstant, integer-valued
polynomial with positive leading term. Suppose that there are infinitely many primes $q$ for which $f$ does not possess a root
modulo $q$, then Ballot and Luca showed that almost all primes $p$ do not divide any number of the form $a^{f(n)}-1$.
Their result was sharpened by Pollack \cite{pollack} who also gave a conjectural asymptotic for the number of primes dividing
the sequence $a^{f(n)}-1$, under the further assumption that the set of primes $q$ for which $f$ does not possess a root
modulo $q$, has positive Dirichlet density.\\

\noindent {\bf 40)} {\it Generalized Titchmarsh divisor problem}. Let $\tau(n)$ denote the number of divisors of the positive integer $n$. In 1931,
Titchmarsh \cite{titchmarsh} established the following estimate, assuming GRH,
$$\sum_{a<p\le x}\tau(p-a)=x\prod_{q|a}\big(1-{1\over q}\Big)\prod_{q\nmid a}\big(1+{1\over q(q-1)}\Big)+
O\Big({x\log \log x\over \log x}\Big).$$
In 1961, Linnik \cite{linnik} established the above asymptotic unconditionally by using his dispersion
method. Later Rodriquez \cite{rod} and independently Halberstam \cite{halber}, by a straightforward application
of the Bombieri-Vinogradov theorem, proved unconditionally the Titchmarsh conjectural
asymptotic formula. In the special case $a = 1$ it is a trivial observation that
$$\sum_{p\le x}\tau(p-1)=\sum_{1\le m\le x-1}\pi(x;m,1)=\sum_{m\ge 1}\pi(x;m,1).$$
Essentially we are counting each prime number $p\le x$ for each occurrence of $m$ such that $p$ splits
completely in $\mathbb Q(\zeta_m)$.\\
\indent Let $\mathfrak F=\{F_m:m\ge 1\}$ be a family of finite Galois extensions of $\mathbb Q$. For each
$m$, let $D_m$ be a union of conjugacy classes of Gal$(F_m/{\mathbb Q})$ and let $\tau_{\mathfrak F}(p)$ be the number
of $m\ge 1$ such that $p$ is unramified in $F_m/\mathbb Q$ and the Artin symbol $\sigma_p$ belongs to $D_m$. Suppose
that $\tau_{\mathfrak F}(p)<\infty$ for each prime $p$. Then we have the following problem:\\
{\it Generalized Titchmarsh divisor problem}:\\
Determine the behaviour of $\sum_{p\le x}\tau_{\mathfrak F}(p)$ as $x$ tends
to infinity.\\
\indent In the above generality, the problem is too unwieldy unless some constraint on the sizes of $D_m$ are imposed. Let $A$ be an abelian variety defined over $\mathbb Q$ and for each $m\ge 1$, let $A[m]$ be the set of torsion points
of $A$ of order dividing $m$. Let $F_m=\mathbb Q(A[m])$ and $D_m=\{Id\}$. Therefore $\tau_{\mathfrak F}(p)$ is the number of $m\geq 1$ such that $p$ splits completely in $\mathbb Q(A[m])$. In this setup the problem specializes to the Titchmarsh divisor problem
for Abelian varieties and it is analogous to the classical case for $\tau(p-1)$. Akbary and Ghioca \cite{akbary2} believe that in this case
\begin{equation}
\label{akkie}
\sum_{p\le x}\tau_{\mathfrak F}(p)={\rm Li}(x)\sum_{n=1}^{\infty}{1\over [\mathbb Q(A[m]):\mathbb Q]}+o\Big({x\over \log x}\Big)
\end{equation}
and under GRH prove a special case. In case $A=E$ is a CM elliptic curve, they unconditionally prove (\ref{akkie}).\\
\indent We note that Felix \cite{Felix} gives an estimate for
$\sum_{p\le x,~p\equiv a({\rm mod~}k)}\tau((p-a)/k)$ and applies it to a question related to
Artin's primitive root conjecture (see \S 9.7.3).\\

\noindent {\tt Remark}. If the reader has not suffered an overdose of primitive
roots after reading this survey, (s)he is recommended to consume them at
industrial strength \cite{DD}.

\section{Acknowledgement} The basis for this paper was a survey lecture I gave at Oberwolfach
concerning Artin's primitive root problem
I am grateful for the privilege
of having visited Oberwolfach several times and would like to than Z. Rudnick for inviting me to give the
lecture on which this survey is based. 
MathSciNet turned out
to be an extremely helpful tool in writing this
survey. 
For close to 20 years now I regularly communicated with Peter Stevenhagen concerning Artin type problems. I
am very grateful for his tremendous help.
Hendrik Lenstra pointed out various inaccuracies in an earlier version.
S.D. Cohen kindly updated me on the status of Golomb's conjectures.
W. Narkiewicz and Ana Zumalac\'arregui provided me with some helpful information concerning
problems 29, respectively 12. 
Kate Petersen fine-tuned the description I gave of her work. Further thanks are due to Joachim von zur Gathen
and the referee.
Many (language) comments were provided by Julie Rowlett; a native English speaker. I am very
grateful for her willingness to invest so much time in this. Comments of A. Akbary, A. Felix and P. Kurlberg on
an earlier version are taken into account. Akbary kindly rewrote my initial version of Problems 6 and 40 a bit.

An earlier and rather shorter version of this paper (with the
same title) is available from the arXiv \cite{survie} (its aim was to
describe the status quo up to 2004). Many thanks are due to Alina
Cojocaru for her contributions on elliptic Artin to the survey and
for her comments and suggested corrections to the rest of survey.
Wojciech Gajda kindly wrote a section on Artin and K-theory. Hester Graves
carefully studied the literature on Artin and Euclidean domains and kindly reported
on that, using my initial version of this as a starting point.
Last but not least I
thank Francesco Pappalardi and Carl Pomerance for urging me to publish this survey. Special thanks are due to Carl for his support
and many helpful comments.

\begin{center}\sc Pieter Moree\\
Max-Planck-Institut f\"ur Mathematik,\\ 
Vivatsgasse 7\\
D-53111 Bonn, Deutschland\\ 
{\em e-mail address:} \rm moree@mpim-bonn.mpg.de\\
\end{center}
\end{document}